
\def\Dist{\roman{Dist}}

\overfullrule=0pt

\documentstyle{amsppt}
\newcount\mgnf\newcount\tipi\newcount\tipoformule\newcount\greco
\tipi=2          
\tipoformule=0   

\global\newcount\numsec\global\newcount\numfor
\global\newcount\numapp\global\newcount\numcap
\global\newcount\numfig\global\newcount\numpag
\global\newcount\numnf
\global\newcount\numtheo

\def\SIA #1,#2,#3 {\senondefinito{#1#2}%
\expandafter\xdef\csname #1#2\endcsname{#3}\else
\write16{???? ma #1,#2 e' gia' stato definito !!!!} \fi}

\def \FU(#1)#2{\SIA fu,#1,#2 }

\def\etichetta(#1){(\veroparagrafo.\veraformula)%
\SIA e,#1,(\veroparagrafo.\veraformula) %
\global\advance\numfor by 1%
\write15{\string\FU (#1){\equ(#1)}}%
\write16{ EQ #1 ==> \equ(#1)  }}

\def\etichettat(#1){\veroparagrafo.\veratheorema:%
\SIA e,#1,{\veroparagrafo.\veratheorema} %
\global\advance\numtheo by 1%
\write15{\string\FU (#1){\thu(#1)}}%
\write16{ TH #1 ==> \thu(#1)  }}

\def\etichettaa(#1){(A\veraappendice.\veraformula)
 \SIA e,#1,(A\veraappendice.\veraformula)
 \global\advance\numfor by 1
 \write15{\string\FU (#1){\equ(#1)}}
 \write16{ EQ #1 ==> \equ(#1) }}
\def\getichetta(#1){Fig. \verafigura
 \SIA g,#1,{\verafigura}
 \global\advance\numfig by 1
 \write15{\string\FU (#1){\graf(#1)}}
 \write16{ Fig. #1 ==> \graf(#1) }}
\def\retichetta(#1){\numpag=\pgn\SIA r,#1,{\verapagina}
 \write15{\string\FU (#1){\rif(#1)}}
 \write16{\rif(#1) ha simbolo  #1  }}
\def\etichettan(#1){(n\verocapitolo.\veranformula)
 \SIA e,#1,(n\verocapitolo.\veranformula)
 \global\advance\numnf by 1
\write16{\equ(#1) <= #1  }}

\newdimen\gwidth
\gdef\profonditastruttura{\dp\strutbox}
\def\senondefinito#1{\expandafter\ifx\csname#1\endcsname\relax}
\def\BOZZA{
\def\alato(##1){
 {\vtop to \profonditastruttura{\baselineskip
 \profonditastruttura\vss
 \rlap{\kern-\hsize\kern-1.2truecm{$\scriptstyle##1$}}}}}
\def\galato(##1){ \gwidth=\hsize \divide\gwidth by 2
 {\vtop to \profonditastruttura{\baselineskip
 \profonditastruttura\vss
 \rlap{\kern-\gwidth\kern-1.2truecm{$\scriptstyle##1$}}}}}
\def\verapagina{
{\romannumeral\number\numcap}.\number\numsec.\number\numpag}}

\def\alato(#1){}
\def\galato(#1){}
\def\veroparagrafo{\number\numsec}\def\veraformula{\number\numfor}
\def\veraappendice{\number\numapp}
\def\verapagina{\number\pageno}\def\veranformula{\number\numnf}
\def\verafigura{{\romannumeral\number\numcap}.\number\numfig}
\def\verocapitolo{\number\numcap}\def\veranformula{\number\numnf}
\def\veratheorema{\number\numtheo}
\def\Eqn(#1){\eqno{\etichettan(#1)\alato(#1)}}
\def\eqn(#1){\etichettan(#1)\alato(#1)}
\def\TH(#1){{\etichettat(#1)\alato(#1)}}
\def\thv(#1){\senondefinito{fu#1}$\clubsuit$#1\else\csname fu#1\endcsname\fi}
\def\thu(#1){\senondefinito{e#1}\thv(#1)\else\csname e#1\endcsname\fi}
\def\ver{\veroparagrafo}
\def\Eq(#1){\eqno{\etichetta(#1)\alato(#1)}}
\def\eq(#1){\etichetta(#1)\alato(#1)}
\def\Eqa(#1){\eqno{\etichettaa(#1)\alato(#1)}}
\def\eqa(#1){\etichettaa(#1)\alato(#1)}
\def\dgraf(#1){\getichetta(#1)\galato(#1)}
\def\drif(#1){\retichetta(#1)}

\def\eqv(#1){\senondefinito{fu#1}$\clubsuit$#1\else\csname fu#1\endcsname\fi}
\def\equ(#1){\senondefinito{e#1}\eqv(#1)\else\csname e#1\endcsname\fi}
\def\graf(#1){\senondefinito{g#1}\eqv(#1)\else\csname g#1\endcsname\fi}
\def\rif(#1){\senondefinito{r#1}\eqv(#1)\else\csname r#1\endcsname\fi}
\def\bib[#1]{[#1]\numpag=\pgn
\write13{\string[#1],\verapagina}}

\def\include#1{
\openin13=#1.aux \ifeof13 \relax \else
\input #1.aux \closein13 \fi}

\openin14=\jobname.aux \ifeof14 \relax \else
\input \jobname.aux \closein14 \fi
\openout15=\jobname.aux
\openout13=\jobname.bib


\ifnum\tipoformule=1\let\Eq=\eqno\def\eq{}\let\Eqa=\eqno\def\eqa{}
\def\equ{}\fi


{\count255=\time\divide\count255 by 60 \xdef\hourmin{\number\count255}
        \multiply\count255 by-60\advance\count255 by\time
   \xdef\hourmin{\hourmin:\ifnum\count255<10 0\fi\the\count255}}

\def\oramin{\hourmin }

\def\data{\number\day/\ifcase\month\or january \or february \or march \or
april \or may \or june \or july \or august \or september
\or october \or november \or december \fi/\number\year;\ \oramin}

\def\titdate{ \ifcase\month\or January \or February \or March \or
April \or May \or June \or July \or August \or September
\or October \or November \or December \fi \number\day, \number\year;\ \oramin}


\newcount\pgn \pgn=1
\def\foglio{\number\numsec:\number\pgn
\global\advance\pgn by 1}
\def\foglioa{A\number\numsec:\number\pgn
\global\advance\pgn by 1}

\footline={\rlap{\hbox{\copy200}}\hss\tenrm\folio\hss}


\global\newcount\numpunt

\magnification=\magstephalf
\baselineskip=16pt
\parskip=8pt

\voffset=2.5truepc
\hoffset=0.5truepc
\hsize=6.1truein
\vsize=8.4truein 
{\headline={\ifodd\pageno\rightheadline \else \leftheadline \fi}}
\def\rightheadline{\it  {Potential theory on the hypercube}\hfil\tenrm\folio}
\def\leftheadline{\tenrm \folio \hfil\it  {Section $\ver$}}

\def\a{\alpha}
\predefine\barunder{\b}
\redefine\b{\beta}
\def\d{\delta}
\def\e{\epsilon}

\def\f{\phi}
\def\g{\gamma}

\def\s{\sigma}
\def\t{\tau}

\def\o{\omega}
\def\D{\Delta}
\def\L{\Lambda}
\def\G{\Gamma}
\def\O{\Omega}

\def\Th{\Theta}
\def\del #1{\frac{\partial^{#1}}{\partial\l^{#1}}}

\def\1{{1\kern-.25em\roman{I}}}
\def\eu{{1\kern-.25em\roman{I}}}
\def\f1{{1\kern-.25em\roman{I}}}

\def\R{{\Bbb R}}  
\def\N{{\Bbb N}}  
\def\P{{\Bbb P}}  
\def\Q{{\Bbb Q}}  
\def\C{{\Bbb C}}  
\def\E{{\Bbb E}}  

\def\del{\partial}

\def\dist{\,\roman{dist}}

\let\cal=\Cal
\def\AA{{\cal A}}

\def\EE{{\cal E}}

\def\GG{{\cal G}}
\def\HH{{\cal H}}
\def\II{{\cal I}}

\def\KK{{\cal K}}

\def\MM{{\cal M}}

\def\OO{{\cal O}}

\def\QQ{{\cal Q}}
\def\RR{{\cal R}}
\def\SS{{\cal S}}

\def\VV{{\cal V}}
\def\UU{{\cal U}}
\def\VV{{\cal V}}
\def\WW{{\cal W}}

\def\A{{\cal A}}

\def\chap #1#2{\line{\ch #1\hfill}\numsec=#2\numfor=1\numtheo=1}

\def\wt{\widetilde}
\def\wh{\widehat}


\def\note#1{\footnote{#1}}

\def\frac#1#2{{#1\over #2}}
\def\sfrac#1#2{{\textstyle{#1\over #2}}}

\def\text#1{\quad{\hbox{#1}}\quad}
\def\newpage{\vfill\eject}
\def\proposition #1{\noindent{\thbf Proposition #1}}

\def\theo #1{\noindent{\thbf Theorem {#1} }}

\def\lemma #1{\noindent{\thbf Lemma {#1} }}
\def\definition #1{\noindent{\thbf Definition {#1} }}

\def\corollary #1{\noindent{\thbf Corollary #1 }}
\def\proof{{\noindent\pr Proof: }}
\def\proofof #1{{\noindent\pr Proof of #1: }}
\def\endproof{$\diamondsuit$}
\def\remark{\noindent{\bf Remark: }}
\def\thanks{\noindent{\bf Acknowledgements: }}

\font\pr=cmbxsl10

\font\thbf=cmbxsl10 scaled\magstephalf

\font\ch=cmbx12

\font\it=cmti10
\font\bf=cmbx10
\font\sm=cmr7



\input epsf

\def\capa{\hbox{\rm cap}}

\font\tit=cmbx12
\font\aut=cmbx12

\def\s{\char'31}
\centerline{\tit Elementary potential theory on the hypercube.}
\vskip.2truecm
\vskip1truecm
\vskip.5cm

\centerline{\aut G\'erard Ben Arous\note{
EPFL, INR 011 (B\^atiment INR) Station 14 CH-1015 Lausanne Switzerland.\hfill\break
email:
gerard.benarous\@epfl.ch
}
\note{and
Courant Institute for the Mathematical Sciences,
New York university, 251 Mercer Street, New York, NY 10012-1185
}
, V\'eronique  Gayrard\note{CPT-CNRS, Luminy, Case 907, F-13288 Marseille Cedex 9, France
\hfill\break
e-mail: veronique\@gayrard.net,
Veronique.Gayrard\@cpt.univ-mrs.fr}
\note{Work partially supported by the Swiss National Fund.}
}


\vskip.5cm

\vskip0.5truecm\rm
\def\s{\sigma}
\noindent {\bf Abstract:}
This work addresses potential theoretic questions for the standard
nearest neighbor random walk on the hypercube $\{-1,+1\}^N$.
For a large class of subsets $A\subset\{-1,+1\}^N$ we give precise estimates for
the harmonic measure of $A$, the mean hitting time of $A$, and the Laplace transform
of this hitting time.
In particular, we give precise sufficient conditions for the harmonic measure to be
asymptotically
uniform, and for the hitting time to
be asymptotically exponentially distributed, as $N\rightarrow\infty$.
Our approach relies on a $d$-dimensional extension of the Ehrenfest urn scheme
called lumping
and covers the case where $d$ is allowed to diverge with $N$
as long as $d\leq  \a_0\frac{N}{\log N}$ for some constant $0<\a_0<1$.




 \vfill
$ {} $

\newpage

\chap{1. Introduction }1

\bigskip
\line{\bf 1.1. Motivation.\hfill}

This work addresses potential theoretic questions for the standard
nearest neighbor random walk on the hypercube $\{-1,+1\}^N$
(or equivalently on $\{0,+1\}^N$). We will write
$\SS_N\equiv \{-1,+1\}^N$ and generically call $\s=(\s_1,\dots,\s_N)$
an element of $\SS_N$. This random walk $\left(\s_N(t)\right)_{t\in \N}$
is a Markov chain and is described by the following
transition probabilities: for $\s,\s'\in \SS_N$,
$$
p_N(\s,\s'):=\P(\s_N(t+1)=\s'\mid \s_N(t)=\s)=\frac 1N
\Eq(I.1)
$$
if and only if $\s$ and $\s'$ are nearest neighbor on $\SS_N$, i\.e\.
if and only if the Hamming distance
$$
\Dist(\s,\s'):=\#\Bigl\{i\in\{1,\dots,N\}\,:\,\s_i\neq\s'_i\Bigr\}
\Eq(I'.2)
$$
is equal to one.
The questions we are interested in  are related to processes of
Hamming distances on $\SS_N$.
For a non empty subset $L\in\{1,\dots,N\}$ define the Hamming distance in $L$
by
$$
\Dist_L(\s,\s'):=\#\{i\in L\,:\,\s_i\neq\s'_i\}
\Eq(I'.3)
$$
Let $\L$ be a partition of $\{1,\dots,N\}$ into
$d$ classes, that is non-empty disjoint subsets $\L_1,\dots,\L_d$, $1\leq d\leq N$, satisfying
$\L_1\cup\dots\cup\L_d=\{1,\dots,N\}$. We will often call such a partition a $d$-partition.
Given a $d$-partition  $\L$ and a point $\xi\in\SS_N$, we
define the associated "multi-radial" process, i\.e\. the process of distances
$$
D^{\L,\xi}(\s_N(t))=
\left(D^{\L,\xi}_1(\s_N(t)),\dots,D^{\L,\xi}_d(\s_N(t))\right)
\Eq(I'.4)
$$
where, for each $1\leq k\leq d$,
$$
D^{\L,\xi}_k(\s)=\Dist_{\L_k}(\s,\xi)
\,,\quad\s\in\SS_N
\Eq(I'.5)
$$
$D^{\L,\xi}(\s_N(t))$ is a Markov chain on a
subset of $\{0,\dots,N\}^d$ that has cardinality smaller than $2^N$.
The main goal of this paper is to give a detailed analysis of the behavior of
this chain asymptotically, when $N\rightarrow\infty$,
with minimal assumptions on the sizes of the sets $\L_1,\dots,\L_d$ and on the
number $d$ of such sets.

The case where $d=1$ and $\L$ is the trivial partition, i\.e\. where
$D^{\L,\xi}(\s_N(t))$ simply is Hamming distance,
$D^{\L,\xi}(\s_N(t))=\Dist_{\L}(\s_N(t),\xi)$, as been
extensively studied. This process can be traced back to Ehrenfest model
of heat exchange
(we refer to [DGM] for a survey of the early literature).
More recently it was used as an important tool
to understand, for instance,  the rate at which the random walk $\s_N(t)$
approaches equilibrium and the associated ``cut-off phenomenon''
(see Aldous [A1-A2], Aldous and Diaconis [AD1-AD2], Diaconis [D], Diaconis et al\. [DGM],
Saloff-Coste [SaCo],  Matthews [M2-M3], Voit  [V]).
In [D]-[DGM] a major role was played by the Fourier-Krawtchouk transform
(i\.e\. harmonic analysis on the group $\{0,1\}$). We will not rely on
this powerful tool for our study of the case $d>1$ (though it might turn
out to be useful for improving our very rough Theorem \thv(T.theo.2)).

A main
motivation for the study of \eqv(I'.4) with $d>1$ comes from statistical
mechanics of mean-field spin glasses. In this context
the maps $D^{\L,\xi}(\s_N(t))$ are used in an equivalent form, namely, we set
$$
\g^{\L,\xi}(\s_N(t))=
\left(\g^{\L,\xi}_1(\s_N(t)),\dots,\g^{\L,\xi}_d(\s_N(t))\right)
\Eq(I'.6)
$$
where, for each $1\leq k\leq d$,
$$
\g^{\L,\xi}_k(\s)
=\frac{1}{|\L_k|}\sum_{i\in\L_k}\s_i\xi_i
=1-\frac{2}{|\L_k|}\Dist_{\L_k}(\s,\xi)
\,,\quad\s\in\SS_N
\Eq(I'.7)
$$
The chain $\g^{\L,\xi}(\s_N(t))$ now takes value in a discrete grid
$\G_{N,d}$ of  $[-1,1]^d$ that contains the set $\SS_d=\{-1,1\}^d$.
This $d$-dimensional process was exploited  for the study of the dynamics of the random field Curie-Weiss model in [BEGK1],
of the the Random Energy Model (REM) in [BBG1,BBG2], and in [G2] for the study of the dynamics of the
Generalized Random Energy Model (GREM). While some of the results presented
here are refinements of results previously obtained in [BBG1],
the present paper should answer all the needs of the more demanding study
of the GREM dynamics.
(Lumping was also used in the context of large deviation theory in
 to treat  the Hopfield model of a mean-field spin glass \cite{KP,G1,BG}).


Note that in statistical mechanics the map $\g^{\L,\xi}$
has a very natural interpretation: it is a coarse graining map that replaces
detailed information on the $N$ microscopic spins, $\s_i$,
by information on a smaller number $d$ of macroscopic blocks-magnetization,
$\g^{\L,\xi}_k(\s)$.
This type of construction, that maps the chain $\s_N(t)$ into a new Markov
chain $\g^{\L,\xi}(\s_N(t))$ whose state space  $\G_{N,d}$
has smaller cardinality is called {\it lumping} in [KS], and the chain
$X^{\L,\xi}_{N}(t)=\g^{\L,\xi}(\s_N(t))$ is called the  {\it lumped chain}.



\vfill\eject
\noindent{\bf The lumped chain.}
Let us now give an informal description of some of the result we obtain
for the lumped chain $X^{\L,\xi}_{N}(t)$
(or equivalently $D^{\L,\xi}(\s_N(t))$).
The behavior of this chain is better understood if one sees it as a
discrete analogue of a diffusion in a convex potential which is very steep
near its boundary.
This potential is given by the entropy produced by the map $\g^{\L,\xi}$,
i\.e\. by
$\psi^{\L,\xi}_N(x)=-\frac{1}{N}\log{(\g^{\L,\xi})^{}}^{-1}(x)+\log 2$.
It achieves its minimum at the origin and its maximum on $\SS_d=\{-1,1\}^d$,
the $2^d$ vertices of $[-1,1]^d$. We give precise sufficient
conditions for the hitting time of subsets $I$ of $\SS_d$ to be
asymptotically exponentially distributed, and for the hitting
distribution to be uniform.

These conditions essentially require
that $I$ should be sparse enough (see Definition \thv(I.def.3)), and that the partition $\L$
does not contain too many small boxes $\L_k$ (which would give flat directions to the potential).
In order to prove these facts we rely
on the following scenario, which would be classical in any large deviation
approach {\it \`a la} Freidlin and Wentzell [FW]: the lumped chain starts by
going down the potential well; it reaches the origin before reaching any
vertex, then returns many time to the origin before finding a vertex
$x\in\SS_d$ with almost uniform distribution.

To implement this scenario we use
two key ingredients
given in Theorem \thv(P.prop.1)
and Theorem \thv(P.prop.2) respectively. In Theorem \thv(P.prop.1), we
consider the
probability that, starting from  a point in $\SS_d$, the lumped chain
reaches the origin without returning to its starting point. We give sharp
estimates  that show that this ``no return probability''
is close to one. In Theorem \thv(P.prop.2) we give an upper
bound on the probability of the (typically) rare event that consists in hitting a
point $x\in\SS_d$ before hitting the origin, when the chain starts from
an arbitrary point $y$ in $\G_{N,d}$. This bound is given as a function
$F_{N,d}$ of the distance between $x$ and $y$ on the grid $\G_{N,d}$.
This function $F_{N,d}$ is explicit but, unfortunately, pretty involved.
It will be used to describe the necessary sparseness of the sets $I$.


Our approach is based on the tools developped in [BEGK1] and [BEGK2]
for the study of metastability of reversible Markov chains with discrete
state space.
An important fact is that they allow us to deal
with the case where $d$ diverges with $N$, as long as
$d\leq d_0(N):=\lfloor \a_0\frac{N}{\log N}\rfloor$ for some
constant $0<\a_0<1$
\note{Here $\lfloor t\rfloor$, $t\in\R$, denotes the largest integer dominated by $t$.}
\note{
The constant $\a_0$ (which initially arises from Theorem \thv(P.prop.1)) is chosen is such a way
that the $d$-partition $\L$ is log-regular (see Definition \thv(I.def.2)), i\.e\. that
in total, the volume of subsets $\L_k$ of size smaller that $10\log N$ is at most
$N/2$. This condition partly motivates the appearance of the
logarithm in the definition of $d_0(N)$ }.
(This condition can be slightly relaxed but to no great gain.)
A large deviation approach {\it \`a la}
Freidlin and Wentzell would seem problematic at least when $d\geq \sqrt N$.
In this case the fact that the lumped chain lives on a discrete grid cannot
be ignored.

\noindent{\bf The random walk on the hypercube.}
Let us see how the lumped chain can be used to solve potential theoretic
questions for
some
subsets of the hypercube. Given a subset $A$ of
$\SS_N$ consider the hitting time
$
\t_{A}:=\inf\left\{t>0\mid  \s_{N}(t)\in A\right\},
$
and  the hitting point
$
\s_{N}(\t_{A})
$
for the random walk $\s_{N}(t)$ started in $\s$. We wish to find
sufficient conditions that ensure that the distribution of the
hitting time is asymptotically exponentially distributed,
and the distribution of the hitting point asymptotically
uniform, as $N\rightarrow\infty$.

When $A$ contains a single point, or a pair of points,
Matthews [M1] showed that the distribution of the hitting time is
asymptotically exponentially distributed when $N\rightarrow\infty$.
For related results when $A$ is a random set of points in $\SS_N$
see [BBG1] Proposition 2.1, and [BC] Proposition 6.
%
%
Our study
of the lumped chain will enable us to tackle these potential theoretic
questions for a special class of sets $A\subset\SS_N$.

\definition{\TH(I.def.1)}{\it  A subset $A$ of $\SS_N$ is called
$(\L,\xi)$-compatible if and only if $\g^{\L,\xi}(A)\subseteq\SS_d$.
}

Since each point in $\SS_d$ has only one pre-image by $\g^{\L,\xi}$ then
obviously, when $A$ is $(\L,\xi)$-compatible, the hitting time
$
\t_{A}:=\inf\left\{t>0\mid  \s_{N}(t)\in A\right\}
$
is equal to the hitting time $\t_{\g^{\L,\xi}(A)}$ for the lumped chain,
$$
\t_{\g^{\L,\xi}(A)}:=\inf\left\{t>0\,\big|\,X^{\L,\xi}_{N}(t)\in \g^{\L,\xi}(A)\right\},
\Eq(I'.9)
$$
and $\s_{N}(\t_{A})$ is the unique point in
${(\g^{\L,\xi})^{}}^{-1}(\g^{\L,\xi}(\s_{N}(\t_{A})))$.
Our results for hitting times and hitting points for the lumped chain will
thus be transferred  directly  to the random walk on the hypercube.

It is thus important to understand what sets $A\subset\SS_N$ can be described
as $(\L,\xi)$-compatible. For a given pair $(\L,\xi)$ define the set
$B(\L,\xi)\in\SS_N$ by
$$
\s\in B(\L,\xi)
\iff
\Dist_{\L_k}(\s,\xi)\in \{1,\dots,|\L_k|\}\,\text{for all} k\in\{1,\dots,d\}
\Eq(I'.10)
$$
The set $B(\L,\xi)$ is thus made of the $2^d$ points in $\SS_N$ obtained by
a global change of sign of the coordinates of $\xi$ in each of the subsets $\L_1,\dots,\L_d$.
$B(\L,\xi)$ can also be seen as the orbit of $\xi$
under the action of the abelian group of isometries of the hypercube generated
by the $(s^{|\L_k|})_{1\leq k\leq d}$, with
$$
s^{|\L_k|}(\s)_i=\cases +\s_i, &\hbox{if}\,\,\,i\notin\L_k\cr
 -\s_i, &\hbox{if}\,\,\,i\in\L_k\cr
\endcases
\Eq(I'.11)
$$
A set $A\subset\SS_N$ is then $(\L,\xi)$-compatible if and only if it is
included in the orbit $B(\L,\xi)$.

It is easy to see that any set is $(\L,\xi)$-compatible for the partition $\L$
where $\L_k=\{k\}$ for each $1\leq k\leq d$. But in this case $d=N$ and our
results on the lumped chain obviously do not apply.
On the other extreme it is easy to see that small sets (say sets of cardinality $|A|$
smaller than $\log_2 N$) are compatible with partitions into $d$ classes for
$d\leq 2^{|A|}$ (see Lemma \thv(A4.lemma.1) of appendix A4).


\bigskip
\line{\bf 1.2. A selection of results.\hfill}

In the remainder of this introduction we give a more detailed
account of some of the results we can obtain for these potential
theoretic questions taking the view point of the random walk on the hypercube.
In the body of the paper most results will be stated both
for the lumped chain $X^{\L,\xi}_{N}(t)$ and the random walk
$\s_N(t)$ on $\SS_N$.



We start by making precise the condition of {\it sparseness} under which our results obtain,
and introduce the so-called {\it Hypothesis $H$} -- a  minimal distance assumption between
the points of subsets of $\SS_N$.
To introduce the notion of sparseness of a set we need the following definition.
If $A$ is a subset of $\SS_N$  define
$$
\UU_{N,d}(A)
:=\cases \max_{\eta\in A}\sum_{\s\in A\setminus\eta}F_{N,d}(\Dist(\eta,\s)),
  &\text{if} |A|>1\cr
0,&\text{if} |A|\leq 1
\endcases
\Eq(I.12)
$$
where $F_{N,d}$ is a function depending on $N$ and $d$, whose
definition is stated in \eqv(P.p2.17bis)-\eqv(P.p2.01)
of Chapter 3,
and whose properties are analyzed in detail in Appendix A3 (see in particular Lemma \thv(A3.lemma.1)).

\definition{\TH(I.def.3) {\bf (Sparseness)}}{\it A set $A\subset\SS_N$ is called
{\it $(\e,d)$-sparse} if there exists $\e>0$ such that
$$
\UU_{N,d}(A)\leq\e
\Eq(I.13)
$$
}


In Appendix A4 we give explicit estimates on $\UU_{N,d}(A)$ that allow to quantify
the sparseness of $(\L,\xi)$-compatible sets that are either small enough
(Lemma \thv(I.lemma.3)) or whose elements satisfy a minimal distance assumption
(Lemma \thv(A4.lemma.3)).
We now give a few selected examples of sparseness estimates to illustrate our
quantitative notion of sparseness, first in the simplest possible case i\.e\. the
case of equipartition, then a few more examples to show that our notion of
sparseness does not prevent the possibility for a sparse set to have many nearby
points or even many nearest neighbors. Finally we show that if the minimal distance
is large a set can be sparse and still grow exponentially.
Without loss of generality we take $\xi=(1,\dots,1)$.

\noindent{\it Example 1: Equipartition.} Let $d\leq\frac{N}{\log N}$. Assume that $N/d\in\N$
and let $\L$ be any $d$-partition satisfying $|\L_k|=N/d$ for all $1\leq k\leq d$. Then
for all $(\L,\xi)$-compatible set $A$ there exists $\sqrt{2/e}\leq\varrho<1$ such that
$$
\UU_{N,d}(A)\leq \varrho^{N/d}
\Eq(I.ex1)
$$

\noindent{\it Example 2: Many nearby points.} Assume that $d$ satisfies
$\frac{d^{2+\d_0}}{N}<1$ for some $\d_0>0$. Let $\L$ be any $d$-partition satisfying
$|\L_k|=1$ for all but one index $k'\in\{1,\dots,d\}$, and $|\L_{k'}|=N-d+1$.
Then for all $\d\leq \d_0$ and for all  $(\L,\xi)$-compatible set $A$,
$$
\UU_{N,d}(A)\leq \left(\frac{d^{2+\d_0}}{N}\right)^{6/\d}
\Eq(I.ex2)
$$

\noindent{\it Example 3: Many nearest neighbors.} Let $d\leq(\e N)^{3/4}$ for some $\e>0$.
Given $\eta^1\in\SS_N$, let  $\{\eta^2,\dots,\eta^d\}$
be $d-1$ nearest neighbors of $\eta^1$ (i\.e\. for each $1< k\leq d$, $\Dist(\eta^1,\eta^k)=1$)
and set $A=\{\eta^1,\dots,\eta^d\}$. Then
$$
\UU_{N,d}(A)\leq \e
\Eq(I.ex3)
$$
(Note that $A$ is compatible with a partition of the type described in Example 2.)

\noindent{\it Example 4: Many far away points.} Let $0<\d_0, \d_1<1$ be constants
chosen in such a way that the set $A\subset\SS_N$ satisfies the following conditions:
(1) $|A|\geq e^{\d_0 N}$,
(2) $A$ is compatible with a partition  into $d=\lfloor \d_1 N\rfloor$ classes, and
(3)
$
\inf_{\eta\in A}\Dist(\eta,A\setminus\eta)\geq C d
$
for some $C\geq 1$ (\note{To construct such a set start from an equipartition into $d=\lfloor \d_1 N\rfloor$ classes
and select all points compatible with this partition that satisfy the third condition.}).
Then there exists $\varrho<1$ and $0<\d_3<1$ such that
$$
\UU_{N,d}(A)\leq \varrho^{\d_3 N}
\Eq(I.ex4)
$$

The bounds \eqv(I.ex1) and \thv(I.ex2) are easily worked out using the estimates
on $F_{N,d}$ from Lemma \thv(A3.lemma.1) of Appendix A3. Example 3 is derived from
Lemma \thv(I.lemma.3) and Example 4 from Lemma \thv(A4.lemma.3).



Our next condition is concerned with the `local' geometric structure of sets $A\subset\SS_N$:

\definition{\TH(I.def.4){\bf (Minimal distance assumption or hypothesis $H$)}}{\it
We will say that a set
$A\in\SS_N$ obeys hypothesis $H(A)$ (or simply that $H(A)$ is satisfied) if
$$
\inf_{\s\in A}\Dist(\s,A\setminus \s)>3
\Eq(P.H)
$$
}

We will treat both the cases of sets that obey, and of sets that do not obey
this assumption.
When $H(A)$ is not satisfied, our results will be affected by the `local'
geometric structure of a given set $A$. Thus, although our techniques allow in
principle to work out accurate estimates for the objects we are interested in
in this situation also, this must be done case by case.
This local effect is lost as soon as \eqv(P.H) is satisfied and, for arbitrary such sets,
we obtain accurate general results.
Let us stress that this local effects are not (only) a byproduct of our techniques
(see Theorem \thv(TL.theo.2) of Chapter 7 and formulae (3.6), (3.7) in [M1]).
It is not clear however whether the minimal distance in \eqv(P.H) should
not be $2$ or $1$ rather than $3$.

We now proceed to state our results. Let us state here once and for all that all
of them must be preceded by the sentence:
\item{}{\it ``There exists a constant $0<\a_0<1$ such that,
setting $d_0(N):=\lfloor \a_0\frac{N}{\log N}\rfloor$, the following holds.''}

\noindent To further simplify the presentation we will only consider the case where
$\xi$ in \eqv(I'.6) is the point whose coordinates are all equal to one.
We accordingly suppress all dependence on $\xi$ in the
notation. In particular, subsets $A$ of $\SS_N$
that are $(\L,\xi)$-compatible will be called
$\L$-compatible (or compatible with $\L)$.
Finally the symbols $\L$ and $\L'$ will
always designate partitions of $\{1,\dots,N\}$
into, respectively, $d$ and $d'$ classes. Unless otherwise specified,
we assume that $d\leq d_0(N)$ and $d'\leq d_0(N)$.
Statements of the form

\item{}{\it ``Assume that $A\subset\SS_N$ is $\L$-compatible''}

\noindent must thus be understood as

\item{}{\it``Assume that  $A\subset\SS_N$ is compatible
with some partition $\L$
of $\{1,\dots,N\}$
into $d\leq d_0(N)$ classes''.}

\noindent We can now summarize our results as follows.

\noindent{\bf The harmonic measure.}
Throughout this paper we make the following (slightly abusive) notational convention for hitting times:
given a subset $A\subset\SS_N$ and $\s\in\SS_N$ we let
$$
\eqalign{
\t^{\s}_{A}:=\t_{A} \text{for the chain started in $\s$}
}
\Eq(I.3)
$$
This will enable us to write $\P(\t_{A}=t\mid  \s_N^{}(0)=\s)\equiv\P(\t^{\s}_{A}=t)$
and, more usefully
$$
\eqalign{
\P\left(\t^{\s}_{\eta}<\t^{\s}_{A\setminus\eta}\right)
&\equiv
\P\left(\s_N^{}(\t_A)=\eta\mid  \s_N^{}(0)=\s\right)\,,\quad \eta\in A
}
\Eq(I.4)
$$

Now let $A\subset\SS_N$, and let  $H_A(\s,\eta)$ denote the harmonic measure of A starting from
$\s\in\SS_N\setminus A$, i\.e\.
$$
H_A(\s,\eta)
:=\P\left(\t^{\s}_{\eta}<\t^{\s}_{A\setminus\eta}\right)\,,\quad \eta\in A
\Eq(I.5)
$$

\theo{\TH(I.theo.1)} {\it
Assume that  $A\subset\SS_N$ is $\L$-compatible.
Then, there exist constants $0<c^-,c^+<\infty$ such that the following holds:
for all $0\leq\rho\leq N$, for all $\s$ satisfying $\Dist(\s,A)>\rho$, and for all $\eta\in A$,
we have:
$$
\frac{1}{|A|}(1-c^-\vartheta_{N,d}(A,\rho))
\leq
H_A(\s,\eta)
\leq
\frac{1}{|A|}(1+c^+\vartheta_{N,d}(A,\rho))
\Eq(I.16)
$$
where
$$
\vartheta_{N,d}(A,\rho)=\max\left\{\UU_{N,d}(A), |A|F_{N,d}(\rho+1)\right\}
\Eq(I.17')
$$
}

Together with the explicit estimates on $F_{N,d}$ established in
Appendix A3 and Appendix A4, Theorem \thv(I.theo.1) enable us to deduce that,
asymptotically, for all $\L$-compatible set $A$ which is
sparse enough, the harmonic measure is close to the uniform measure
provided that the starting point $\s$ is located outside some suitably chosen balls
centered at the points of $A$. More precisely define
$$
\WW(A,M):=\left\{\s\in\SS_N\mid \Dist(\s,A)\geq\rho(M)\right\}
\Eq(I.14)
$$
where $\rho(M)\equiv\rho_{N,d}(M)$ is any function defined on the integers
(possibly depending on $N$ and $d$) that satisfies
$$
MF_{N,d}(\rho(M)+1)=o(1)\,,\quad N\rightarrow\infty
\Eq(I.15)
$$
It follows from Lemma \thv(A3.lemma.1) (see also the simpler and more
explicit Lemma \thv(I.lemma.3) and Lemma \thv(A4.lemma.2)) that, under the assumptions
on $A$ of Theorem \thv(I.theo.1), we may always choose $\rho(M)$ in such a way that
\eqv(I.15) holds true for $M=|A|$ while at the same time
$\WW(A,M)\neq\emptyset$.
Thus, for all $\s\in\WW(A,|A|)$, $\vartheta_{N,d}(A,\rho)$ decays to zero as $N\rightarrow\infty$
whenever
$$
\UU_{N,d}(A)=o(1)\,,\quad N\rightarrow\infty
\Eq(I.15bis)
$$
Of course $\WW(A,|A|)$ simply is $\SS_N\setminus A$ for all sets $A$ such that \eqv(I.15)
holds true with $\rho(|A|)=0$. This observation and
Corollary \thv(I.cor.1) trigger the next result.

\corollary{\TH(I.cor.theo.1)} {\it  Let $A\subset\SS_N$ be such that
$2^{|A|}\leq d_0(N)$. Then, for all $\s\notin A$, the harmonic
measure of $A$ starting from $\s$ is, asymptotically,
the uniform measure on $A$: there exist constants $0<c^-,c^+<\infty$ such that,
for all $\eta\in A$,
$$
\frac{1}{|A|}\left(1-\frac{c^-}{(\log N)^2}\right)
\leq
H_A(\s,\eta)
\leq
\frac{1}{|A|}\left(1+\frac{c^+}{(\log N)^2}\right)\,,\quad
\Eq(I.16bis)
$$
}

\noindent{\bf Hitting times.}
Our next theorem is concerned with the mean hitting time of subsets $A\subset\SS_N$.
We will see that the precision of our result depends on whether or not assumption
$H$ is satisfied.

\theo {\TH(I.theo.3)}{\it
Let $d'\leq 2d_0(N)$ and assume that $A\subset\SS_N$ is compatible with a $d'$-partition.

Then for all $\s\notin A$
there exists an integer $d$  satisfying $d'<d\leq 2d'$
such that, if $\UU_{N,d}(A\cup\s)\leq\frac{1}{4}$,
$$
\KK^-
\leq
\E(\t^{\s}_{A})
\leq
\KK^+
\Eq(I.21)
$$
and,  for all $\eta\in A$,
$$
\KK^-
\leq
\E\left(\t^{\s}_{\eta}\mid\t^{\s}_{\eta}<\t^{\s}_{A\setminus\eta}\right)
=
\E\left(\t^{\s}_{A}\mid\t^{\s}_{\eta}<\t^{\s}_{A\setminus\eta}\right)
%
%
\leq
\KK^+
\Eq(I.22)
$$
where $\KK^{\pm}$ are defined as follows:
there exist
constants $0<c^-,c^+<\infty$ such that,
$$
\KK^{\pm}=\frac{2^N}{|A|}\Bigl(1+\frac{1}{N}\Bigr)\Bigl(1\pm c^{\pm}\max\left\{\UU_{N,d}(A\cup\s), \frac{1}{N^k}\right\}\Bigr)
\Eq(I.22bis)
$$
where
$$
k=
\cases 2, &\text{if} H(A\cup\s) \text{is satisfied}\cr
1, &\text{if} H(A\cup\s) \text{is not satisfied.}
\endcases
$$

}

\smallskip
\noindent{\bf Laplace transforms of Hitting times.}
We finally
give estimates for the Laplace transforms
of hitting times. By looking at the
object
$$
\E e^{s\t^\s_{A}}\1_{\{\t^\s_{\eta}<\t^\s_{A\setminus\eta}\}}
\,,\quad s\leq 0
\Eq(I.23)
$$
for $\eta\in A$ and $\s\in\SS_N$,
we will also obtain the joint asymptotic behavior of hitting time
and hitting distribution.


\theo{\TH(I.theo.4nonsparse)}{\it Let $d'\leq d_0(N)/2$ and assume that
$A\subset\SS_N$ is compatible with a $d'$-partition.
Then, for all $0\leq\rho\leq N$, for all $\s$ satisfying $\Dist(\s,A)>\rho$,
there exists an integer $d$  satisfying $d'<d\leq 2d'$
such that, if
$$
\UU_{N,d}(A\cup\s)\leq\frac{\d}{4}
\Eq(I''.24)
$$
for some $0<\d<1$ then the following holds for all $\eta\in A$:
for all $\e\geq\d$, there exists a constant $0<c_\e<\infty$
(independent of $\s, |A|, N$, and $d$) such that, for all $s$ real satisfying
$-\infty<s<1-\e$, and all $N$ large enough we have,
$$
\left|
\E\left(
e^{s{\t^{\s}_{A}}/\E\t^{\s}_{A}}
\1_{\{\t^{\s}_{\eta}<\t^{\s}_{A\setminus\eta}\}}\right)
-\frac{1}{|A|}\frac{1}{1-s}
\right|\leq
\frac{1}{|A|}c_{\e}
\tilde\vartheta_{N,d}(A\cup\s,\rho,k)
\Eq(I''.24')
$$
where
$$
\tilde\vartheta_{N,d}(A\cup\s,\rho,k)=\max\Big\{\UU_{N,d}(A\cup\s),\frac{1}{N^k},|A|F_{N,d}(\rho+1)\Big\}
\Eq(I''.24'')
$$
and
$$
k=
\cases 2, &\text{if} H(A\cup\s) \text{is satisfied}\cr
1, &\text{if} H(A\cup\s) \text{is not satisfied.}
\endcases
\Eq(I''.25)
$$
}

The quantity $\tilde\vartheta_{N,d}(A\cup\s,\rho,k)$ defined in \eqv(I''.24') is very similar to
the quantity $\vartheta_{N,d}(A,\rho)$ that appears in \eqv(I.17') of Theorem \thv(I.theo.1).
Reasoning just as in \eqv(I.14)-\eqv(I.15bis) one can show that
there exists $\rho(M)$ satisfying \eqv(I.15) such that, for all $\s\in\WW(A,|A|)$,
$\tilde\vartheta_{N,d}(A\cup\s,\rho,k)$ decays to zero as $N\rightarrow\infty$
whenever
$$
\UU_{N,d}(A)=o(1)\,,\quad N\rightarrow\infty
$$
From this it will follow that, as $N\rightarrow\infty$,
\item{(i)} ${\t^{\s}_{A}}/\E\t^{\s}_{A}$
converges in distribution to an exponential random variable of mean value one, and that
\item{(ii)} for any finite collection $A_1, A_2,\dots,A_n$ of non empty disjoint subsets of
$A$, the random variables $({\t^{\s}_{A_i}},1\leq i\leq n)$ become
asymptotically independent.

\noindent The specialization of Theorem \thv(I.theo.4nonsparse) to the case where $\s\in\WW(A,|A|)$
and  $\UU_{N,d}(A)=o(1)$ is stated in Section 7 as Theorem \thv(I.theo.4).
Just as Corollary \thv(I.cor.theo.1) was deduced from
Theorem \thv(I.theo.1) we will deduce from it the following result:


\corollary{\TH(I.cor.theo.4)} {\it Let $A\subset\SS_N$ be such that
$2^{|A|}\leq d_0(N)$.
Then, for all $\s\notin A$, the following holds:
for all $\e>0$, there exists a constant $0<c_\e<\infty$
(independent of $\s, |A|, N$, and $d$) such that, for all $s$ real satisfying
$-\infty<s<1-\e$, for $N$ large enough, we have,
$$
\left|
\E\left(
e^{s{\t^{\s}_{A}}/\E\t^{\s}_{A}}
\1_{\{\t^{\s}_{\eta}<\t^{\s}_{A\setminus\eta}\}}\right)
-\frac{1}{|A|}\frac{1}{1-s}
\right|\leq
\frac{1}{|A|}
\frac{c_{\e}}{(\log N)^2}
\Eq(I.26)
$$
}

We will also see in Theorem \thv(TL.theo.2) of Section 7 that,
as was established by Matthews [M1],
a sharper result can be obtained when the set $A$ reduces to a single point.

The rest of this paper is organized as follows.
In Chapter 2 we briefly introduce our
notation and the basic facts about our lumping procedure.
In Chapter 3 we study the two key ingredients needed for the analysis of
the lumped chain, namely Theorem \thv(P.prop.1) and Theorem \thv(P.prop.2),
and introduce the important function $F_{N,d}$.
In Chapter 4 we deduce estimates for the hitting probabilities of the
lumped chain from Theorem \thv(P.prop.1) and Theorem \thv(P.prop.2).
In Chapter 5 we show how the results of Chapter 4 can be lifted to the
hypercube, and give estimates for the harmonic measure of $(\L,\xi)$-compatible
subsets that are sparse enough. In Chapter 6 we give estimates for
the expectation of hitting times and in Chapter 7 for the distribution of
these hitting times (through their Laplace transform).
We also give sufficient conditions for
hitting times and hitting points to be asymptotically independent.

\smallskip
\vfill\eject

\bigskip
\vfill\eject
\chap{2. Lumping.}2

Let $1\leq d<N$. Given a point $\xi\in\SS_N$ and a $d$-partition $\L$
(i\.e\. a partition of $\{1,\dots,N\}$ into $d$ classes, $\L_1,\dots,\L_d$),
let $\g^{\L,\xi}$ be the map  defined in \eqv(I'.7), and let
$\left\{X_N(t)\right\}_{t\in\N}$ be the lumped chain
$X^{\L,\xi}_{N}(t)=\g^{\L,\xi}(\s_N(t))$.

\noindent{\bf Notation and conventions.}
The following notation and assumptions will hold throughout the rest of the paper.
For the sake of brevity we will keep
the dependence on $\L$ and on $\xi$ implicit.
We thus write $\g^{\L,\xi}\equiv\g$
and call this map a {\it $d$-lumping}.
Without loss of generality we may and will assume that $\xi$
is the point whose coordinates are all equal to one.
We will then simply say that the $d$-lumping $\g$ is
{\it generated by the $d$-partition $\L$}
if needs be to refer to the underlying partition $\L$ explicitly.
Similarly, we will write $X^{\L,\xi}_{N}(t)\equiv X_{N}(t)$.
This chain evolves on the grid $\G_{N,d}:=\g(\SS_N)$.
Note that the origin of $\R^d$ does not necessarily belong to $\G_{N,d}$.
This happens if and only if all classes of the partition $\L$ have even
cardinality, in which case the potential function
$\psi_N(x)=-\frac{1}{N}\log\g^{-1}(x)+\log 2$ is minimized at the origin.
By convention we will denote by $0$ (and call {\it zero} or {\it the origin}) any point chosen from the
set where $\psi_N(x)$ achieves its global minimum.
The superscript ${}^{\circ}$ will be used to distinguished
objects defined in the lumped chain setting from their counterparts on the hypercube.
Hence we will denote by $\P^{\circ}$ the law of the lumped chain and by $\E^{\circ}$
the corresponding expectation.
Unless otherwise specified $d$ is any integer satisfying $d< N$.


The next two lemmata are quoted from [BBG1] where their proofs can be found.
The first lists a few basic properties of $\g$.

\lemma {\TH(L.2) (Lemma 2.2 of [BBG1])} {\it
The range of $\g$,
$\G_{N,d}:= \g(\SS_N)$, is a discrete subset of the
$d$-dimensional cube $[-1,1]^d$ and may be described as follows. Let $\{u_k\}_{k=1}^d$ be
the canonical basis of $\R^d$. Then,
$$
x\in\G_{N,d}\Longleftrightarrow
x=\sum_{k=1}^d \left(1-2\frac{n_k}{|\L_k|}\right)u_k
\Eq(2.1.10)
$$
where, for each $1\leq k\leq d$, $0\leq n_k\leq|\L_k|$ is an integer.
Moreover, for each $x\in\G_{N,d}$,
$$
|\{\s\in\SS_N\mid\g(\s)=x\}|
=
\prod_{k=1}^{d}
\binom{|\L_k|}{|\L_k|\frac{1+x_k}{2}},
\Eq(2.1.110)
$$
}


To $\G_{N,d}$ we associate an undirected graph, $\GG(\G_{N,d})=(V(\G_{N,d}),E(\G_{N,d}))$,
with set of vertices $V(\G_{N,d})=\G_{N,d}$ and set of edges:
$$
E(\G_{N,d})=
\left\{(x,x')\in \G_{N,d}\mid
\exists_{k\in\{1,\dots,d\}},\exists_{s\in\{-1,1\}}\,:\,
x'-x=s\sfrac{2}{|\L_k|}u_k\right\}
\Eq(2.1.12)
$$
In the next lemma we summarize the main properties of the lumped chain
$\left\{X_N(t)\right\}$.

\lemma {\TH(L.3) (Lemma 2.3 of [BBG1])}{\it
\item{i)}
$\left\{X_N(t)\right\}$ is Markovian no matter how the initial
distribution $\pi^{0}$ of $\{\s_N^{}(t)\}$ is chosen.
\item{ii)} Set $\Q^{}_N = \mu^{}_N\circ \g^{-1}$ where
$$
\mu^{}_N(\s)=2^{-N},\,\,\,\,\s\in\SS_N
$$
denotes the unique reversible invariant measure for the chain $\{\s_N(t)\}$. Then
$\Q^{}_N$ is the unique reversible invariant measure for the chain
$\left\{X_N(t)\right\}$. In explicit form, the density of $\Q^{}_N$ reads:
$$
\Q^{}_N(x)=
\frac{1}{2^N}|\{\s\in\SS_N\mid\g(\s)=x\}|
,\,\,\,\,\forall x\in\G_{N,d}
\Eq(2.1.13)
$$

\item{iii)}  The transition probabilities
$
r^{}_N(x,x'):=\P^{\circ}(X_N(t+1)=x'\mid X_N(t)=x)
$
of $\left\{X_N(t)\right\}$ are given by
$$
r^{}_N(x,x')=
\cases
\frac{|\L_k|}{N}\frac{1-s x_k}{2}\,
&\text{if} (x,x')\in E(\G_{N,d}))
\text{and} x'-x=s\sfrac{2}{|\L_k|}u_k
\cr
                           0,&\text{otherwise}
\endcases
\Eq(2.1.14)
$$
}


For us the key observation is
the following lemma,
which will allow us to express hitting probabilities, mean hitting times, and Laplace transforms
on the hypercube in terms of their lumped chain counterparts.

\lemma {\TH(L.5)} {\it If $A\subset\SS_N$ is $(\L,\xi)$-compatible then
$$
\t_{A}:=\inf\left\{t>0\mid  \s_{N}(t)\in A\right\}
=\inf\left\{t>0\mid X_{N}(t)\in \g(A)\right\},
\Eq(L.5.1)
$$
and $\s_N(\t_{A})$ is the unique point in $\g^{-1}(X_{N}(\t_{A}))$.
%
}

%
%

\proof The content of this lemma was in fact already sated and proven in
the paragraph following Definition \thv(I.def.1) (see \eqv(I'.9)).
Let us repeat the main line of argument: for each $t\in\N$,
$\s_N(t)\in A$ if and only if $X_{N}(t)\in \g(A)$, which implies that
$\s_N(t)\in\g^{-1}(\g(A))$, and since $A\subset\SS_N$ is $(\L,\xi)$-compatible $\g^{-1}(A)=A$.
\endproof


The next two lemmata are elementary consequences of Lemma \thv(L.6) whose proofs we skip.
Recall that $\P^{\circ}$ denote the law of the lumped chain and  $\E^{\circ}$
the corresponding expectation.

\lemma{\TH(L.6)} {\it Let $A,B\subset\SS_N$ be such that $A\cap B=\emptyset$.
Then, for all $d$-lumping $\g$ compatible with $A\cup B$,
$$
\P\left(\t^{\s}_{A}<\t^{\s}_{B}\right)
=\P^{\circ}\left(\t^{\g(\s)}_{\g(A)}<\t^{\g(\s)}_{\g(B)}\right)\,,
\text{for all}\s\in\SS_N
\Eq(L.6.1)
$$
}


\lemma{\TH(L.7)} {\it Let $A\subset\SS_N$ and $\eta\in A$.
Then for all $\s\in\SS_N$ and all $d$-lumping $\g$ compatible with $A\cup\s$,
if $|A|>1$,
$$
\E\left(
e^{s{\t^{\s}_{A}}/\E\t^{\s}_{A}}
\1_{\{\t^{\s}_{\eta}< \t^{\s}_{A\setminus\eta}\}}\right)
=
\E^{\circ}\left(
e^{s{\t^{\g(\s)}_{\g(A)}}/\E^{\circ}\t^{\g(\s)}_{\g(A)}}
\1_{\{\t^{\g(\s)}_{\g(\eta)}< \t^{\g(\s)}_{\g(A)\setminus\g(\eta)}\}}\right)
\Eq(L.7.1)
$$
and if $A=\{\eta\}$,
$$
\E\left(e^{s{\t^{\s}_{\eta}}/\E\t^{\s}_{\eta}}\right)
=
\E^{\circ}\left(
e^{s{\t^{\g(\s)}_{\g(\eta)}}/\E^{\circ}\t^{\g(\s)}_{\g(\eta)}}
\right)
\Eq(L.7.1bis)
$$
}



We finally state and prove a lemma that will be needed in Section 6 and Section 7.

\lemma {\TH(L.4)} {\it
$
{\displaystyle{
\E\t^0_0=\prod_{k=1}^{d}\sqrt{\frac{\pi}{2}|\L_k|}
\left(1+O\left({|\L_k|^{-1}}\right)\right)
}}$.
}

\proof Since $\left\{X_N(t)\right\}$ is an irreducible
chain on a finite state space whose invariant measure
$\Q^{}_N(x)$ satisfies $\Q(0)>0$, it follows from Kac's Lemma that
$\E\t^0_0{\Q^{}_N(0)}=1$. Lemma \thv(L.4) then follows from
\eqv(2.1.110), \eqv(2.1.13), and Stirling's formula. \endproof


\bigskip
\vfill\eject

\bigskip

\chap{3. The lumped chain: key probability estimates.}3


This chapter centers on the lumped chain.
As noted earlier this chain is a random walk in a simple, convex, potential:
the entropy produced by the lumping procedure gives rise through
\eqv(2.1.13) to a potential
$\psi_N(x)=-\frac{1}{N}\log\g^{-1}(x)+\log 2$
and,
by  Lemma \thv(L.2), this potential
is convex and takes on its global maximum on the set $\SS_d$,
its global minimum being achieved at zero\note{
Recall that by convention the point denoted by $0$ and called {\it zero} or {\it the origin} is any given
point chosen from the set where $\psi_N(x)$ achieves its global minimum: this set
reduces to the actual zero of $\R^d$ if and only if all classes of the partition $\L$ have even cardinality.
}.
Following the strategy developed in [BEGK1], where such chains were investigated,
we will decompose all events
at visits of the chains to zero.
The aim of this chapter is to provide probability estimates for the key events
that will emerge from such decompositions.




Theorem \thv(P.prop.1)
and Theorem \thv(P.prop.2) can be viewed as the two building blocs of this strategy.
Theorem \thv(P.prop.1) establishes that starting from a point $x\in\SS_d$,
the chain finds zero before returning to its starting point
with a probability close to one.

%
\theo{\TH(P.prop.1)}{\it There exists a constant $0<\a_0\leq1/20$ such that if
$d\leq\a_0\frac{N}{\log N}$ then, for all $x\in\SS_d$,
$$
0\leq
\left(1-\frac{1}{N}\right)-\P^{\circ}(\t^x_0<\t^x_x)
\leq \frac{3}{N^2}(1+O(1/{\sqrt N}))
\Eq(P.p1.2)
$$
}


Theorem \thv(P.prop.2) gives an upper bound on the probability that
starting from an arbitrary point $y\in\G_{N,d}$, the chain finds a point
$x\in\SS_d$ before visiting zero. This bound is expressed as
a function $F_{N,d}$ of the  distance between $x$ and $y$ on the grid
$\G_{N,d}$ which guarantees, in particular, that for small enough $d$
this probability decays to zero as $N$ diverges.
Unfortunately, though explicit, the function $F_{N,d}$
looks rather terrible and is not easy to handle.
%
For this reason we state the theorem first and give its definition next.

Given two points $x,y\in \G_{N,d}$, we denote by $\dist(x,y)$ the graph distance
between $x$ and $y$, namely, the number of edges of the
shortest path connecting $x$ to $y$ on the graph $\GG(\G_{N,d})$
(see \eqv(P.l5.4) for the formal definition of a path):
$$
\dist(x,y)\equiv\sum_{k=1}^d\frac{|\L_k|}{2}|x^k-y^k|
\Eq(N.0)
$$
Define $d_0(N)$ as the largest integer
dominated by $\a_0\frac{N}{\log N}$,
where $\a_0$ is the constant appearing in Theorem \thv(P.prop.1).

\theo{\TH(P.prop.2)}{\it Let $d\leq d_0(N)$.
%
%
Then, for all $x\in\SS_d$ and $y\in\G_{N,d}\setminus x$, we have,
with the notation of Definition \thv(P.def.1),
$$
\P^{\circ}\left(\t^y_x<\t^y_0\right)\leq F_{N,d}(\dist(x,y))
\Eq(P.p2.1)
$$
}

From now on (and except in the statement and proofs of the main results of Sections 5 - 7)
we will drop the indices $N$ and $d$ and write $F\equiv F_{N,d}$.
Let us now give the definition of this function.
To this aim let $\del_n x$ be the sphere of radius $n$ centered at $x$,
$$
\del_n x=\{y\in\G_{N,d}\mid \dist(x,y)=n\}\,,\quad n\in\N
\Eq(N.1)
$$

\definition{\TH(P.def.1)}{\it Let
$F, F_1, F_2$, and $\kappa$ be functions, parametrised by $N$ and $d$,
defined on $\{1,\dots,N\}\subset\N$ as follows: let $I(n)$, $n\in$ be the set
defined through
$$
I(n)\equiv\left\{m\in\N^* \mid \exists\, 0\leq p\in\N \: m+2p=n+2\right\}\,;
\Eq(P.p2.17bis)
$$
then
$$
F(n)\equiv F_1(n)+F_2(n)
\Eq(P.p2.02)
$$
where
$$
F_1(n)\equiv \kappa(n)\frac{n!}{N^n}\,,\quad
F_2(n)\equiv  \kappa^2(n+2) \frac{(n+2)!}{N^{(n+2)}}
\sum_{m\in I(n)}\frac{N^{(n+2-m)/2}}{[(n+2-m)/2]!}\left|\del_{m}x\right|
\Eq(P.p2.03)
$$
and
$$
\kappa(n)\leq{
\cases \kappa_0&\text{if $n$ is independent of $N$}\cr
       N&\text{if $n$ is an increasing function of $N$}\cr
\endcases}
\Eq(P.p2.01)
$$
where $1\leq \kappa_0\leq \infty$ is a numerical constant.
}

Lemma \thv(A3.lemma.1) of Appendix A3 contains a detailed analysis of the
large $N$ behavior of the function $F_2$ from \eqv(P.p2.03). There,
we strove to give explicit, easy to handle, expressions that should
meet our needs for all practical purposes.

The rest of this chapter revolves around the proofs of
Theorem \thv(P.prop.1) and Theorem \thv(P.prop.2). However, while the
probabilities dealt with in these two theorems will suffice to express bounds on the
harmonic measure, more general `no-return before hitting' probabilities  than that of
Theorem \thv(P.prop.1)
will
enter a number of our estimates (see e\.g\. the formulae \eqv(T.5),\eqv(T.6),
for hitting times).
%
%
%
%
Therefore, anticipating our future needs, we divide the chapter in five sections as follows.
We first establish upper bounds on `no return before hitting' probabilities of
the general form $\P^{\circ}(\t^x_J<\t^x_x)$ for $J\subset\SS_d$ and
$x\in \G_{N,d}\setminus J$ (Lemma \thv(P.lemma.3)), from which we will deduce the upper
bound on $\P^{\circ}(\t^x_0<\t^x_x)$ with $x\in\SS_d$ (Corollary \thv(P.cor.4)) needed
to prove Theorem \thv(P.prop.1).
This is the content of Section 3.1.
In Section 3.2 we prove a lower bound on `no return' probabilities
of the form $\P^{\circ}(\t^x_0<\t^x_x)$ (Lemma \thv(P.lemma.7)) which is rather rough
but holds uniformly for $x\in\G_{N,d}\setminus 0$. This general a priori upper
bound will be needed in the proof of Theorem \thv(P.prop.2),
carried out in Section 3.3. We will in fact prove a slightly stronger version of
Theorem \thv(P.prop.2), namely Theorem \thv(P.prop.2stronger), valid for
all under the only
assumption that the partition $\L$ is log-regular (see Definition \thv(I.def.2)).
Theorem \thv(P.prop.2) is in turn  needed to obtain
the sharp upper bound on  $\P^{\circ}(\t^x_0<\t^x_x)$
of Theorem \thv(P.prop.1), which we next prove in Section 3.4.


\bigskip
\line{\bf 3.1. Upper bounds on `no return before hitting' probabilities.\hfill}

Given a a subset $J\subset\SS_d$ and a point $x\in \G_{N,d}\setminus J$, consider
the probability $\P^{\circ}(\t^x_J<\t^x_x)$ that the lumped chain hits $J$
before returning to its starting point.
Our general strategy to bound  these `no return'  probabilities
is drawn from [BEGK1,3]
and summarized in Appendix A1. It hinges
on the fact that they admit of a  variational representation
(stated in Lemma \thv(A.lemma.1)) which is nothing but an analogue for our reversible
Markov chains of the classical Dirichlet principle from potential theory.
This  variational representation enables us to derive upper bounds in a very simple
way, simply guessing the minimizer. It will also allow us to obtain lower bounds
by
comparing the initial problem
to a sum of
one dimensional problems (Lemma \thv(A.lemma.2))
that, as we will see in Section 3.2,
can be worked out explicitly with good precision.


We now focus on the upper bounds problem for $d>1$ only,
the case $d=1$ being covered in Lemma \thv(A.lemma.3). These bounds
will be obtained under the condition that the set
$J\cup x$ obeys hypothesis $H^{\circ}$, namely, under the condition
that
$$
\inf_{z\in J\cup x}\dist(z,(J\cup x)\setminus z)>3
\Eq(P.H'lump)
$$
This is a transposition in the lumped chains setting
of  hypothesis $H$ initially defined in \eqv(P.H) for subsets $A$
of the hypercube $\SS_N$
(we will see in Chapter 5 that for certain sets, conditions
\eqv(P.H'lump) and \eqv(P.H) are equivalent).
Naturally we will say that $J\cup x$ obeys hypothesis
$H^{\circ}(J\cup x)$ (or simply that $H^{\circ}(J\cup x)$ is satisfied)
whenever \eqv(P.H'lump) is satisfied.


\lemma{\TH(P.lemma.3)}{\it  Let $J\subset\SS_d$, $x\in \G_{N,d}\setminus J$,
and set
$$
\a_J\equiv \frac{\Q_N(J)}{\Q_N(x)}\,,\quad
\b\equiv\frac{\Q_N(\del_1 x)}{\Q_N(x)}-1\,,\quad
\d_J=\frac{|\del_1x\cap J|}{N}
\Eq(P.l3.1)
$$

\item{i)} If
$H^{\circ}(J\cup x)$ is satisfied then
$$
\P^{\circ}(\t^x_J<\t^x_x)\leq
\frac
{\a_J\left(1-\frac{1}{N}\right)}
{1+\a_J\left(1-\frac{1}{N}\right)\left(1+\frac{1}{\b}\right)}
\Eq(P.l3.2)
$$

\item{ii)} If $H^{\circ}(J\cup x)$ is not satisfied,
and if $x\in\SS_d\setminus J$, then
$$
\P^{\circ}(\t^x_J<\t^x_x)\leq
\frac{\a_J}{1+\a_J}\Big(1+\frac{2\d_J}{1+\a_J}+O\Big(\frac{\d^2_J}{\a_J}\Big)\Big)
\Eq(P.l3.3')
$$
}

\remark The condition \eqv(P.H'lump)
is the weakest condition we could find that yields a very accurate upper
bound\note{We will actually work out the corresponding lower bound
(see \eqv(P.t0.11) of  Theorem \thv(P.theo.0)).} on
$\P^{\circ}(\t^x_J<\t^x_x)$ that is independent of the geometry of the set
$J\cup x$. In contrast, \eqv(P.l3.3') is the roughest bound we could derive but
depends only in a mild way on the  geometry of $J$ near $x$. As we will
explain, our techniques allow to work out better (albeit often
inextricable) bounds.

\remark Observe that Lemma \thv(P.lemma.3) holds with no assumption on the
cardinality of $J$.

\remark Also observe that since $|\del_1x|=d$,
$
\d_J
\leq\frac{d}{N}\1_{\{|J|\geq d\}}+\frac{|J|}{N}\1_{\{|J|\leq d\}}
\leq 1
$.
Since for $x\in\SS_d$, $\a_J=|J|$, we have
$$
\frac{\d_J}{\a_J}\leq
\frac{d}{N|J|}\1_{\{|J|\geq d\}}+\frac{1}{N}\1_{\{|J|\leq d\}}
\leq \frac{1}{N}\,,
\Eq(P.l3.3'')
$$
and \eqv(P.l3.3') may be bounded by
$$
\P^{\circ}(\t^x_J<\t^x_x)\leq
\left(1-\sfrac{1}{|J|+1}\right)\left(1+O(\sfrac{1}{N})\right)
\Eq(P.l3.3''')
$$


\proof Assume first that $H^{\circ}(J\cup x)$ is satisfied.
Using the
variational representation of Lemma \thv(A.lemma.1) we may write
$$
\P^{\circ}(\t^x_J<\t^x_x) =
\Q_N(x)^{-1}\inf_{h\in \HH^y_J}\Phi_{N,d}(h)
\leq \Q_N(x)^{-1}\Phi_{N,d}(h)\,,\quad \forall\,h\in\HH^x_J
\Eq(P.l3.3)
$$
where  $\Phi_{N,d}$ is defined in \eqv(A.6). We then choose
$$
h(y)=\cases 0,&\text{if}\,y=x\cr
                1,&\text{if}\,y\in J\cr
                a,&\text{if}\,y\in \del_1 J\cr
                c,&\text{if}\,y\in \del_1 x\cr
                b,&\text{if}\,y\notin (J\cup x)\cup(\del_1 J\cup \del_1 x)\cr
         \endcases
\Eq(P.l3.4)
$$
where $a, b$, and $c$ are still to be determined.
Inserting this ansatz into $\Phi_{N,d}$, we see that
$$
\eqalign{
\Phi_{N,d}(h)
&\equiv \sum_{y\in J}\left\{
\sum_{y'\in\del_1 y}\Q_N(y)r^{\circ}_N(y,y')(1-a)^2
+\sum_{y'\in\del_1 y}
\sum_{y''\in(\del_1 y')\setminus y}
\Q_N(y')r^{\circ}_N(y',y'')(a-b)^2
\right\}
\cr
&+
\left\{
\sum_{y'\in\del_1 x}\Q_N(x)r^{\circ}_N(x,y')(c-0)^2
+\sum_{y'\in\del_1 x}
\sum_{y''\in(\del_1 y')\setminus x}
\Q_N(y')r^{\circ}_N(y',y'')(b-c)^2
\right\}
}
\Eq(P.l3.5)
$$
To evaluate the various sums in the last formula simply observe that, for all
$z\in\G_{N,d}$,
$$
\eqalign{
\sum_{z'\in\del_1 z}\Q_N(z)r^{\circ}_N(z,z')
=& \Q_N(z)\sum_{z'\in\del_1 z}r^{\circ}_N(z,z')
=\Q_N(z)
\cr
\sum_{z'\in\del_1 z}\sum_{z''\in(\del_1 z')\setminus z}
\Q_N(z')r^{\circ}_N(z',z'')
=&
\sum_{z'\in\del_1 z}\Q_N(z')
\sum_{z''\in(\del_1 z')\setminus z}r^{\circ}_N(z',z'')
\cr
=&
\sum_{z'\in\del_1 z}\Q_N(z')
\Bigl(\sum_{z''\in\del_1 z'}r^{\circ}_N(z',z'')-r^{\circ}_N(z',z)\Bigr)
\cr
=&
\sum_{z'\in\del_1 z}\Q_N(z')
\left(1-r^{\circ}_N(z',z)\right)
\cr
=&
\left(\Q_N(\del_1 z)-\Q_N(z)\right)
\cr
}
\Eq(P.l3.6)
$$
where the last line follows from reversibility. Also observe that when
$y\in\SS_d$,
$$
\Q_N(\del_1 y)=N\Q_N(y)
\Eq(P.l3.7)
$$
Then, using \eqv(P.l3.6) and \eqv(P.l3.7) in \eqv(P.l3.5),
we get
$$
\eqalign{
\Phi_{N,d}(h)
=& \Q_N(J)\left[(1-a)^2+(N-1)(a-b)^2\right]
+\Q_N(x) c^2+\left(\Q_N(\del_1 x)-\Q_N(x)\right)(b-c)^2
}
\Eq(P.l3.8)
$$
and by \eqv(P.l3.3), for $\a\equiv\a_J$ and $\b$ defined
in \eqv(P.l3.1),
$$
\P^{\circ}(\t^x_J<\t^x_x)\leq
\a\left[(1-a)^2+(N-1)(a-b)^2\right]+c^2+\left(\b-1\right)(b-c)^2
\Eq(P.l3.9)
$$
This allows us to determine $a,b$, and $c$ by maximizing the right hand side of
\eqv(P.l3.9): one easily finds that the maximum is attained at $a=a^*$, $b=b^*$,
$c=c^*$, where
$$
\eqalign{
a^*=&1-\frac{c^*}{\a}\cr
b^*=&c^*\left(\sfrac{1+\b}{\b}\right)\cr
c^*=&\a\left(\sfrac{1}{N-1}+1+\a\sfrac{1+\b}{\b}\right)^{-1}\cr
}
\Eq(P.l3.9')
$$
Plugging these values into \eqv(P.l3.9) then yields \eqv(P.l3.2).

If now $H^{\circ}(J\cup x)$ is not satisfied,
the test function $h(y)$ of
\eqv(P.l3.4) is no longer suitable: we can either add extra parameters to handle the pairs
$y',y''\in J\cup x$ that are such that $\dist(y',y'')\leq 3$, or simplify the form of $h(y)$ by, e\.g\.,
suppressing all but one of the parameters. Clearly the first option should yield more accurate
results, but these results will strongly depend on the local structure of the $J\cup x$,
and in practice this will be tractable only when this structure is given explicitly.
Instead, we choose the one parameter test function
$$
h(y)=\cases 0,&\text{if}\,y=x\cr
                1,&\text{if}\,y\in J\cr
                a,&\text{otherwise}\cr
     \endcases
\Eq(P.l3.4')
$$
Eq\. \eqv(P.l3.8) then becomes
$$
\eqalign{
\Phi_{N,d}(h)
=&
\left[\Q_N(J)-\Q_N(y)\d_J\right](1-a)^2
+\Q_N(x)(1-\d_J)a^2+\Q_N(y)\d_J
\cr
=&\Q_N(x)\left[
\left(\a_J-\d_J\right)(1-a)^2
+(1-\d_J)a^2+\d_J
\right]
}
\Eq(P.l3.8')
$$
where we used in the last line that since $x,y\in\SS_d$, $\Q_N(y)=\Q_N(x)$.
One then readily gets that the maximum in \eqv(P.l3.8') is attained at
$a=a^*\equiv (\a_J-\d_J)/(1+\a_J-\d_J)$,
and takes the value
$$
\Q_N(x)
\frac{\a_J}{1+\a_J}
\frac{1-\d_J^2/\a_J}{1-2\d_J/(1+\a_J)}
$$
From here \eqv(P.l3.3') follows immediately.
%
\endproof

A few immediate consequences of Lemma \thv(P.lemma.3) are collected below.

\corollary{\TH(P.cor.4)}{\it Set
$
\b_0\equiv 2d\left(1-\frac{1}{d}
\sum_{k=1}^d\frac{1}{|\L_k|/2+1}\right)-1
$.
Then, for all $1\leq d<N$,
\item{i)} for all $J\subset\SS_d$ and $x\in\SS_d\setminus J$,
if $H^{\circ}(J\cup x)$ is satisfied,
$$
\P^{\circ}(\t^x_J<\t^x_x)\leq
\left(1-\sfrac{1}{|J|+1}\right)\left(1-\sfrac{1}{N}\right)\,,
\Eq(P.c4.1)
$$
whereas, if $H^{\circ}(J\cup x)$ is not satisfied,
$$
\P^{\circ}(\t^x_J<\t^x_x)\leq
\left(1-\sfrac{1}{|J|+1}\right)\left(1+O(\sfrac{1}{N})\right)
\Eq(P.c4.1bis)
$$
\item{ii)} for all $J\subset\SS_d$,
$$
\P^{\circ}(\t^0_J<\t^0_0)\leq |J|2^{-N}\left(1-\sfrac{1}{N}\right)
\left(1-|J|2^{-N}\left(1-\sfrac{1}{N}\right)\left(1+\sfrac{1}{\b_0}\right)
+\OO(|J|^22^{-2N})
\right)
\Eq(P.c4.2)
$$
\item{iii)} for all $x\in\SS_d$,
$$
\P^{\circ}(\t^x_0<\t^x_x)\leq \left(1-\sfrac{1}{N}\right)
\left(1-2^{-N}\left(1-\sfrac{1}{N}\right)\left(1+\sfrac{1}{\b_0}\right)
+\OO(2^{-2N})
\right)
\Eq(P.c4.3)
$$
}

\proofof{Corollary \thv(P.cor.4)} All we have to do is to evaluate the
coefficients $\a_J$, $\b$, and $\d_J$ of \eqv(P.l3.1), and to decide which of the
formula \eqv(P.l3.2) or \eqv(P.l3.3') to use.
Clearly, \eqv(P.c4.1) and \eqv(P.c4.2) of the corollary satisfy assumption (i) of
Lemma \thv(P.lemma.3), so that \eqv(P.l3.2) apply in both these cases, while
\eqv(P.c4.1bis) satisfy assumption (ii).
Now, when $J\subset\SS_d$ and $x\in\SS_d\setminus J$,
$$
\a_J= |J|\,,\quad
\b= N-1\,,\quad
\frac{\d_J}{\a_J}\leq \frac{1}{N}
$$
where the bound on ${\d_J}/{\a_J}$ was established in \eqv(P.l3.3'').
Inserting these values in \eqv(P.l3.2), respectively \eqv(P.l3.3'),
yields \eqv(P.c4.1), respectively  \eqv(P.c4.1bis). This proves
assertion (i). To prove (ii) note that when $x=0$,
$$
\Q_N(\del_1 x)=2d\Q_N(0)
\left(1-\frac{1}{d}\sum_{k=1}^d\frac{1}{|\L_k|/2+1}\right)
\Eq(P.c4.4)
$$
and hence
$$
\b\equiv\b_0=2d\left(1-\frac{1}{d}\sum_{k=1}^d\frac{1}{|\L_k|/2+1}\right)-1
\Eq(P.c4.5)
$$
On the other hand, for $J\subset\SS_d$ and $x=0$, $\a_J=|J|2^{-N}$.
Then, plugging these values into \eqv(P.l3.2) and setting
$u=|J|2^{-N}\left(1-\sfrac{1}{N}\right)\left(1+\sfrac{1}{\b_0}\right)$
yields
$$
\P^{\circ}(\t^0_J<\t^0_0)\leq(1+\sfrac{1}{\b_0})^{-1}\frac{u}{1+u}
\Eq(P.c4.6)
$$
This and the bound $\frac{1}{1+u}=1-u+\frac{u^2}{1+u}\leq 1-u+u^2$ for
$u>0$ proves assertion (ii) of the corollary. To prove the
last assertion we use reversibility to write
$$
\P^{\circ}(\t^x_0<\t^x_x)
=\frac{\Q_N(0)}{\Q_N(x)}
\P^{\circ}(\t^0_x<\t^0_0)
=\a_{\{x\}}^{-1}\P^{\circ}(\t^0_x<\t^0_0)
\Eq(P.c4.7)
$$
The bound \eqv(P.c4.3) then follows from \eqv(P.c4.2) with $J=\{x\}$ and
$\a_{\{x\}}=2^{-N}$. This concludes the proof of Corollary \thv(P.cor.4).
\endproof

The last assertion of Corolloray \thv(P.cor.4) proves the upper bound of
\eqv(P.p1.2). We may now turn to the corresponding lower bound.

\vfill\eject
\bigskip
\line{\bf 3.2. Lower bound on general probabilities of
`no return' before hitting $0$. 
\hfill}


This subsection culminates in Lemma \thv(P.lemma.7)
which provides an upper bound on the probability
that the lumped chain hits the origin
(i\.e\.the global minimum of the potential well
$\psi_N(y)=-\frac{1}{N}\log\Q^{}_N(y)$))
without returning to its starting point $x$, for arbitrary
$x\in\G_{N,d}\setminus 0$. While so far we made no assumption on the
partition $\L$,  Lemma \thv(P.lemma.7) holds provided that $\L$ be log-regular
(see Definition \thv(I.def.2)), i\.e\. that it does not
contain too many small boxes $\L_k$ (which would give flat
directions to the potential).
We will see, comparing \eqv(P.p1.2) and \eqv(P.l7.5),
that the latter bound is rather rough
and can at best yield the correct leading order when $x\in\SS_d$.


The proof of  Lemma \eqv(P.lemma.7) proceeds as follows.
Using Lemma \thv(A.lemma.2) of Appendix A1 we can bound the
`no return' probability
$\varrho_{N,d}(x)\equiv\P^{\circ}(\t^x_0<\t^x_x)$ of a $d$-dimensional chain
with the help of similar quantities, $\varrho_{|\L_k|,d=1}(x^k)$, $1\leq k\leq d$,
but defined in a $1$-dimensional setting. This
is the content of Lemma \thv(P.lemma.5). The point of
doing this is that, as stated in Lemma \thv(A.lemma.3) of Appendix A2,
such one dimensional probabilities
can be worked out explicitly with good precision.
It will then only remain to combine the results of the previous two lemmata.
This is done in Lemma \thv(P.lemma.7) under the assumption that
the partition $\L$ is log-regular .


%
\def\mo{\hbox{\sm mod}}
%

\lemma{\TH(P.lemma.5)}{\it Set
$$
\varrho_{N,d}(x)\equiv\P^{\circ}(\t^x_0<\t^x_x)
\,,\quad x\in \G_{N,d}
\Eq(P.l5.1)
$$
Then, writing $x=(x^1,\dots,x^d)$,
$$
\varrho_{N,d}(x)\geq \sum_{\mu=1}^d
\left[
%
\sum_{\nu=0}^{d-1}
\varepsilon^{(\mu)}_{\nu}(x)
\frac{N}{|\L_{(\mu+\nu)\mo_d}|}
\varrho^{-1}_{|\L_{(\mu+\nu)\mo_d}|, 1}(x^{(\mu+\nu)\mo_d})
\right]^{-1}
\Eq(P.l5.2)
$$
where
$$
\eqalign{
&\varepsilon^{(\mu)}_{\nu}(x)
\equiv\cases 1,&\text{if}\,\nu=0\cr
\displaystyle{
  \prod_{k=1}^{\mu+\nu}q(|\L_{k \mo_d}|, x^{k \mo_d}) }
         ,&\text{if}\,1\leq \nu\leq d-1\cr
  \endcases
\cr
&\cr
&q(|\L_{k}|, x^{k})\equiv
{\binom{|\L_{k}|}{|\L_{k}|\frac{1+x^{k}}{2}}}\Big/
{\binom{|\L_{k}|}{|\L_{k}|\frac{1}{2}}}
}
\Eq(P.l5.3)
$$
and, for $d=1$, $\varrho_{N,1}(0)\equiv0$.
}

\proofof{Lemma \thv(P.lemma.5)}  An $L$-steps path $\o$ on $\G_{N,d}$,
 beginning at  $x$ and ending at  $y$ is defined as sequence of $L$ sites
$\o=(\o_0,\o_1,\dots,\o_L)$, with $\o_0=x$, $\o_L=y$, and
$\o_l=(\o_l^k)_{k=1,\dots,d}\in V(\G_{N,d})$ for all $1\leq l\leq d$,
that satisfies:
$$
(\o_l,\o_{l-1})\in E(\G_{N,d}), \,\,\,\,\text{ for all} \,l=1,\dots,L
\Eq(P.l5.4)
$$
(We may also write $|\o|=L$ to denote the length of $\o$.)

Since $\varrho_{N,1}(0)\equiv0$ we may assume without loss of generality
that in \eqv(P.l5.1), $x=(x^1,\dots,x^d)$ is such that $x^k>0$ for all
$1\leq k\leq d$. There is no loss of generality either in assuming that, for
each $1\leq k\leq d$, $|\L_k|$ is even.
With this we introduce $d$ one-dimensional paths in $\G_{N,d}$,
connecting $x$ to the endpoint $0$, each being of length
$$
L=L_1+\dots+L_d\,,\quad L_k\equiv\frac{|\L_k|}{2}x^k\,.
\Eq(P.l5.5)
$$

\definition {\TH(P.def)} {\it
 Set $L_0\equiv 0$ and
let $\o=(\o_0,\dots,\o_n,\dots,\o_L)$, $\o_n=(\o_n^k)_{k=1}^d$, be the path
defined by
$$
\o_0= x
\Eq(P.l5.6)
$$
and, for
$L_0+\dots+L_i+1\leq n \leq L_0+\dots+L_{i+1}$,
$0\leq i \leq d-1$,
$$
\o_n^k=\cases
0,
&\text{if} k<i+1\cr
x^k-\frac{2}{|\L_k|}n,
&\text{if} k=i+1
\cr
x^k &\text{if}  k>i+1\cr
\endcases
\Eq(P.l5.7)
$$
For $1\leq\mu\leq d$ let
$\pi_{\mu}$ be the permutation of \{1\,\dots,d\} defined by
$\pi_{\mu}(k)=(\mu+k-1)\mo_d$. Then, for each
$1\leq\mu\leq d$,  we denote by
$\o(\mu)=(\o_0(\mu),\dots,\o_n(\mu),\dots,\o_L(\mu))$,
$\o_n(\mu)=(\o_n^k(\mu))_{k=1}^d$, the path in $\G_{N,d}$ defined
through
$$
\o_n^k(\mu)=\o_n^{\pi_{\mu}(k)}(\mu)
\Eq(P.l5.8)
$$
for $L_0+L_{\pi_{\mu}(1)}+\dots+L_{\pi_{\mu}(i)}+1\leq n \leq
L_0+L_{\pi_{\mu}(1)}+\dots+L_{\pi_{\mu}(i+1)}$,
$0\leq i \leq d-1$.
}

Thus, the path $\o$ defined by \eqv(P.l5.6) and \eqv(P.l5.7) consists of a
sequence of straight pieces along the coordinate axis, starting with the
first and ending with the last one, and deacreasing each coordinate to zero
(all steps in the path "goes down".) In the same way, the path $\o(\mu)$
of \eqv(P.l5.8) follows the axis but, this time, in the permuted order
$\pi_{\mu}(1), \pi_{\mu}(2), \dots, \pi_{\mu}(d)$.


Now, for each $1\leq\mu\leq d$, let  $\D_{\mu}$ the subgraph of
$\GG(\G_{N,d})$ ``generated'' by
the path $\o(\mu)$, i.e., having for set of vertices the set
$V(\D_{\mu})=\{x'\in\G_{N,d}\mid
\exists_{0\leq n\leq L}\,:\,x'=\o_n(\mu)\}$. Clearly the collection
$\D_{\mu}$, $1\leq\mu\leq d$,
verifies conditions \eqv(A.2) and \eqv(A.3) of Lemma \thv(A.lemma.2).
It then follows from the latter that
$$
\P^{\circ}\left(\t^x_0<\t^x_x\right)\geq
\sum_{\mu=1}^d
\wt\P^{\circ}_{\D_{\mu}}
\left(\t^{\o_0(\mu)}_{\o_L(\mu)}<\t_{\o_0(\mu)}^{\o_0(\mu)}\right)
\Eq(P.l5.9)
$$
so that Lemma \thv(P.lemma.5) will be proven if we can
establish that:

\lemma{\TH(P.lemma.6)}{\it
Under the assumptions and with the notation of Lemma \thv(P.lemma.5) and
Definition \eqv(P.def), for each $1\leq\mu\leq d$,
$$
\wt\P^{\circ}_{\D_{\mu}}
\left(\t^{\o_0(\mu)}_{\o_L(\mu)}<\t_{\o_0(\mu)}^{\o_0(\mu)}\right)
=
\left[
%
\sum_{\nu=0}^{d-1}
\varepsilon^{(\mu)}_{\nu}(x)
\frac{N}{|\L_{(\mu+\nu)\mo_d}|}
\varrho^{-1}_{|\L_{(\mu+\nu)\mo_d}|, 1}(x^{(\mu+\nu)\mo_d})
\right]^{-1}
\Eq(P.l5.10)
$$
}

\proofof{Lemma \thv(P.lemma.6)}
Without loss of generality we may assume that
$\mu=1$, in which case the path $\o(\mu)$ coincides with $\o$, and
\eqv(P.l5.10) reads
$$
\wt\P^{\circ}_{\D_{1}}
\left(\t^{\o_0}_{\o_L}<\t_{\o_0}^{\o_0}\right)
=
\left[
\sum_{\nu=0}^{d-1}
\varepsilon^{(1)}_{\nu}(x)
\frac{N}{|\L_{\nu}|}
\varrho^{-1}_{|\L_{\nu}|, 1}(x^{\nu})
\right]^{-1}
\Eq(P.l5.11)
$$
where, $\varepsilon^{(1)}_{0}(x)\equiv 1$ and for $1\leq \nu\leq d-1$,
$
\varepsilon^{(1)}_{\nu}(x)
\equiv
\prod_{k=1}^{\nu}q(|\L_{k}|, x^{k})\,,
$
and $q(|\L_{k}|, x^{k})$ is defined in \eqv(P.l5.3).
As we have stressed several times already, the point of \eqv(P.l5.9) is that
each of the $d$ chains appearing in the r\.h\.s\. evolves in a one dimensional
state space, and that in dimension one last passage probabilities are well
known (see e\.g\. [Sp]). We recall that
$$
\wt\P^{\circ}_{\D_{1}}
\left(\t^{\o_0}_{\o_L}<\t_{\o_0}^{\o_0}\right)
=
\left[
\sum_{n=1}^{L}
\frac{\wt\Q^{\circ}_{\D_{1}}(\o_0)}{\wt\Q^{\circ}_{\D_{1}}(\o_n)}
\frac{1}{\wt r^{\circ}_{\D_{1}}(\o_n,\o_{n-1})}
\right]^{-1}
\Eq(P.l5.13)
$$
which we may also write, using reversibility together with the
definitions of  $\wt r^{\circ}_{\D_{1}}$ and
$\wt\Q^{\circ}_{\D_{1}}$ (see Appendix A1),
$$
\eqalign{
\wt\P^{\circ}_{\D_{\mu}}
\left(\t^{\o_0}_{\o_L}<\t_{\o_0}^{\o_0}\right)
=
\left[
\sum_{n=0}^{L-1}
\frac{\Q_N(\o_0)}{\Q_N(\o_n)}
\frac{1}{r_N(\o_n,\o_{n+1})}
\right]^{-1}
=
\left[
\sum_{\nu=0}^{d-1}
A_{\nu}
\right]^{-1}
\cr
}
\Eq(P.l5.14)
$$
where, setting $L_0=0$,
$$
A_{\nu}=\sum_{n=L_{\nu}}^{L_{\nu+1}-1}
\frac{\Q_N(\o_0)}{\Q_N(\o_n)}
\frac{1}{r_N(\o_n,\o_{n+1})}
\Eq(P.l5.15)
$$
Now for
$L_{\nu}\leq n\leq L_{\nu+1}-1$ we have, on the one hand,
$$
\eqalign{
\frac{\Q_N(\o_0)}{\Q_N(\o_n)}
=&
\prod_{k=1}^{d}
\frac{\Q_{|\L_k|}(\o_0^{k})}{\Q_{|\L_k|}(\o_n^{k})}
\cr
=&
\left(
\prod_{k=1}^{\nu}
\frac{\Q_{|\L_k|}(x^{k})}{\Q_{|\L_k|}(0)}
\right)
\frac{\Q_{|\L_{\nu+1}|}(\o_0^{\nu+1})}{\Q_{|\L_{\nu+1}|}(\o_n^{\nu+1})}
\left(
\prod_{k=\nu+2}^{d}
\frac{\Q_{|\L_k|}(x^{k})}{\Q_{|\L_k|}(x^{k})}
\right)
\cr
=&
\left(
\prod_{k=1}^{\nu}q(|\L_{k}|, x^{k})
\right)
\frac{\Q_{|\L_{\nu+1}|}(\o_0^{\nu+1})}{\Q_{|\L_{\nu+1}|}(\o_n^{\nu+1})}
\cr
=&
\varepsilon^{(1)}_{\nu}(x)
\frac{\Q_{|\L_{\nu+1}|}(\o_0^{\nu+1})}{\Q_{|\L_{\nu+1}|}(\o_n^{\nu+1})}
\cr
}
\Eq(P.l5.16)
$$
where the one before last line follows from  \eqv(2.1.110), \eqv(2.1.13), and
the definition \eqv(P.l5.3) of $q(|\L_{k}|, x^{k})$.
On the other hand, 
$$
r_N(\o_n,\o_{n+1})
=
\frac{|\L_{\nu+1}|}{N}r_{|\L_{\nu+1}|}(\o_n^{\nu+1},\o_{n+1}^{\nu+1})
\Eq(P.l5.17)
$$
where
$r_{|\L_{\nu+1}|}(\,.\,,\,.\,)$ are the rates of the one
dimensional lumped chain $X_{|\L_{\nu+1}|}(t)$ with state space
$\G_{|\L_{\nu+1}|,1}$.
Inserting \eqv(P.l5.16) and \eqv(P.l5.17) in \eqv(P.l5.15) yields
$$
\eqalign{
A_{\nu}
&=
\frac{N}{|\L_{\nu+1}|}
\varepsilon^{(1)}_{\nu}(x)
\sum_{n=L_{\nu}}^{L_{\nu+1}-1}
\frac{\Q_{|\L_{\nu+1}|}(\o_0^{\nu+1})}{\Q_{|\L_{\nu+1}|}(\o_n^{\nu+1})}
\frac{1}
{r_{|\L_{\nu+1}|}(\o_n^{\nu+1},\o_{n+1}^{\nu+1})}
\cr
}
\Eq(P.l5.18)
$$
and, in view of formula \eqv(P.l5.13) (or equivalently  \eqv(P.l5.14))
$$
\eqalign{
A_{\nu}
=&
\frac{N}{|\L_{\nu+1}|}
\varepsilon^{(1)}_{\nu}(x)
\varrho^{-1}_{|\L_{\nu+1}|, 1}(x^{\nu+1})
\cr
}
\Eq(P.l5.19)
$$
where, with the notation of \eqv(P.l5.1), $\varrho^{-1}_{|\L_{\nu+1}|, 1}(x')$
is the last passage probability
$
\varrho^{-1}_{|\L_{\nu+1}|, 1}(x')
\equiv\P^{\circ}(\t^{x'}_0<\t^{x'}_{x'})
$
for the one dimensional lumped chain $X_{|\L_{\nu+1}|}(t)$.
%
Inserting \eqv(P.l5.19) in \eqv(P.l5.14) proves \eqv(P.l5.11).
Lemma \thv(P.lemma.6) is thus proven.\endproof

Inserting \eqv(P.l5.10) in \eqv(P.l5.9) yields \eqv(P.l5.2), and
concludes the proof of Lemma \thv(P.lemma.5).\endproof






Combining Lemma \thv(P.lemma.5) and the one dimensional estimates of
Lemma \thv(A.lemma.3) readily yields upper bounds on last passage probabilities.
We expect these bounds to be reasonably good when the contribution of the
terms $\varepsilon^{(\mu)}_{\nu}(x)$ with  $\nu>0$ in  \eqv(P.l5.2) remains
negligible. Inspecting \eqv(P.l5.3), one sees that this will be the case when
$x$ is far enough from zero (i\.e\. away from the global minimum of the potential
$\psi_N(x)=-\frac{1}{N}\log\Q^{}_N(x)$) and thus, even more so when $x$ is close to
$\SS_d$, provided however that the partition $\L$ does not contain too many small
boxes $\L_k$, i\.e\, provided that the partition $\L$ is log-regular. We now make this
condition precise:

\definition{\TH(I.def.2){\bf (Log-regularity)}}{\it
A partition $\L$ into $d$ classes $(\L_1,\dots,\L_d)$ is called log-regular
if there exists $0\leq\a\leq 1/2$ such that
$$
\sum_{\mu=1}^d|\L_{\mu}|\1_{\{|\L_{\mu}|<10\log N\}}\leq \a N
\Eq(P.good.1)
$$
}
We will call $\a$ the rate of regularity. Note that if $\L$ is log-regular
there exists at least one index $1\leq \mu\leq d$
such that $|\L_{\mu}|\geq 10\log N$ (since supposing that  $|\L_{\mu}|< 10\log N$
for all $1\leq \mu\leq d$ implies that $\sum_{\mu=1}^d|\L_{\mu}|<N$).
Also note that a necessary condition for $\L$ to be  log-regular is that
$d<\a'N$ for some $1\geq \a'\geq\a$ (more precisely $\a'\equiv\a'(N)\sim\a$ as
$N\uparrow\infty$)
while, clearly, all partitions into $d\leq\frac{\a}{10}\frac{N}{\log N}$ classes are
log-regular with rate $\a$.

%
%
%

\lemma{\TH(P.lemma.7)}{\it
For all fixed integer $n$, for all
$x\in\G_{N,d}$ such that $\dist(x,\SS_d)=n$,
if the partition $\L$ is log-regular with rate $\a$,
then
$$
\P^{\circ}(\t^x_0<\t^x_x)\geq 1-\a-\frac{C}{\log N}
\Eq(P.l7.5)
$$
where $0<C<\infty$ is a numerical constant.
Moreover, for all $x\in\G_{N,d}\setminus 0$,
$$
\P^{\circ}(\t^x_0<\t^x_x)\geq
\frac{c}{N}\left[\frac{1}{d}\sum_{\nu=1}^{d}\frac{1}{\sqrt{|\L_{\nu}|}}\right]^{-1}
\Eq(P.l7.5')
$$
where $0<c<\infty$ is a numerical constant.
}

\remark
Eq\. \eqv(P.l7.5') implies that
$\P^{\circ}(\t^x_0<\t^x_x)\geq\frac{c}{N}$. In the case
when $|\L_{\mu}|=\frac{N}{d}(1+o(1))$ ({\it i\.e\.} when
all boxes have roughly same size) \eqv(P.l7.5') implies
that $\P^{\circ}(\t^x_0<\t^x_x)\geq c\frac{1}{\sqrt{dN}}$.

\remark We see, comparing \eqv(P.p1.2) and \eqv(P.l7.5), that
choosing
$\a=o(1)$
in  \eqv(P.l7.5) yields the correct  leading term.

\proof Eq\. \eqv(P.l7.5) is a byproduct of the following more general statement:
for $\L$ a log-regular $d$-partition  with  rate $\a$,
set
$
\II\equiv\{k\in\{1,\dots,d\}\mid |\L_k|\geq 10\log N\}
$
(hence $|\II|\neq\emptyset$); then, defining the set
$$
\wt\G_{N,d}\equiv
\left\{x\in \G_{N,d}\,\Big|\,
\sup_{k\in\II}q(|\L_{k}|, x^{k})
=o\left(\sfrac{1}{N^5}\right)
\right\}
\Eq(P.l7.2)
$$
we have, for all $x\in\wt\G_{N,d}$,
$$
\P^{\circ}(\t^x_0<\t^x_x)\geq
\left(1-o\left(\sfrac{1}{N^{2}}\right)\right)
{\textstyle{
\left(1-\sum_{\mu\notin\II}\frac{|\L_{\mu}|}{N}\right)
}}
\inf_{\mu\in\II}\varrho_{|\L_{\mu}|, 1}(x^{\mu})
\Eq(P.l7.3)
$$
Let us first show that this result implies \eqv(P.l7.5).
Let $k\in\II$. Recall that
$q(|\L_{k}|, x^{k})\equiv\Q_{|\L_k|}(x^{k})/\Q_{|\L_k|}(0)$
where
$\Q_{|\L_k|}(x^{k})=2^{-|\L_k|}\binom{|\L_{k}|}{|\L_{k}|\frac{1+x^{k}}{2}}$.
By Stirling formula,
$$
q(|\L_{k}|, x^{k})=\Q_{|\L_k|}(x^{k})\sqrt{2\pi |\L_k|}
\left(1+O\left(\sfrac{1}{|\L_k|}\right)\right)
\Eq(P.l7.6)
$$
Assume now that $x$ is such that  $\dist(x,\SS_d)=n$ for some integer $n$
independent of $N$. Then,
$
\Q_{|\L_k|}(x^{k})\leq 2^{-|\L_k|}|\L_k|^n
$
and, for $k\in\II$ and for $N$ large enough,
$
\Q_{|\L_k|}(x^{k})\leq e^{-6\log N}
$.
From this and \eqv(P.l7.6) we conclude that $x\in\wt\G_{N,d}$. It
remains to evaluate \eqv(P.l7.3). By \eqv(A.21) of Lemma \thv(A.lemma.3),
$\inf_{\mu\in\II}\varrho_{|\L_{\mu}|, 1}(x^{\mu})\geq
\inf_{\mu\in\II}(1-O(\sfrac{1}{|\L_k|}))
\geq (1-c\sfrac{1}{\log N})
$
for some constant $0<c<\infty$,
whereas $\sum_{\mu\notin\II}|\L_{\mu}|\leq\a N$
(since, by assumption, $\L$ is a log-regular $d$-partition).
It thus follows from \eqv(P.l7.3)
that
$
\P^{\circ}(\t^x_0<\t^x_x)\geq
\left(1-o\left(\sfrac{1}{N^{2}}\right)\right)
\left(1-\a\right)
\left(1-c\sfrac{1}{\log N}\right)
$,
which yields the first line of \thv(P.l7.5).

Let us now prove \eqv(P.l7.3). Let $x\in\wt\G_{N,d}$. By \eqv(P.l5.2),
$$
\varrho_{N,d}(x)\geq \sum_{\mu\in\II}
\left[
\sum_{\nu=0}^{d-1}
\varepsilon^{(\mu)}_{\nu}(x)
\frac{N}{|\L_{(\mu+\nu)\mo_d}|}
\varrho^{-1}_{|\L_{(\mu+\nu)\mo_d}|, 1}(x^{(\mu+\nu)\mo_d})
\right]^{-1}
\Eq(P.l7.7)
$$
From now on let $\mu\in\II$. Since $q(|\L_{\nu}|, x^{\nu})\leq 1$,
$$
\varepsilon^{(\mu)}_{\nu}(x)
\leq q(|\L_{\mu}|, x^{\mu})
\leq o\left(\sfrac{1}{N^{5}}\right)
\Eq(P.l7.8)
$$
where we used in the last bound that $x\in\wt\G_{N,d}$. Using in addition that,
for all $1\leq \nu\leq d$,
$\varrho_{|\L_{\nu}|, 1}(x^{\nu})\leq 1$ and that, by \eqv(A.24)
of Lemma \thv(A.lemma.3),
$\varrho^{-1}_{|\L_{\nu}|, 1}(x^{\nu})\leq C\sqrt{|\L_{\nu}|}\leq C\sqrt N$
for some constant $0<C<\infty$, we get
$$
\eqalign{
&\left[
\sum_{\nu=0}^{d-1}
\varepsilon^{(\mu)}_{\nu}(x)
\frac{N}{|\L_{(\mu+\nu)\mo_d}|}
\varrho^{-1}_{|\L_{(\mu+\nu)\mo_d}|, 1}(x^{(\mu+\nu)\mo_d})
\right]^{-1}
\cr
\geq &
\frac{|\L_{\mu}|}{N}
\varrho_{|\L_{\mu}|, 1}(x^{\mu})
\left[1+Cq(|\L_{\mu}|, x^{\mu})|\L_{\mu}|
\sum_{\nu=1}^{d-1}
\frac{1}{\sqrt{|\L_{(\mu+\nu)\mo_d}|}}
\right]^{-1}
\cr
\geq &
\frac{|\L_{\mu}|}{N}
\varrho_{|\L_{\mu}|, 1}(x^{\mu})
\left[1+CdNq(|\L_{\mu}|, x^{\mu})
\right]^{-1}
\cr
= &
\frac{|\L_{\mu}|}{N}
\varrho_{|\L_{\mu}|, 1}(x^{\mu})
\left[1+o\left(\sfrac{1}{N^{2}}\right)\right]^{-1}
}
\Eq(P.l7.9)
$$
Inserting
\eqv(P.l7.9) in \eqv(P.l7.7),
$$
\eqalign{
\varrho_{N,d}(x)
&\geq
\left(1-o\left(\sfrac{1}{N^{2}}\right)\right)
{\textstyle{\sum_{\mu\in\II}\frac{|\L_{\mu}|}{N}}}
\varrho_{|\L_{\mu}|, 1}(x^{\mu})
\cr
&\geq
{\textstyle{
\left(1-o\left(\sfrac{1}{N^{2}}\right)\right)
\Bigl(\sum_{\mu\in\II}\sfrac{|\L_{\mu}|}{N}\Bigr)
\inf_{\mu\in\II}\varrho_{|\L_{\mu}|, 1}(x^{\mu})
}}
}
\Eq(P.l7.10)
$$
But this is \eqv(P.l7.3).

It remains to prove the last line of \thv(P.l7.5'). Reasoning as in \eqv(P.l7.9)
to bound $\varrho_{|\L_{\nu}|, 1}(x^{\nu})$
but using that $\varepsilon^{(\mu)}_{\nu}(x)\leq 1$ for all $\mu, \nu$,
$$
\left[
\sum_{\nu=0}^{d-1}
\varepsilon^{(\mu)}_{\nu}(x)
\frac{N}{|\L_{(\mu+\nu)\mo_d}|}
\varrho^{-1}_{|\L_{(\mu+\nu)\mo_d}|, 1}(x^{(\mu+\nu)\mo_d})
\right]^{-1}
\geq
\left[CN\sum_{\nu=1}^{d}
\frac{1}{\sqrt{|\L_{\nu}|}}
\right]^{-1}
\Eq(P.l7.11)
$$
From this and \eqv(P.l5.2) we get
$$
\P^{\circ}(\t^x_0<\t^x_x)\geq \frac{c}{N}
\left[\frac{1}{d}\sum_{\nu=1}^{d}
\frac{1}{\sqrt{|\L_{\nu}|}}
\right]^{-1}
\Eq(P.l7.12)
$$
Lemma \thv(P.lemma.7) is now proven.
\endproof


\bigskip
\line{\bf 3.3. Proof of Theorem \thv(P.prop.2).\hfill}

Theorem \thv(P.prop.2) will in fact be deduced from the following stronger statement

\theo{\TH(P.prop.2stronger)}{\it Assume that the $d$-partition $\L$ is log-regular.
Then, for all $x\in\SS_d$ and $y\in\G_{N,d}\setminus x$, we have,
with the notation of Definition \thv(P.def.1),
$$
\P^{\circ}\left(\t^y_x<\t^y_0\right)\leq F(\dist(x,y))
\Eq(P.p2.1stronger)
$$
}

The following useful identity is cited from Corollary 1.6 of [BEGK1]. It will be needed in the proof
of Theorem{ \thv(P.prop.2stronger)} and in different places in the next sections.

\lemma{\TH(P.renew) (Corollary 1.6 [BEGK1])}{\it Let $I\subset\G_{N,d}$. For all $y\notin I$,
$$
\P^{\circ}\left(\t^{y}_{x}<\t^y_I\right)
=
\frac{\P^{\circ}(\t^y_x<\t^y_{I\cup y})}
{\P^{\circ}\left(\t^{y}_{I\cup x}<\t^y_y\right)}
\Eq(P.renew')
$$
}

\proofof{Theorem{ \thv(P.prop.2stronger)}}
Given an integer  $0\leq n\leq N$, let $y$ be a point in $\G_{N,d}$ such that
$\dist(x,y)=n$ where, without loss of generality, we again assume
that $x=(x^1,\dots,x^d)\in\SS_d$ is the vertex whose components all take the value one:
$x=(1,\dots,1)$.
Our starting point then is the relation \eqv(P.renew'),
$$
\P^{\circ}\left(\t^y_x<\t^y_0\right)
= \frac{\P^{\circ}\left(\t^y_x<\t^y_{y\cup 0}\right)}
{\P^{\circ}\left(\t^{y}_{x\cup 0}<\t^y_y\right)}
\Eq(P.p2.10)
$$
To bound the denominator simply note that, by Lemma \thv(P.lemma.7),
$$
\P^{\circ}\left(\t^{y}_{x\cup 0}<\t^y_y\right)
\geq\P^{\circ}\left(\t^{y}_{0}<\t^y_y\right)
\geq \kappa^{-1}(\dist(y,\SS_d))
\geq \kappa^{-1}(n)
\Eq(P.p2.11)
$$
where $\kappa(n)$ is defined in  \eqv(P.p2.01)
(this requires choosing $\kappa_0^{-1}\leq 1-\a(1+o(1))$ for large enough $N$,
which is guaranteed by e\.g\. choosing $\kappa_0^{-1}\leq 1/4$).
To deal with the numerator we first use reversibility to write
$$
\P^{\circ}\left(\t^y_x<\t^y_{y\cup 0}\right)=
\frac{\Q_N(x)}{\Q_N(y)}\P^{\circ}\left(\t^x_y<\t^x_{x\cup 0}\right)
\Eq(P.p2.12)
$$
Since $\dist(x,y)=n$, we may decompose the probability in the
r\.h\.s\. of \eqv(P.p2.12) as
$$
\P^{\circ}\left(\t^x_y<\t^x_{x\cup 0}\right)=
\P^{\circ}\left(n=\t^x_y<\t^x_{x\cup 0}\right)
+\P^{\circ}\left(n<\t^x_y<\t^x_{x\cup 0}\right)
\Eq(P.p2.13)
$$
and set
$$
\eqalign{
f_1&\equiv\P^{\circ}\left(n=\t^x_y<\t^x_{x\cup 0}\right)\cr
f_2&\equiv\P^{\circ}\left(n<\t^x_y<\t^x_{x\cup 0}\right)\cr
}
\Eq(P.p2.14)
$$
Thus, by \eqv(P.p2.10), \eqv(P.p2.11), and \eqv(P.p2.12), defining
$$
F_i\equiv
\frac{\Q_N(x)}{\Q_N(y)}\P^{\circ}\left(\t^{y}_{0}<\t^y_y\right)^{-1}f_i\,,\quad i=1, 2
\Eq(P.p2.21)
$$
Eq\. \eqv(P.p2.13) yields
$$
\P^{\circ}\left(\t^y_x<\t^y_{0}\right)=F_1+F_2
\Eq(P.p2.22)
$$

Let us note here for later use that, by \eqv(P.p2.12),
$
\P^{\circ}\left(\t^y_x<\t^y_{y\cup 0}\right)\leq\Q_N(x)/\Q_N(y)
$
and, combining with \eqv(P.p2.10) and \eqv(P.p2.11),
$$
\P^{\circ}\left(\t^y_x<\t^y_0\right)
\leq \kappa(\dist(y,\SS_d))\frac{\Q_N(x)}{\Q_N(y)}
\Eq(P.p2.12bis)
$$

We now want to bound the terms $f_i$, $i=1, 2$. $f_1$ is the easiest:
For  $z,z'\in\G_{N,d}$, let $r^{(n)}_N(z,z')\equiv\P^{\circ}_z\left(X_N(n)=z'\right)$
be the $n$-steps transition probabilities of the chain $X_{N}$. Then, because
$y$ is exactly $n$ steps away from $x$, we have
$$
f_1
\leq \P^{\circ}\left(n=\t^x_y<\t^x_{x}\right)
=r^{(n)}_N(x,y)
\Eq(P.p2.23)
$$

To bound the term $f_2$
we will decompose it according  to the distance between
the position of the chain at time $n+2$ and its starting point.
Namely, defining the events
$$
\AA_m\equiv\left\{\dist(x,X_N(n+2))=m\right\}\,,\quad m\in\N
\Eq(P.p2.15)
$$
we write
$$
f_2
=\P^{\circ}\left(n+2\leq\t^x_y<\t^x_{x\cup 0}\right)
=\sum_{0<m\leq n+2}\P^{\circ}\left(\left\{n+2\leq\t^x_y<\t^x_{x\cup 0}\right\}
\cap\AA_m\right)
\Eq(P.p2.16)
$$
The only non zero terms in the sum above are those for which $m$ has the same parity as
$n$ or, equivalently, those for which $m$ belongs to the set
$$
I(n)\equiv\left\{m\in\N^* \mid \exists\, 0\leq p\in\N \: m+2p=n+2\right\}
\Eq(P.p2.17)
$$
Next observe that
by the Markov property, if $m\neq n$,
$$
\eqalign{
&\P^{\circ}\left(\left\{n+2\leq\t^x_y<\t^x_{x\cup 0}\right\}\cap\AA_m\right)
\cr
=&
\sum_{z\in\del_mx}
\P^{\circ}_x\left( \left\{X_N(n+2)=z\right\}
\cap\left\{\t^x_{x\cup 0\cup y}>n+2\right\}\right)
\P^{\circ}\left(\t^z_y<\t^z_{x\cup 0}\right)
\cr
\leq&
\sum_{z\in\del_mx}
r^{(n+2)}_N(x,z)
\P^{\circ}\left(\t^z_y<\t^z_{x\cup 0}\right)
\cr
}
\Eq(P.p2.24)
$$
while if $m=n$,
$$
\eqalign{
&\P^{\circ}\left(\left\{n+2\leq\t^x_y<\t^x_{x\cup 0}\right\}\cap\AA_n\right)
\cr
=&
\P^{\circ}\left(n+2=\t^x_y<\t^x_{x\cup 0}\right)
\cr
&+
\sum_{z\in(\del_nx)\setminus y}
\P^{\circ}_x\left( \left\{X_N(n+2)=z\right\}
\cap\left\{\t^x_{x\cup 0\cup y}>n+2\right\}\right)
\P^{\circ}\left(\t^z_y<\t^z_{x\cup 0}\right)
\cr
\leq&
r^{(n+2)}_N(x,y)
+
\sum_{z\in(\del_nx)\setminus y}
r^{(n+2)}_N(x,z)
\P^{\circ}\left(\t^z_y<\t^z_{x\cup 0}\right)
\cr
}
\Eq(P.p2.25)
$$


\lemma{\TH(P.lemma.10)}{\it Let $r^{(n)}_N(\,.\,,\,.\,)$ denote the $n$-steps
transition
probabilities of the chain $X_{N}$.
\item{(i)}
For all $0<n\leq N$,
$$
r^{(n)}_N(x,z)=\frac{n!}{N^n}\frac{\Q_N(z)}{\Q_N(x)}\,,\quad\text{for all}
x\in\SS_d\,,\, z\in\del_nx
\Eq(P.p2.26)
$$
\item{(ii)}  Let $m\in I(n)$ and set $p=(n+2-m)/2$. Then,
$$
r^{(n+2)}_N(x,z)\leq r^{(m)}_N(x,z)\frac{1}{N^p}\frac{(m+2p)!}{m!\,p!}
\,\quad\text{for all}x\in\SS_d\,,\, z\in\del_mx
\Eq(P.p2.27)
$$
}

\proofof{Lemma{ \thv(P.lemma.10)}} We first prove \eqv(P.p2.26).
Without loss of generality we may  assume that $x^k=1$ for all $1\leq k\leq d$.
Given an integer $0\leq n\leq N$, choose $z\in\del_nx$.
Clearly
the set
$\del_nx\subset \G_{N,d}$
is in a one-to-one correspondence with the set $\QQ^{\G}_d(n)$
of integer solutions of the constrained equation
$$
\eqalign{
& n_1+\dots+n_d=n\,,\quad 0\leq n_k\leq |\L_k| \text{for all} 1\leq k\leq d
\cr
}
\Eq(P.p2.31)
$$
Indeed, to each $z=(z^1,\dots,z^d)\in\G_{N,d}$ corresponds the
$d$-tuple of integers $(n_1,\dots,n_d)$ where, by \eqv(2.1.10),
$n_k\equiv\frac{|\L_k|}{2}(1-z^k)=\frac{|\L_k|}{2}(x^k-z^k)$
is the distance from $x$ to $z$ along the $k$-th coordinate axis,
and, in view \eqv(N.0), $n_1+\dots+n_d=n$.
Thus, a path going from $x$ to $z$ in exactly $n$
steps is a path that takes $n_k$ successive ``downwards'' steps along the
$k$-th coordinate axis, for each $1\leq k\leq d$.
Now the number of such paths simply is the multinomial number
$$
\frac{n!}{n_1!\dots n_d!}
\Eq(P.p2.28)
$$
(see {\it e\.g\.} [Be]) and, by \eqv(2.1.14) of Lemma \thv(L.3), all these paths have
same probability, given by 
$$
\eqalign{
&
\prod_{l_1=1}^{n_1}\frac{|\L_1|}{N}\left(1-\frac{l_1}{|\L_1|}\right)
\dots
\prod_{l_d=1}^{n_d}\frac{|\L_d|}{N}\left(1-\frac{l_d}{|\L_d|}\right)
\cr
=&
\frac{1}{N^n}\frac{|\L_1|!}{(|\L_1|-n_1)!}\dots\frac{|\L_d|!}{(|\L_d|-n_d)!}
}
\Eq(P.p2.29)
$$
Therefore,
$$
r^{(n)}_N(x,z)
=\frac{n!}{N^n}
\binom{|\L_1|}{n_1}\dots\binom{|\L_d|}{n_d}
=\frac{n!}{N^n}\frac{\Q_N(z)}{\Q_N(x)}
\Eq(P.p2.30)
$$
where the last equality follows from \eqv(2.1.110) and \eqv(2.1.13).
This proves \eqv(P.p2.26). To prove \eqv(P.p2.27) we likewise begin by
counting the number of paths of length $n+2$ going from $x$ to a point
$z\in\del_mx$, given $m\in I(n)$.
Since choosing a point $z\in\del_mx$ is equivalent to
choosing an element $(m_1,\dots,m_d)\in\QQ^{\G}_d(m)$,
each path going from  $x$ to $z$ must at least contain the steps of a path
going $x$ to $z$ in exactly $m$ steps (namely, for each $1\leq k\leq d$,
$m_k$ successive ``downwards'' steps along the $k$-th coordinate axis).
Denoting by $\QQ_d(p)$ be the set of integer
solutions of the unconstrained equation
$$
p_1+\dots+p_d=p\,,\quad p_i\geq 0 
\Eq(P.p2.31bis)
$$
the remaining $n+2-m=2p$ extra steps can be distributed over the different axis
according to any element $(p_1,\dots,p_d)\in\QQ_d(p)$. More precisely, for each
such element, the total number of steps taken along the $k$-th coordinate axis
is $m_k+2p_k$,
of which the chains traces back exactly $p_k$\note{These steps of course need not
be distinct (the chain can trace back $l_k\leq p_k$ times a same step) and, clearly,
they need not either be consecutive.}. Therefore, for fixed
$(p_1,\dots,p_d)\in\QQ_d(p)$, the number of paths going from $x$ to
$z$ in $\sum_{k=1}^d(m_k+2p_k)=m+2p$ steps is bounded above by
$$
\frac{(m+2p)!}{(m_1+2p_1)!\dots(m_d+2p_d)!}
\binom{m_1+2p_1}{p_1}\dots\binom{m_d+2p_d}{p_d}
\Eq(P.p2.32)
$$
and their probability is of the form
$A_{(m_1,\dots,m_d)}B_{(p_1,\dots,p_d)}$ where
$$
\eqalign{
A_{(m_1,\dots,m_d)}
&=
\prod_{l_1=1}^{m_1}\frac{|\L_1|}{N}\left(1-\frac{l_1}{|\L_1|}\right)
\dots
\prod_{l_d=1}^{m_d}\frac{|\L_d|}{N}\left(1-\frac{l_d}{|\L_d|}\right)
\cr
B_{(p_1,\dots,p_d)}
&=\prod_{j_1=1}^{p_1}\frac{|\L_1|}{N}\left(1-\frac{a_{j_1}}{|\L_1|}\right)
\frac{|\L_1|}{N}\left(\frac{a_{j_1}}{|\L_1|}\right)
\dots
\prod_{j_d=1}^{p_d}\frac{|\L_d|}{N}\left(1-\frac{a_{j_d}}{|\L_d|}\right)
\frac{|\L_d|}{N}\left(\frac{a_{j_d}}{|\L_d|}\right)
\cr
}
\Eq(P.p2.33)
$$
and where, for each $1\leq k\leq d$,  $\left(a_{j_k}\right)_{1\leq j_k\leq p_k}$
is a family of integers having the following properties:
$$
1\leq a_{j_k}\leq m_k+p_k\,,
\Eq(P.p2.34)
$$
and
$$
\eqalign{
&\text{at most one of the} {a_{j_k}}'s \text{takes the value} m_k+p_k,\cr
&\text{at most two of the} {a_{j_k}}'s \text{take the value} m_k+p_k-1,\cr
&\dots\cr
&\text{at most} p_k \text{of the} {a_{j_k}}'s \text{take the value} m_k+1.\cr
}
\Eq(P.p2.35)
$$
To reason this out simply note that...
Thus
$$
\prod_{j_k=1}^{p_k}\frac{|\L_k|}{N}\left(1-\frac{a_{j_k}}{|\L_k|}\right)
\frac{|\L_k|}{N}\left(\frac{a_{j_k}}{|\L_k|}\right)
\leq
\prod_{j_k=1}^{p_k}\frac{|\L_k|}{N}
\frac{a_{j_k}}{N}
\leq \left(\frac{|\L_k|}{N}\right)^{p_k}\frac{1}{N^{p_k}}\frac{(m_k+p_k)!}{m_k!}
\Eq(P.p2.36)
$$
Inserting this bound in \eqv(P.p2.33) and making use of \eqv(P.p2.30) we get:
$$
\eqalign{
A_{(m_1,\dots,m_d)}B_{(p_1,\dots,p_d)}
&\leq
\sum_{(p_1,\dots,p_d)\in\QQ_d(p)}
r^{(n)}_N(x,z)
\frac{1}{p_1!\dots p_d!}
\left(\frac{|\L_1|}{N}\right)^{p_1}\dots\left(\frac{|\L_d|}{N}\right)^{p_d}
\frac{(m+2p)!}{N^{p} m!}
\cr
&=r^{(n)}_N(x,z)\frac{1}{N^p}\frac{(m+2p)!}{p!\, m!}
\left(\frac{|\L_1|}{m_1}+\dots+\frac{|\L_d|}{m_d}\right)^{p_1+\dots+p_d}
\cr
&=r^{(n)}_N(x,z)\frac{1}{N^p}\frac{(m+2p)!}{p!\, m!}
\cr
}
\Eq(P.p2.37)
$$
which proves \eqv(P.p2.27).\endproof

\lemma{\TH(P.lemma.11)}{\it For $F_i$ defined in \eqv(P.p2.21), with the
notation of Lemma \thv(P.lemma.10).
$$
\eqalign{
F_1&\leq\kappa(n)\frac{n!}{N^n}
\cr
%
F_2&\leq \kappa^2(n+2) \frac{(n+2)!}{N^{n+2}}
\sum_{m\in I(n)}\frac{N^p}{p!} |\del_mx|
\cr
}
\Eq(P.p2.38)
$$
}

\proof
Inserting \eqv(P.p2.26) in \eqv(P.p2.23) and, combining the result
with \eqv(P.p2.21) for $i=1$ and with \eqv(P.p2.11), proves the bound
\eqv(P.p2.38) on $F_1$.
We next bound $F_2$. Assuming first that $m\neq n$, and combining
the results of lemma \thv(P.lemma.10) with \eqv(P.p2.24), we get
$$
\eqalign{
\P^{\circ}\left(\left\{n+2\leq\t^x_y<\t^x_{x\cup 0}\right\}\cap\AA_m\right)
\leq &
\sum_{z\in\del_mx} \frac{m!}{N^m}\frac{1}{N^p}\frac{(m+2p)!}{m!\,p!}
\frac{\Q_N(z)}{\Q_N(x)}\P^{\circ}\left(\t^z_y<\t^z_{x\cup 0}\right)
\cr
=&\sum_{z\in\del_mx} \frac{(n+2)!}{N^{n+2}}\frac{N^p}{p!}
\frac{\Q_N(z)}{\Q_N(x)}\P^{\circ}\left(\t^z_y<\t^z_{x\cup 0}\right)
}
\Eq(P.p2.39)
$$
Observing that
$
\P^{\circ}\left(\t^z_y<\t^z_{x\cup 0}\right)
\leq
\P^{\circ}\left(\t^z_y<\t^z_{0}\right),
$
it follows from \eqv(P.p2.12bis) that, for $z\in\del_mx$,
$$
\P^{\circ}\left(\t^z_y<\t^z_{x\cup 0}\right)
\leq
\kappa(\dist(z,\SS_d))\frac{\Q_N(y)}{\Q_N(z)}
\leq
\kappa(m)\frac{\Q_N(y)}{\Q_N(z)}
\Eq(P.p2.40)
$$
Therefore
$$
\eqalign{
\P^{\circ}\left(\left\{n+2\leq\t^x_y<\t^x_{x\cup 0}\right\}\cap\AA_m\right)
\leq &
\kappa(m)\frac{\Q_N(y)}{\Q_N(x)}\frac{(n+2)!}{N^{n+2}}\frac{N^p}{p!} |\del_mx|
}
\Eq(P.p2.41)
$$
In the same way it follows from \eqv(P.p2.25) that, for $m=n$,
$$
\eqalign{
\P^{\circ}\left(\left\{n+2\leq\t^x_y<\t^x_{x\cup 0}\right\}\cap\AA_n\right)
\leq &
\frac{\Q_N(y)}{\Q_N(x)}\frac{(n+2)!}{N^{n+2}}\frac{N^p}{p!}\Bigl[
1+\kappa(n)|(\del_nx)\setminus y|
\Bigr]
\cr
\leq &
\kappa(n)\frac{\Q_N(y)}{\Q_N(x)}\frac{(n+2)!}{N^{n+2}}\frac{N^p}{p!} |\del_nx|
}
\Eq(P.p2.42)
$$
where we used that $\kappa(n)\geq 1$ for all $n$. Finally, by \eqv(P.p2.21),
inserting \eqv(P.p2.41) and
\eqv(P.p2.42) in \eqv(P.p2.16) yields
$$
F_2
\leq
\kappa(n)\sum_{m\in I(n)}
\kappa(m)\frac{(n+2)!}{N^{n+2}}\frac{N^p}{p!} |\del_mx|
\Eq(P.p2.43)
$$
This and the fact that
$\kappa(m)\leq\kappa(n)\leq\kappa(n+2)$ for all $0<m\leq n+2$
proves the bound on $F_2$ of \eqv(P.p2.38).
\endproof

The proof of Theorem  \thv(P.prop.2stronger) is now complete.\endproof

\proofof{Theorem  \thv(P.prop.2)}
Since all partitions into $d\leq\frac{\a}{10}\frac{N}{\log N}$ classes are
log-regular with rate $\a$, Theorem  \thv(P.prop.2) follows from
Theorem  \thv(P.prop.2stronger) by choosing $d\leq\a_0\frac{N}{\log N}$
provided that $\a_0\leq \frac{\a}{10}$ for some $0\leq\a\leq 1/2$, i\.e\.
provided that $\a_0\leq 1/20$.
\endproof

We are now ready to prove the lower bound on $\P^{\circ}(\t^x_0<\t^x_x)$
of Theorem  \thv(P.prop.1) (i\.e\. the upper bound of Theorem \thv(P.prop.1)).

\bigskip
\line{\bf 3.4. Proof of the upper bound of Theorem \thv(P.prop.1).\hfill}

The upper bound of Theorem \thv(P.prop.1) will easily be deduced from the following lemma.

\lemma{\TH(P.lemma.8)}{\it
Assume that the $d$-partition $\L$ is log-regular.
There exists a constant $0<\a'<\infty$ such that if
$d\leq \a'\frac{N}{\log N}$ then, for  all $x\in\SS_d$,
%
%
$$
\P^{\circ}(\t^x_0<\t^x_x)\geq 1-\frac{1}{N}-\frac{3}{N^2}(1+O(1/{\sqrt N}))
\Eq(P.l8.1)
$$
}

\proof
Let $\O_x$ be the set of
paths on $\G_{N,d}$ that start in $x$ and return to $x$ without visiting 0:
$$
\O_x\equiv\bigcup_{L\geq 2}
\bigl\{\o=(\o_0,\o_1,\dots,\o_L)\:
\o_0=\o_L=x,\,\forall_{0<k<L}\,\o_k\neq 0 \bigr\}
\Eq(P.l8.2)
$$
We want to classify these paths according to their canonical projection
on the coordinate axis. For simplicity, we will place the origin
of the coordinate system at $x$  and, as usual, we set $x=(1,\dots,1)$.
Thus, for each $1\leq k\leq d$, recalling the notation $y=(y^1,\dots,y^d)$, we let
$\pi_k$ be the map
$y\mapsto\pi_k y =((\pi_k y)^1,\dots,(\pi_k y)^d)$,
where
$$
(\pi_k y)^k=y^k\,,\text{and} (\pi_k y)^{k'}=1 \text{for all} k'\neq k
$$
Given a path $\o\in\O_x$ we then call the set $\pi\o$ defined by
$$
\pi\o\equiv\bigcup_{1\leq k\leq d}\pi_k\o\,,\quad
\pi_k\o\equiv\{\pi_k\o_0,\pi_k\o_1,\dots,\pi_k\o_L\}
\Eq(P.l8.3)
$$
the projection of this path.
Now observe that the set of projections of all the
paths of $\O_x$ simply is the set
$$
\{\o\in\O_x \: \pi\o\}
=\bigcup_{1\leq m\leq N}\bigcup_{(m_1,\dots,m_d)\in\QQ^{\G}_d(m)}\EE_m(m_1,\dots,m_d)
\Eq(P.l8.4)
$$
where, given an integer $1\leq m\leq N$ and an element
$(m_1,\dots,m_d)\in\QQ^{\G}_d(m)$ (see \eqv(P.p2.31)),
$$
\eqalign{
&\EE_m(m_1,\dots,m_d)\equiv\bigcup_{1\leq k\leq d}\EE_m(m_k),\,
\cr
&\EE_m(m_k)\equiv\{x, x-\sfrac{2}{|\L_k(I)|}u_k,\dots, x-\sfrac{2m_k}{|\L_k(I)|}u_k\}
}
\Eq(P.l8.5)
$$
With this notation in hand we may rewrite the quantity
$1-\P^{\circ}(\t^x_0<\t^x_x)=\P^{\circ}(\t^x_x<\t^x_0)$ as
$\P^{\circ}(\t^x_x<\t^x_0)=\P^{\circ}(\O_x)$
which, for a given fixed $M$ (to be chosen later as a function of $N,d$),
we may decompose according to the cardinality of the set $\pi\o\setminus x$ into
three terms,
$$
\P^{\circ}(\O_x)=\RR_1+\RR_2+\RR_3
\Eq(P.l8.11)
$$
where
$$
\eqalign{
\RR_1
&=
\P^{\circ}(\o\in\O_x\: |\pi\o\setminus x|=1)
\cr
\RR_2
&=
\sum_{m=2}^{M}\P^{\circ}(\o\in\O_x\: |\pi\o\setminus x|=m)
\cr
\RR_3
&=
\P^{\circ}(\o\in\O_x\: |\pi\o\setminus x|>M)
\cr
}
\Eq(P.l8.12)
$$
We will now estimate the three probabilities of \eqv(P.l8.12) separately.
Firstly, note that the set $\{\o\in\O_x\: |\pi\o\setminus x|=1\}$ can only contain
paths of length $|\o|=2$. Thus
$$
\eqalign{
\RR_1
&=
\P^{\circ}\left(2=\t^x_x<\t^x_0\right)
\cr
&=
\sum_{l=1}^d
r_N(x, x-\sfrac{2u_l}{|\L_l|})r_N(x-\sfrac{2u_l}{|\L_l|}, x)
\cr
&=
\sum_{l=1}^d\frac{|\L_l|}{N}\frac{1}{N}
\cr
&=
\frac{1}{N}
\cr
}
\Eq(P.l8.13)
$$
Let us turn to $\RR_2$. Our strategy here is to enumerate the paths
of the set $\{\o\in\O_x\: |\pi\o\setminus x|=m\}$ and to bound the probability of
each path. To do so we associate to each $\EE_m(m_1,\dots,m_d)$ the set of edges:
$$
\eqalign{
&E(\EE_m(m_1,\dots,m_d))\equiv\bigcup_{1\leq k\leq d}E(\EE_m(m_k)),\,
\cr
& E(\EE_m(m_k)) \equiv
\left\{(x,x')\in \EE_m(m_k)\mid\exists_{s\in\{-1,1\}}\,
:\,x'-x=s\sfrac{2}{|\L_k(I)|}u_k\right\}
}
\Eq(P.l8.7)
$$
Next, choose $\o\in\{\o\in\O_x\: |\pi\o\setminus x|=m\}$ and, for
each $1\leq k\leq d$, let us denote by $l_k$ be the number of steps of
$\o$ that project onto the $k^{th}$ axis: namely, if $|\o|=L$,
$$
l_k(\o)=\sum_{0\leq l< L}\1_{\{(\pi_k\o_l, \pi_k\o_{l+1})\in\EE_m(m_k)\}}
\Eq(P.l8.8)
$$
A step along the $k^{th}$ axis that decreases (increases) the value of the
$k^{th}$ coordinate will be called a downward (upward) step. Clearly, because
a path $\o\in\O_x$ ends where it begins, each (non oriented) edge $(\o_l,\o_{l-1})$
of $\o$ must be stepped through an equal number of times upward and downward.
As a result such paths must be of even length.
Therefore, setting
$$
L=2n \text{and} l_k=2n_k \text{for all} 1\leq k\leq d
\Eq(P.l8.9)
$$
we have
$$
\eqalign{
& n_1+\dots+n_d=n\,,\quad  n_k\geq m_k \text{for all} 1\leq k\leq d
\cr
}
\Eq(P.l8.10)
$$
Then
$$
\eqalign{
&\P^{\circ}(\o\in\O_x\: |\pi\o\setminus x|=m)
\cr
=&\sum_{L=2m}^{\infty}\P^{\circ}(\o\in\O_x\: |\o|=L\,, |\pi\o\setminus x|=m)
\cr
=&
\sum_{n=m}^{\infty}
\sum_{(m_1,\dots,m_d)\in\QQ^{\G}_d(m)}
\sum_{
{n_1\geq m_1,\dots, n_d\geq m_d\:}\atop{n_1+\dots+n_d=n}
}
\P^{\circ}(\o\in\O_x\: \forall_{k=1}^d\,\pi_k\o\in\EE_m(m_k)\,,l_k(\o)=2n_k)
}
\Eq(P.l8.14)
$$
The probabilities appearing in the last line above are easily bounded.
On the one hand the number of paths in
$\{\o\in\O_x\: \forall_{k=1}^d\,\pi_k\o\in\EE_m(m_k)\,,l_k(\o)=2n_k\}$
is bounded above by the number of ways to arrange the
$|\o|=2n$ steps of the path into sequences of
$n_k$ upward steps and $n_k$ downward steps along each of the $1\leq k\leq d$ axis,
disregarding all constraints. This yields:
$$
\frac{(2n)!}{(n_1!)^2\dots (n_d!)^2}
\Eq(P.l8.15)
$$
On the other hand, as used already in the proof of Lemma \thv(P.lemma.10),
the probability to step up and down a given edge $(\o_l,\o_{l+1})$ along,
say, the $k^{th}$ coordinate axis, only depends on its projection on this axis
(see e\.g\. \eqv(P.p2.33)) and the probability of a path in
$\{\o\in\O_x\: \forall_{k=1}^d\,\pi_k\o\in\EE_m(m_k)\,,l_k(\o)=2n_k\}$
is easily seen to be of the form
$$
\prod_{k=1}^d\prod_{j_k=1}^{m_k}\left\{\frac{|\L_k|}{N}\left(1-\frac{j_k}{|\L_k|}\right)
\frac{|\L_k|}{N}\left(\frac{j_k}{|\L_k|}\right)\right\}^{a^k_{j_k}}
\Eq(P.l8.17)
$$
where, for  $\QQ_d(\,.\,)$ defined in \eqv(P.p2.31bis),
$(a^k_1,\dots,a^k_{m_k})$ is an element of $\QQ_d(m_k)$ having the property that
$a^k_{j_k}\geq 1$ for all $1\leq j_k\leq m_k$.
The quantity \eqv(P.l8.17) may thus be bounded above by
$$
\eqalign{
&\prod_{k=1}^d\prod_{j_k=1}^{m_k}
\left\{
\frac{|\L_k|}{N}\left(1-\frac{j_k}{|\L_k|}\right)
\frac{|\L_k|}{N}\left(\frac{j_k}{|\L_k|}\right)
\left(\frac{|\L_k|}{N}\frac{m_k}{N}\right)^{a^k_{j_k}-1}
\right\}
\cr
= &
\frac{1}{N^m}\prod_{k=1}^d
\frac{1}{N^{m_k}}\frac{m_k!|\L_k|!}{(|\L_k|-m_k)!}
\left(\frac{|\L_k|}{N}\frac{m_k}{N}\right)^{n_k-m_k}
}
\Eq(P.l8.18)
$$
Inserting \eqv(P.l8.15) and \eqv(P.l8.18) in \eqv(P.l8.14)
 we have to evaluate the resulting sum. To do so note first that
$$
\eqalign{
\frac{(2n)!}{\prod_{k=1}^d(n_k!)^2}
=&
\frac{1}{\prod_{k=1}^d(m_k!)^2}
\frac{(2n)!}{(2(n-m))!}
\frac{(2(n-m))!}{\prod_{k=1}^d(2(n_k-m_k))!}
\prod_{k=1}^d\frac{(2(n_k-m_k))!}{((n_k-m_k)!)^2}
\binom{n_k}{m_k}^{-2}
\cr
\leq&
\frac{2^{2(n-m)}}{\prod_{k=1}^d(m_k!)^2}
\frac{(2n)!}{(2(n-m))!}
\frac{(2(n-m))!}{\prod_{k=1}^d(2(n_k-m_k))!}
}
\Eq(P.l8.19)
$$
Thus
$$
\eqalign{
&\sum_{{n_1\geq m_1,\dots, n_d\geq m_d\:}\atop{n_1+\dots+n_d=n}}
\frac{(2n)!}{(n_1!)^2\dots (n_d!)^2}
\prod_{k=1}^d
\left(\frac{|\L_k|}{N}\frac{m_k}{N}\right)^{n_k-m_k}
\cr
&\leq
\frac{2^{2(n-m)}}{\prod_{k=1}^d(m_k!)^2}
\frac{(2n)!}{(2(n-m))!}
\left(
\sum_{k=1}^d\sqrt{\frac{|\L_k|}{N}\frac{m_k}{N}}
\right)^{2(n-m)}
\cr
&\leq
\frac{2^{2(n-m)}}{\prod_{k=1}^d(m_k!)^2}
\frac{(2n)!}{(2(n-m))!}
\left(\frac{m}{N}\right)^{n-m}
}
\Eq(P.l8.20)
$$
where the last line follows from Schwarz's inequality.
From now on we assume that there exists a constant $0\leq C<1$
such that $4\frac{m}{N}\leq C$ for all $m\leq M$.
Using that
$
\frac{(2n)!}{(2(n-m))!}x^{2(n-m)}
=\frac{d^{2(n-m)+1}}{dx^{2(n-m)+1}}x^{2n}
$
we get,
$$
\sum_{n=m}^{\infty}
\frac{(2n)!}{(2(n-m))!}
\left(\frac{m}{N}\right)^{n-m}
=\frac{(2m-1)!}{\left(1-\sqrt{4\frac{m}{N}}\right)^{2m}}
\Eq(P.l8.21)
$$
Finally,
$$
\eqalign{
&\P^{\circ}(\o\in\O_x\: |\pi\o\setminus x|=m)
\cr
\leq
&
\frac{1}{N^m}
\sum_{(m_1,\dots,m_d)\in\QQ^{\G}_d(m)}
\frac{(2m-1)!}{\left(1-\sqrt{4\frac{m}{N}}\right)^{2m}}
\frac{1}{\prod_{k=1}^d(m_k!)^2}
\prod_{k=1}^d
\frac{1}{N^{m_k}}\frac{m_k!|\L_k|!}{(|\L_k|-m_k)!}
\cr
}
\Eq(P.l8.22)
$$
and since
$$
\eqalign{
&\sum_{(m_1,\dots,m_d)\in\QQ^{\G}_d(m)}\binom{m}{m_1,\dots,m_d}
\prod_{k=1}^d
\frac{1}{N^{m_k}}\frac{|\L_k|!}{(|\L_k|-m_k)!}
\cr
&\leq
\sum_{(m_1,\dots,m_d)\in\QQ^{\G}_d(m)}\binom{m}{m_1,\dots,m_d}
\prod_{k=1}^d
\left(\frac{|\L_k|}{N}\right)^{m_k}
\cr
&\leq
\left(\sum_{k=1}^d\frac{|\L_k|}{N}\right)^m
\cr
&=1
}
\Eq(P.l8.23)
$$
we obtain
$$
\P^{\circ}(\o\in\O_x\: |\pi\o\setminus x|=m)
\leq
\frac{a(m)}{N^m}\,,\quad
a(m)\equiv
\frac{(2m-1)!}{m!}
\left(1-\sqrt{4\sfrac{m}{N}}\right)^{-2m}
\Eq(P.l8.24)
$$
To bound the term $\RR_2$ from \eqv(P.l8.12) we still have to sum
\eqv(P.l8.24) over $2\leq m\leq M$. Writing
$$
\RR_2=\frac{a(2)}{N^2}+\frac{a(3)}{N^3}\left[1+\sum_{m=4}^M\frac{a(m)/a(3)}{N^{m-3}}\right]
\Eq(P.l8.25)
$$
we have, using Stirling's formula that, for some constant $c>0$,
$$
1+\sum_{m=4}^M\frac{a(m)/a(3)}{N^{m-3}}
\leq 1+\sum_{m=4}^Me^{-c(m-3)}\left(\frac{m}{N}\right)^{m-3}
\leq 1+\sum_{m=4}^M\left(\frac{M}{N}\right)^{m-3}
\leq \left(1-\frac{M}{N}\right)^{-1}
\Eq(P.l8.26)
$$
and, since
$
{a(2)}/{N^2}=({3}/{N^2})(1+O(1/{\sqrt N})),
$
we get, for all $M$ such that $4\frac{M}{N}\leq C<1$,
$$
\RR_2\leq\frac{3}{N^2}\left(1+O(1/{\sqrt N})\right)
\Eq(P.l8.28)
$$

It now remains to bound $\RR_3$. Observe that all paths
in $\{\o\in\O_x\: |\pi\o\setminus x|>M\}$ must hit the set
$\MM\equiv\del_{\lfloor\frac{M}{d}\rfloor}x$ (here $\lfloor\frac{M}{d}\rfloor$
denotes the integer part of $\frac{M}{d}$). Assume indeed that
it is not the case. Since $\pi\o\in\EE_m(m_1,\dots,m_d)$ for some
$m>M$ and $(m_1,\dots,m_d)\in\QQ^{\G}_d(m)$, this would in particular
imply that $\max_{1\leq k\leq d}m_k<\lfloor\frac{M}{d}\rfloor$.
But this in turn implies that $m=m_1+\dots+m_d<M$, which is a contradiction.
We are thus lead to write
$$
\eqalign{
\RR_3
\leq &
\P^{\circ}(\t^x_{\MM}<\t^x_x<\t^x_0)
\cr
=&
\sum_{z\in\MM}
\P^{\circ}\left(\t^x_z<\t^x_{x\cup 0\cup\MM\setminus z}\right)
\P^{\circ}\left(\t^z_x<\t^z_0\right)
\cr
\leq &
\max_{z\in\MM}\P^{\circ}\left(\t^z_x<\t^z_0\right)
\sum_{z\in\MM}
\P^{\circ}\left(\t^x_z<\t^x_{x\cup 0\cup\MM\setminus z}\right)
\cr
\leq &
\max_{z\in\MM}\P^{\circ}\left(\t^z_x<\t^z_0\right)
\cr
}
\Eq(P.l8.29)
$$
and, under the assumption that
the partition $\L$ is log-regular,
by Theorem  \thv(P.prop.2stronger),
$$
\RR_3\leq \max_{z\in\MM}F(\dist(x,z))
\leq F(\lfloor M/d\rfloor)
\Eq(P.l8.30)
$$

Plugging \eqv(P.l8.13), \eqv(P.l8.28), and \eqv(P.l8.30) in
\eqv(P.l8.11) we have proven the following statement:

\lemma{\TH(P.lemma.12)}{\it
Assume that the $d$-partition $\L$ is log-regular.
Then,
for all $x\in\SS_d$, for all $M$ such that $4\frac{M}{N}\leq C$ where
$0\leq C<1$ is a numerical constant, and for large enough $N$,
$$
\P^{\circ}(\t^x_x<\t^x_0)
=\frac{1}{N}
+\frac{3}{N^2}(1+O(1/{\sqrt N}))
+F(\lfloor M/d\rfloor)
\Eq(P.l8.31)
$$
}

It is easy to check that there exists a constant $0<\a'<1$
such that for all $d\leq \a'\frac{N}{\log N}$,
choosing $M=\frac{C}{4} N$ (that is $\frac{M}{d}\geq \frac{C}{4\a'}\log N$) the bounds
of Lemma \thv(A3.lemma.1) imply
that $F(\lfloor M/d\rfloor)=O(1/N^{\frac{5}{2}})$.
This concludes the proof of Lemma \thv(P.lemma.8).\endproof

As in the proof of Theorem \thv(P.prop.2)
observing that all
partitions into $d\leq\frac{\a}{10}\frac{N}{\log N}$ classes are
log-regular with rate $\a$,
the upper bound of Theorem \thv(P.prop.1)
follows from  Lemma \thv(P.lemma.8)
by choosing $d\leq\a_0\frac{N}{\log N}$ with e\.g\.
$\a_0\leq\inf\left\{\a',\frac{1}{20}\right\}$.


The proof of  Theorem \thv(P.prop.1) is now complete.


\bigskip
\vfill\eject

\bigskip
\chap{4. Estimates on hitting probabilities for the lumped chain.}4

In this section we pursue the investigation of the lumped chains
initiated in Chapter 3. Using the probability estimates established therein we will
prove sharp estimates on the harmonic measure and on `no return before hitting'
probabilities.
As a warm up to these proofs we begin in Section 4.1 by drawing
some simple consequences of Theorem \thv(P.prop.1) and Theorem \thv(P.prop.2)
(Corollary \thv(P.cor.1)).
Doing so, we will show how the bound \eqv(P.p2.1)
of Theorem \thv(P.prop.2) gives rise to the sparsity condition.
The procedure  described in Section 1 will be used repeatedly in the follow-up sections to prove
estimates on: the Harmonic measure starting from zero (Section 4.2);
the Harmonic measure with arbitrary starting point (Section 4.3);  no `return
before hitting probabilities' of the general form
$\P^{\circ}(\t^x_{J\setminus x}<\t^x_x)$ for $J\subset\SS_d$ and $x\in \G_{N,d}$ (Section 4.4).

Let us finally mention that while the results on the harmonic measure of Section 4.2 and 4.3
will be needed both in section 6 and 7,
Corollary \thv(P.cor.1) of Section 4.1 will be crucial for the investigation of the
Laplace transforms carried out in Chapter 7 and,
as mentioned earlier, Theorem \thv(P.theo.0) of Section 4.4 will play a key role in Chapter 6
for the analysis
of hitting times.


\bigskip
\line{\bf 4.1. Generalization of Theorems \thv(P.prop.1) and \thv(P.prop.2):
emergence of the sparsity  condition\hfill}

We begin with some notation and definitions.
Recall from \eqv(N.1) that, given two points $x,y\in \G_{N,d}$, $\del_n x$ denotes the sphere centered
at $x$ and of radius $n$ for the graph distance.
For $x\in\SS_d$ and arbitrary $y\in \G_{N,d}$ define
$$
\phi_x(n)=\max_{y\in\del_n x}
\P^{\circ}\left(\t^y_x<\t^y_0\right)
\Eq(N.2)
$$

\lemma{\TH(N.lemma.1)}{\it For all $x\in\G_{N,d}$, $\phi_x$ is non increasing.}

\proof Let $x\in\SS_d$ be fixed. For all $n\geq 1$ and $y\in\del_{n+1}x$,
$$
\eqalign{
\P^{\circ}\left(\t^y_x<\t^y_0\right)
\leq &
\P^{\circ}(\t^y_{\del_n x}<\t^y_x<\t^y_0)
\cr
=&
\sum_{z\in\del_n x}
\P^{\circ}\left(\t^y_z<\t^y_{x\cup 0\cup(\del_n x)\setminus z}\right)
\P^{\circ}\left(\t^z_x<\t^z_0\right)
\cr
\leq &
\max_{z\in\del_n x}\P^{\circ}\left(\t^z_x<\t^z_0\right)
\sum_{z\in\del_n x}
\P^{\circ}\left(\t^y_z<\t^y_{x\cup 0\cup(\del_n x)\setminus z}\right)
\cr
\leq &
\max_{z\in\del_n x}\P^{\circ}\left(\t^z_x<\t^z_0\right)\P^{\circ}(\t^y_{\del_n x}<\t^y_{x\cup 0})
\cr
\leq &\phi_x(n)
}
\Eq(N.l1.1)
$$
Thus, taking the maximum over $y\in\del_{n+1}x$,
$$
\phi_x(n+1)\leq\phi_x(n)\,,\quad n\geq 1
\Eq(N.l1.2)
$$
which proves the claim of the lemma.\endproof


From Theorem \thv(P.prop.2stronger)
we immediately deduce that:

\lemma{\TH(N.lemma.2)}{\it
Assume that the $d$-partition $\L$ is  log-regular.
Then, for all $x\in\SS_d$,
$$
\phi_x(n)\leq F(n)\,,\quad n\in\N
\Eq(N.3)
$$
}
\proof Just note that
$
\displaystyle
\phi_x(n)\leq \max_{y\in\del_n x} F(\dist(x,y))=F(n)\,.
$
\endproof

Now let $J\subset\G_{N,d}$ and $y\in\G_{N,d}$, and define
$$
\eqalign{
V^{\circ}(y,J)&
=\cases\sum_{z\in J\setminus y}\phi_z(\dist(y,z)), &\text{if} J\setminus y\neq\emptyset\cr
0,&\text{otherwise}
\endcases
\cr
U^{\circ}(y,J)&=\cases \sum_{z\in J\setminus y}F(\dist(y,z)), &\text{if} J\setminus y\neq\emptyset\cr
0,&\text{otherwise}
\endcases
\cr
}
\Eq(N.3bis)
$$

Clearly, by Lemma \thv(N.lemma.2), under the assumptions therein,
$$
V^{\circ}(y,J)\leq U^{\circ}(y,J)\,,\quad J\subset\SS_d\,, y\in\G_{N,d}
\Eq(N.4)
$$
implying that
$$
\sum_{z\in J\setminus y}\P^{\circ}\left(\t^y_z<\t^y_0\right)\leq V^{\circ}(y,J)\leq U^{\circ}(y,J)\leq \max_{y\in J} U^{\circ}(y,J)
\Eq(N.4')
$$
Obviously the function $\max_{y\in J} U^{\circ}(y,J)$, $J\subset\G_{N,d}$, strongly resembles
the function $\UU_{N,d}(A)$, $A\subset\SS_N$, introduced in \eqv(I.12) to define
the notion of sparseness of a set. We will see in Chapter 5 that, on appropriate sets,
these two functions indeed coincide. In view of \eqv(N.4') the
sparsity condition will thus serve to
guarantee the smallness of sums of
the form $\sum_{z\in J\setminus y}\P^{\circ}\left(\t^y_z<\t^y_0\right)$.


We now use the previous observations to prove the following generalization  of
Theorem \thv(P.prop.1) and Theorem  \thv(P.prop.2) where the exclusion point
(respectively hitting point) $x$ is replaced by a subset $J\subset \SS_d$.

\corollary{\TH(P.cor.1)}{\it
Let $d\leq d_0(N)$.
Then, for all $J\in\SS_d$
the following holds:
\item{(i)} For all $x\in J$
$$
1-\frac{1}{N}-\frac{c}{N^2}-V^{\circ}(x,J)
\leq
\P^{\circ}(\t^x_0<\t^x_{J})
\leq
1-\frac{1}{N}
\Eq(P.c1.0)
$$
and
$$
1-\frac{1}{N}-\frac{c}{N^2}-
\frac{1}{|J|}\sum_{z\in J}V^{\circ}(z,J)
\leq
\frac{\Q(0)}{\Q(J)}\,\P^{\circ}\left(\t^{0}_{J}<\t^{0}_{0}\right)
\leq
1-\frac{1}{N}
\Eq(P.c1.00)
$$
for some numerical constant $0<c<4$. (Note that $\Q(J)=|J|2^{-N}$)
\item{(ii)} for all $y\notin J$
$$
\P^{\circ}\left(\t^y_J<\t^y_0\right)\leq V^{\circ}(y,J)
\Eq(P.c1.000)
$$

Moreover \eqv(P.c1.0), \eqv(P.c1.00), and \eqv(P.c1.000) remain true with
$V^{\circ}(\,.\,,J)$ replaced by $U^{\circ}(\,.\,,J)$.
}


\proofof{Corollary \thv(P.cor.1)} Note that
$$
\P^{\circ}(\t^x_0<\t^x_J)=
\P^{\circ}(\t^x_0<\t^x_{x})-\P^{\circ}(\t^x_{J\setminus x}<\t^x_0<\t^x_{x})
\Eq(P.c1.1)
$$
where
$$
\eqalign{
\P^{\circ}(\t^x_{J\setminus x}<\t^x_0<\t^x_{x})
=&
\sum_{z\in J\setminus x}\P^{\circ}(\t^x_{z}<\t^x_{(J\setminus z)\cup 0})
             \P^{\circ}(\t^z_0<\t^z_x)
\cr
\leq &\sum_{z\in J\setminus x}\P^{\circ}(\t^x_{z}<\t^x_{(J\setminus z)\cup 0})\cr
\leq &\sum_{z\in J\setminus x}\P^{\circ}(\t^x_{z}<\t^x_0)\cr
\leq &V^{\circ}(x,J)\cr
}
\Eq(P.c1.2)
$$
Thus
$$
\P^{\circ}(\t^x_0<\t^x_{x})-V^{\circ}(x,J)
\leq
\P^{\circ}(\t^x_0<\t^x_{J\cup x})
\leq
\P^{\circ}(\t^x_0<\t^x_{x})
\Eq(P.c1.4)
$$
and, together with Theorem  \thv(P.prop.1), this proves \eqv(P.c1.0).
Next,
$$
\eqalign{
\P^{\circ}\left(\t^{0}_{J}<\t^{0}_{0}\right)
&=\sum_{z\in J}
\P^{\circ}\left(\t^{0}_{z}<\t^{0}_{0\cup(J\setminus z)}\right)\cr
&=\sum_{z\in J}\frac{\Q(z)}{\Q(0)}
\P^{\circ}\left(\t^{z}_{0}<\t^{z}_{J}\right)\cr
&=\frac{1}{2^N\Q(0)}
\sum_{z\in J}\P^{\circ}\left(\t^{z}_{0}<\t^{z}_{J}\right)\cr
}
\Eq(TL.40bis)
$$
where, for each $z\in J$, $\P^{\circ}\left(\t^{z}_{0}<\t^{z}_{J}\right)$
obeys the bounds of \eqv(P.c1.0). Since $\Q(J)=|J|2^{-N}$, \eqv(P.c1.00)
is proven. Finally
$$
\P^{\circ}\left(\t^y_J<\t^y_0\right)
=\sum_{z\in J}\P^{\circ}(\t^y_{z}<\t^y_{(J\setminus z)\cup 0})
\leq\sum_{z\in J}\P^{\circ}(\t^y_{z}<\t^y_0)
\leq V^{\circ}(y,J)
\Eq(P.c1.5)
$$
proving \eqv(P.c1.000).
In view of \eqv(N.4),  \eqv(P.c1.0), \eqv(P.c1.00), and \eqv(P.c1.000) remain
true with $V^{\circ}(\,.\,,J)$ replaced by $U^{\circ}(\,.\,,J)$.
Corollary \thv(P.cor.1) is proven.
\endproof



\vfill\eject
\bigskip
\line{\bf 4.2. The harmonic measure starting from the origin.\hfill}



Given $J\subset\G_{N,d}$ and $y\notin J$, let $H^{\circ}_J(y,x)$ denote the
harmonic measure of the lumped chain, namely,
$$
H^{\circ}_J(y,x)
=
\P^{\circ}\left(\t^y_x<\t^{y}_{J\setminus x}\right)\,,\quad x\in J
\Eq(N.10)
$$

\lemma{\TH(P.lemma.9)}{\it
Let $d\leq d_0(N)$.
Then, for all $J\subset\SS_d$  and all $x\in J$,
$$
\frac{c_N^-}{|J|}\Bigl[1-(1+O(\sfrac{1}{N}))V^{\circ}(x,J)\Bigr]
\leq
H^{\circ}_J(0,x)
\leq
\frac{c_N^+}{|J|}
\Bigl[1-\max_{z\in J} V^{\circ}(z,J)\Bigr]^{-1}
\Eq(P.l9.1)
$$
where, for some numerical constant $0<c<4$,
$$
c_N^{\pm}=
1\pm \frac{c}{N^2}
\Eq(P.l9.2)
$$
Moreover \eqv(P.l9.1) remains true with $V^{\circ}(\,.\,,J)$ replaced by $U^{\circ}(\,.\,,J)$.
}



\proof Again using Lemma \thv(P.renew)
$$
\P^{\circ}\left(\t^0_x<\t^{0}_{J\setminus x}\right)
= \frac{\P^{\circ}\left(\t^0_x<\t^{0}_{(J\setminus x)\cup 0}\right)}
{\P^{\circ}\left(\t^0_J<\t^{0}_0\right)}
= \frac{\P^{\circ}\left(\t^0_x<\t^{0}_{(J\setminus x)\cup 0}\right)}
{\sum_{y\in J}
\P^{\circ}\left(\t^0_{y}<\t^{0}_{(J\setminus y)\cup 0}\right)}
\Eq(P.l9.3)
$$
We basically want to show that this last ratio behaves like
$$
R_x\equiv \frac{\P^{\circ}\left(\t^0_x<\t^{0}_{0}\right)}
{\sum_{y\in J}\P^{\circ}\left(\t^0_y<\t^{0}_{0}\right)}
\Eq(P.l9.11)
$$
Let us first treat the denominator of \eqv(P.l9.3). Observe that
$$
\P^{\circ}\left(\t^0_y<\t^{0}_{(J\setminus y)\cup 0}\right)
=
\P^{\circ}\left(\t^0_y<\t^{0}_{0}\right)
-
\P^{\circ}\left(\t^0_{J\setminus y}<\t^{0}_{y}<\t^{0}_{0}\right)
\Eq(P.l9.4)
$$
and that
$$
\P^{\circ}\left(\t^0_{J\setminus y}<\t^{0}_{y}<\t^{0}_{0}\right)
=
\sum_{z\in J\setminus y}
\P^{\circ}\left(\t^0_{z}<\t^{0}_{(J\setminus z)\cup 0}\right)
\P^{\circ}\left(\t^z_{y}<\t^{z}_{0}\right)
\Eq(P.l9.5)
$$
Then, summing \eqv(P.l9.5) over $y\in J$,
$$
\eqalign{
\sum_{y\in J}\P^{\circ}\left(\t^0_{J\setminus y}<\t^{0}_{y}<\t^{0}_{0}\right)
\leq &
\sum_{z\in J}
\P^{\circ}\left(\t^0_{z}<\t^{0}_{(J\setminus z)\cup 0}\right)
\sum_{y\in J\setminus z}
\P^{\circ}\left(\t^z_{y}<\t^{z}_{0}\right)
\cr
\leq & \sum_{z\in J}\P^{\circ}\left(\t^0_{z}<\t^{0}_{(J\setminus z)\cup 0}\right)
V^{\circ}(z,J)
\cr
\leq &\max_{z\in J}V^{\circ}(z,J)\P^{\circ}\left(\t^0_{J}<\t^{0}_{0}\right)
\cr
}
\Eq(P.l9.6)
$$
Combining  \eqv(P.l9.6) with  \eqv(P.l9.4) and using that
$$
\frac{\P^{\circ}\left(\t^0_J<\t^{0}_{0}\right)}
{\sum_{y\in J}\P^{\circ}\left(\t^0_y<\t^{0}_{0}\right)}
\leq 1
\Eq(P.l9.9)
$$
we get the bounds:
$$
\left(1-\max_{z\in J}V^{\circ}(z,J)\right)
\sum_{y\in J}\P^{\circ}\left(\t^0_{y}<\t^{0}_{0}\right)
\leq
\sum_{y\in J}\P^{\circ}\left(\t^0_{y}<\t^{0}_{(J\setminus y)\cup 0}\right)
\leq
\sum_{y\in J}\P^{\circ}\left(\t^0_{y}<\t^{0}_{0}\right)
\Eq(P.l9.4up)
$$

To bound the numerator of \eqv(P.l9.3) from above we of course simply use
\eqv(P.l9.4), removing the negative term. To get a good lower bound  we do not
use \eqv(P.l9.5) directly. Instead, we use that plugging \eqv(P.l9.4) in
the r\.h\.s\. of  \eqv(P.l9.5) gives
$$
\eqalign{
\P^{\circ}\left(\t^0_{y}<\t^{0}_{(J\setminus y)\cup 0}\right)
=&
\P^{\circ}\left(\t^0_y<\t^{0}_{0}\right)-
\sum_{z\in J\setminus y}
\left\{
\P^{\circ}\left(\t^0_z<\t^{0}_{0}\right)-
\P^{\circ}\left(\t^0_{J\setminus z}<\t^{0}_{z}<\t^{0}_{0}\right)
\right\}
\P^{\circ}\left(\t^z_{y}<\t^{z}_{0}\right)
\cr
\geq & \P^{\circ}\left(\t^0_y<\t^{0}_{0}\right)
-\sum_{z\in J\setminus y}
\P^{\circ}\left(\t^0_z<\t^{0}_{0}\right)
\P^{\circ}\left(\t^z_{y}<\t^{z}_{0}\right)
\cr
\geq &
\Bigl\{
R_y-\max_{z\in J}R_z\sum_{z\in J\setminus y}\P^{\circ}\left(\t^z_{y}<\t^{z}_{0}\right)
\Bigr\}
\sum_{y\in J}\P^{\circ}\left(\t^0_{y}<\t^{0}_{0}\right)
}
\Eq(P.l9.4low)
$$
Now by \eqv(P.p2.10), \eqv(P.p2.11), and \eqv(P.p2.12),
since ${\Q_N(z)}={\Q_N(y)}$,
$$
\P^{\circ}\left(\t^z_{y}<\t^{z}_{0}\right)=
\frac{\P^{\circ}\left(\t^y_z<\t^y_{z\cup 0}\right)}
{\P^{\circ}\left(\t^{z}_{y\cup 0}<\t^z_z\right)}
\leq \frac{\P^{\circ}\left(\t^y_z<\t^y_{0}\right)}
{\P^{\circ}\left(\t^{z}_{0}<\t^z_z\right)}
\leq \frac{N}{N-1}\P^{\circ}\left(\t^y_z<\t^y_{0}\right)
\Eq(P.l9.4')
$$
where the last line follows from Theorem  \thv(P.prop.1). Thus
$$
\sum_{y\in J}\P^{\circ}\left(\t^0_{y}<\t^{0}_{0}\right)
\leq \frac{N}{N-1}\sum_{y\in J}\P^{\circ}\left(\t^y_z<\t^y_{0}\right)
\leq \frac{N}{N-1}V^{\circ}(z,J)
\Eq(P.l9.7)
$$
Plugging this back in \eqv(P.l9.4'low), and collecting both upper and lower bounds,
we have established that
$$
R_y-\max_{z\in J}R_z \frac{N}{N-1}V^{\circ}(y,J)
\leq
\frac{
\P^{\circ}\left(\t^0_{y}<\t^{0}_{(J\setminus y)\cup 0}\right)
}
{
\sum_{y\in J}\P^{\circ}\left(\t^0_{y}<\t^{0}_{0}\right)
}
\leq
R_y
\Eq(P.l9.4'low)
$$

Inserting the bounds \eqv(P.l9.4up) and \eqv(P.l9.4'low) into \eqv(P.l9.3)
we arrive at:
$$
R_y-\max_{z\in J}R_z \frac{N}{N-1}V^{\circ}(x,J)
\leq
\P^{\circ}\left(\t^0_x<\t^{0}_{(J\setminus x)\cup 0}\right)
\leq
R_x
\frac{1}{1-\max_{z\in J}V^{\circ}(z,J)}
\Eq(P.l9.10)
$$

\remark We could of course iterate the use of \eqv(P.l9.4) in \eqv(P.l9.5)
to bound both the numerator and the denominator of \eqv(P.l9.3) but we do not gain
much (the maximum in the r\.h\.s\. of \eqv(P.l9.10) would be raised to some power).

It now remains to estimate the ratios \eqv(P.l9.11). But this is simple
since by reversibility,
$$
R_x=\frac{\Q^{\circ}_N(x)\P^{\circ}\left(\t^x_0<\t^{x}_{x}\right)}
{\sum_{y\in J}\Q^{\circ}_N(y)\P^{\circ}\left(\t^y_0<\t^{y}_{y}\right)}
\Eq(P.l9.12)
$$
and by Theorem  \thv(P.prop.1),
$$
c_N^-\overline R
\leq
R_x
\leq
c_N^+\overline R
\Eq(P.l9.13)
$$
where $c_N^{\pm}$ are defined in \eqv(P.l9.2) and
$$
\overline R\equiv\frac{\Q^{\circ}_N(x)}
{\sum_{y\in J}\Q^{\circ}_N(y)}
\Eq(P.l9.14)
$$
Now since $J\subseteq \g_{I}(I)$, and since $\Q^{\circ}_N(y)=2^{-N}$
for all $y\in \g_{I}(I)$,
$$
\overline R
=\frac{1}{|J|}
\Eq(P.l9.15)
$$
Collecting \eqv(P.l9.10), \eqv(P.l9.13), and \eqv(P.l9.15) yields
\eqv(P.l9.1). By \eqv(N.4), \eqv(P.l9.1) remains true with $V^{\circ}(\,.\,,J)$
replaced by $U^{\circ}(\,.\,,J)$.
This concludes the proof of Lemma \thv(P.lemma.9). \endproof



\bigskip
\line{\bf 4.3. The harmonic measure $H^{\circ}_J(x,y)$.\hfill}

We now turn to the estimate of the general hitting probabilities \eqv(N.10).

%
%

\theo{\TH(P.theo.00)} {\it
Let $d\leq d_0(N)$.
Then, for all $J\in\SS_d$, all $x\in J$, and all $y\in\G_{N,d}\setminus J$,
$$
\frac{c_N^-}{|J|}\Bigl[1-(1+O(\sfrac{1}{N}))V^{\circ}(x,J)\Bigr]\left(1-V^{\circ}(y,J)\right)
\leq
H^{\circ}_J(y,x)
\leq
\frac{c_N^+}{|J|}\Bigl[1-\max_{z\in J} V^{\circ}(z,J)\Bigr]^{-1}+\phi_x(\dist(y,x))
\Eq(P.t0.1)
$$
where $c_N^{\pm}$ are defined in \eqv(P.l9.2).
Moreover \eqv(P.t0.1) remains true with $V^{\circ}(\,.\,,J)$ replaced by $U^{\circ}(\,.\,,J)$.
}


\proofof{Theorem \thv(P.theo.00)}
$$
\eqalign{
\P^{\circ}(\t^y_x<\t^y_{J\setminus x})
=&\P^{\circ}(\t^y_0<\t^y_x<\t^y_{J\setminus x})+\P^{\circ}(\t^y_x<\t^y_{(J\setminus x)\cup 0})
\cr
=&\P^{\circ}(\t^y_0<\t^y_J)\P^{\circ}(\t^0_x<\t^0_{J\setminus x})+\P^{\circ}(\t^y_x<\t^y_{(J\setminus x)\cup 0})
}
\Eq(P.t0.2)
$$
This immediately yields the upper bound
$$
\eqalign{
\P^{\circ}(\t^y_x<\t^y_{J\setminus x})
\leq & \P^{\circ}(\t^0_x<\t^0_{J\setminus x})+\P^{\circ}(\t^y_x<\t^y_{0})
\cr
\leq & H^{\circ}_J(0,x)+\phi_x(\dist(y,x))
}
\Eq(P.t0.3)
$$
To bound  $H^{\circ}_J(x,y)$ from below we use \eqv(P.t0.2) to write that
$$
\eqalign{
\P^{\circ}(\t^y_x<\t^y_{J\setminus x})
\geq &\P^{\circ}(\t^y_0<\t^y_J)H^{\circ}_J(0,x)
}
\Eq(P.t0.4)
$$
which together with
$$
\eqalign{
1-\P^{\circ}(\t^y_0<\t^y_{J})
=&\sum_{z\in J}\P^{\circ}(\t^y_z<\t^y_{0\cup(J\setminus z)})
\cr
\leq &\sum_{z\in J}\P^{\circ}(\t^y_z<\t^y_0)
\cr
\leq &\max_{z\in J}V^{\circ}(z,J)
}
\Eq(P.t0.5)
$$
gives
$$
\eqalign{
\P^{\circ}(\t^y_x<\t^y_{J\setminus x})
\geq \left(1-\max_{z\in J}V^{\circ}(z,J)\right)H^{\circ}_J(0,x)
}
\Eq(P.t0.6)
$$
The bounds \eqv(P.t0.1) then follow from \eqv(P.t0.3) and \eqv(P.t0.6)
and the bounds on $H^{\circ}_J(0,x)$  of Lemma \thv(P.lemma.9). \endproof

\bigskip
\line{\bf 4.4. `No return before hitting' probabilities.\hfill}

In Section 3.1 we proved upper bounds on `no return before hitting' probabilities of
the general form $\P^{\circ}(\t^x_{J\setminus x}<\t^x_x)$ for $J\subset\SS_d$ and
$x\in \G_{N,d}$ (Lemma \thv(P.lemma.3)). We now complement this result with
a lower bound in the case $x\in \SS_d$.

\theo{\TH(P.theo.0)}{\it
Let $d\leq d_0(N)$.
Then, for all $J\in\SS_d$ and all $x\in J$, the following holds:
\item{i)}
$$
\P^{\circ}(\t^x_{J\setminus x}< \t^x_x)\geq
\left(1-\frac{1}{N}-\frac{c}{N^2}-V^{\circ}(x,J)\right)
\left(1-\frac{1}{|J|}\right)
\Eq(P.t0.11)
$$
where $c_N^-$ are defined in \eqv(P.l9.2). Moreover \eqv(P.t0.11)
remains true with $V^{\circ}(\,.\,,J)$ replaced by $U^{\circ}(\,.\,,J)$.
\item{ii)} if $H(J)$ is satisfied
then
$$
\P^{\circ}(\t^x_{J\setminus x}<\t^x_x)
\leq \left(1-\frac{1}{|J|}\right)\left(1-\frac{1}{N}\right)
\Eq(P.t0.7)
$$
otherwise, if $H(J)$ is not satisfied,
$$
\P^{\circ}(\t^x_{J\setminus x}<\t^x_x)
\leq  \left(1-\frac{1}{|J|}\right)\left(1+O\Big(\frac{1}{N}\Big)\right)
\Eq(P.t0.7bis)
$$
}

\proofof{Theorem \thv(P.theo.0)}
The upper bounds \eqv(P.t0.7) and \eqv(P.t0.7bis) where established
in assertion (i) of Corollary \thv(P.cor.4). To prove
\eqv(P.t0.11), we write
$$
\eqalign{
\P^{\circ}(\t^x_{J\setminus x}< \t^x_x)
\geq &\P^{\circ}(\t^x_0<\t^x_{J\setminus x}< \t^x_x)
\cr
=&\P^{\circ}(\t^x_0<\t^x_J)
\P^{\circ}(\t^0_{J\setminus x}< \t^0_x)
\cr
=&
\left[1-\P^{\circ}(\t^x_J<\t^x_0)\right]
\left[1-\P^{\circ}(\t^0_x<\t^0_{J\setminus x})\right]
\cr
=&
\left[1-\P^{\circ}(\t^x_J<\t^x_0)\right]
\left[1-H^{\circ}_J(0,x)\right]
}
\Eq(P.t0.7')
$$
Now
$$
\eqalign{
1-\P^{\circ}(\t^x_J<\t^x_0)
=& \P^{\circ}(\t^x_0<\t^x_x)-\sum_{y\in J\setminus x}\P^{\circ}(\t^x_y<\t^x_{(J\setminus y)\cup 0})
\cr
\geq & \P^{\circ}(\t^x_0<\t^x_x)-\sum_{y\in J\setminus x}\P^{\circ}(\t^x_y<\t^x_0)
\cr
\geq & 1-\frac{1}{N}-\frac{c}{N^2}-V^{\circ}(x,J)
}
\Eq(P.t0.8)
$$
 where the last line follows from \eqv(P.p1.2)
of Theorem  \thv(P.prop.1), and where $0<c<4$ is a numerical constant.
Thus,
$$
\P^{\circ}(\t^x_{J\setminus x}< \t^x_x)\geq
\left(1-\frac{1}{N}-\frac{c}{N^2}-V^{\circ}(x,J)\right)
\left(1-H^{\circ}_J(0,x)\right)
\Eq(P.t0.9)
$$
where, by \eqv(P.l9.1) of Lemma \thv(P.lemma.9),
$$
\eqalign{
1-H^{\circ}_J(0,x)
\geq &
1-\frac{c_N^-}{|J|}\Bigl[1-(1+O(\sfrac{1}{N}))V^{\circ}(x,J)\Bigr]
\cr
= &
1-\frac{1}{|J|}+\frac{c}{|J|N^2}+
\frac{c_N^-}{|J|}(1+O(\sfrac{1}{N}))V^{\circ}(x,J)
\cr
\geq &
1-\frac{1}{|J|}
}
\Eq(P.t0.10)
$$
The lower bound \eqv(P.t0.11)  now follows from \eqv(P.t0.9) and
\eqv(P.t0.10).
\endproof

\bigskip
\vfill\eject

\bigskip
\chap{5. Back to the hypercube $\SS_N$.}5

Let a $d$-lumping $\g$ be given and  consider the corresponding lumped chain.
In this chapter we show how the results of Chapter 4, obtained for such
lumped chains, can be used to obtain
estimates on hitting probabilities for the ordinary random walk on $\SS_N$.
Clearly our key tool will be Lemma \thv(L.6) that states that
$$
\P\left(\t^{\s}_{A}\leq\t^{\s}_{B}\right)=
\P^{\circ}\left(\t^{\g(\s)}_{\g(A)}\leq
\t^{\g(\s)}_{\g(B)}\right)\,,\text{for all}\s\in\SS_N
$$
provided that $A\cup B$ is compatible with $\g$. More precisely,
in analogy with Definition \thv(I.def.1) and with the notation
therein:

\definition{\TH(I.def.1bis)}{\it  A subset $A$ of $\SS_N$ is called
$\g^{\L,\xi}$-compatible if and only if there exists a partition
$\L$ and a point $\xi\in\SS_N$ such that
A is $(\L,\xi)$-compatible.
}

As usual we will drop the superscripts $\L,\xi$ and assume that $\xi$ is
the point whose components are all equal to 1.
Inspecting the expressions of our various bounds on
hitting probabilities for the lumped chain, we see that
the only lumping-dependent quantities\note{
Note that $\phi_x(\dist(y,x))=V^{\circ}(y,x\cup y)$.
}
(i\.e\. $\g$-dependent quantities) are
the functions $V^{\circ}(y,J)$ and $U^{\circ}(y,J)$ defined in \eqv(N.3bis)
for subsets $J$ of the lumped state space $\G_{N,d}$.

The aim of Section 5.1 below is to show that these functions have
equivalent expressions in the hypercube setting. At the same time this
will allow us to draw the correspondence between the notions of
sparseness in these two different spaces. The same question is
addressed for so-called  Hypothesis $H$.
Section 5.2 is then devoted to the statements and proofs of
a number of results for
the random walk on $\SS_N$. It contains in particular the proofs of
Theorem \thv(I.theo.1) and Corollary \thv(I.cor.theo.1) of Chapter 1.

\bigskip
\line{\bf 5.1. Sparseness and Hypothesis $H$: from the hypercube $\SS_N$ to the
grid $\G_{N,d}$.\hfill}

Recall from \eqv(N.0) that, given two points $x,y\in \G_{N,d}$, $\dist(x,y)$ denotes
the graph distance,
$$
\dist(x,y)\equiv\sum_{k=1}^d\frac{|\L_k|}{2}|x^k-y^k|
$$
The following elementary but key lemma states that whenever
the distance is measured from a
vertex $x\in\SS_d$, the lumping function is distance preserving.

\lemma{\TH(N.lemma.3)}{\it For all $x\in\SS_d$ and $y\in \G_{N,d}$,
for all $\s,\eta\in\SS_N$ such that
$$
\g(\s)=x\,,\quad\g(\eta)=y
\Eq(N.l3.1)
$$
we have
$$
\dist(x,y)=\Dist(\s,\eta)
\Eq(N.l3.2)
$$
}

\proof Immediate.\endproof



Recall that for $J\subset\G_{N,d}$ and $y\in\G_{N,d}$,
$$
\eqalign{
V^{\circ}(y,J)&
=\cases\sum_{z\in J\setminus y}\phi_z(\dist(y,z)), &\text{if} J\setminus y\neq\emptyset\cr
0,&\text{otherwise}
\endcases
\cr
U^{\circ}(y,J)&=\cases \sum_{z\in J\setminus y}F(\dist(y,z)), &\text{if} J\setminus y\neq\emptyset\cr
0,&\text{otherwise}
\endcases
\cr
}
\Eq(N.3bis-repetee)
$$
and define, for  $A\subset\SS_N$ and $\s\in\SS_N$,
$$
\eqalign{
V(\s,A)&
=\cases\sum_{\eta\in A\setminus \s}\phi_{\g(\eta)}(\Dist(\s,\eta)), &\text{if} A\setminus \s\neq\emptyset\cr
0,&\text{otherwise}
\endcases
\cr
U(\s,A)&=\cases \sum_{\eta\in A\setminus \s}F(\Dist(\s,\eta)), &\text{if} A\setminus \s\neq\emptyset\cr
0,&\text{otherwise}
\endcases
\cr
}
\Eq(B.1)
$$
(where $\g$ in the definition of $V(\s,A)$ is the same lumping function as that
used in \eqv(N.3bis-repetee)). Note that in \eqv(N.3bis-repetee) we allow for the
possibility that $y\in J$; similarly, in \eqv(B.1), we may have $\s\in A$.

\remark Recall that for the sake of brevity we chose to drop the indices $N$ and $d$
and write $F\equiv F_{N,d}$, except in the statement and proofs of the main results from Sections 1, 5, 6, and 7.
The same notational rule applies to the functions $V^{\circ}$, $U^{\circ}$, $V$, and $U$ from
\eqv(N.3bis-repetee) and \eqv(B.1),
which will gain back the indices $N$ and $d$ whenever $F$ does.
The same again  applies to the functions
$\UU$, $\VV$, $\UU^{\circ}$, and $\VV^{\circ}$ that will shortly be defined
(see \eqv(I.12repetee)-\eqv(I.12bis)).

As a first consequence of Lemma \thv(N.lemma.3) we have:

\lemma{\TH(B.lemma.2)}{\it
For all $\g$-compatible subset $A\subset\SS_N$,
%
for all pairs of points $y\in\G_{N,d}$ and $\s\in\SS_N$ such that $\g(\s)=y$,
we have, setting $J=\g(A)\subset\SS_d$,
$$
\eqalign{
V^{\circ}(y,J)&=V(\s,A)\cr
U^{\circ}(y,J)&=U(\s,A)\cr
}
\Eq(B.2)
$$
}

\proof Immediate using  \eqv(N.3bis-repetee), \eqv(B.1),
and Lemma \thv(N.lemma.3).\endproof

Note now that  among the quantities defined in \eqv(N.3bis-repetee), \eqv(B.1),
the only one that does not depend on the underlying lumping function $\g$ is
$U(\s,A)$. The next lemma shows how to pass from $V(\s,A)$ to $U(\s,A)$.

\lemma{\TH(B.lemma.3)}{\it
Assume that $\g$ is generated by a
log-regular $d$-partition. Then, for all $\g$-compatible subset $A\subset\SS_N$
and for all  $\s\in\SS_N$,
$$
V(\s,A)\leq U(\s,A)
\Eq(B.3)
$$
}

\proof Set $J=\g(A)$ and $y=\g(\s)$. By assumption $J\subset\SS_d$ and $y\in\G_{N,d}$.
We proved
in \eqv(N.4) that if the $d$-partition generating $\g$ is
log-regular then,
$$
V^{\circ}(y,J)\leq U^{\circ}(y,J)\,\text{for all}
J\subset\SS_d\,,\, y\in\G_{N,d}
\Eq(N.4repetee)
$$
But this and \eqv(B.2) prove \eqv(B.3).\endproof

We now want to relate the sparsity condition, defined in \eqv(I.13)
for subsets $A\subset\SS_N$, to corresponding quantities in the lumped state
space $\G_{N,d}$.
To this aim
recall that for $ A\subset\SS_N$,
$$
\UU(A)= \max_{\eta\in A}U(\eta,A)
\Eq(I.12repetee)
$$
and set
$$
\VV(A)=\max_{\eta\in A}V(\eta,A)
\Eq(I.12new)
$$
Similarly, for $J\subset\SS_d$, define
$$
\UU^{\circ}(J)=\max_{x\in J}U^{\circ}(x,J)
\Eq(I.12ter)
$$
$$
\VV^{\circ}(J)=\max_{x\in J}V^{\circ}(x,J)
\Eq(I.12bis)
$$

\remark Again, among the quantities \eqv(I.12repetee), \eqv(I.12ter),
\eqv(I.12bis), and \eqv(I.12new),
the only one that does not depend on $\g$ is \eqv(I.12repetee).


\vfill\eject

\lemma{\TH(B.4)}{\it
%
For all $\g$-compatible subset $A\subset\SS_N$,
setting $J=\g(A)\subset\SS_d$,
%
%
%
$$
\UU^{\circ}(J)=\UU(A)
\Eq(TL.l4.5)
$$
$$
\VV^{\circ}(J) =\VV(A)
\Eq(TL.l4.5bis)
$$
}

\proof This follows from Lemma \thv(B.lemma.2) and the
definitions \eqv(I.12repetee)-\eqv(I.12bis).\endproof

Naturally, we will say that a subset $J\subset\G_{N,d}$ is
$(\e,d)$-sparse if there exists $\e>0$ such that
$
\UU^{\circ}(J)\leq \e
$.
Thus, \eqv(TL.l4.5) entails that

\corollary{\TH(B.cor.1)}{\it
For all $\g$-compatible subset $A\subset\SS_N$,
setting $J=\g(A)\subset\SS_d$,
%
$A$ is $(\e,d)$-sparse
if and only $J$ is $(\e,d)$-sparse.
}

As in Lemma \thv(B.lemma.3) the following lemma will be used to pass to the
($\g$ independent) function $\UU$.

\lemma{\TH(TL.lemma.4)}{\it Assume that $\g$ is generated by
a log-regular $d$-partition. Then,
for all $\g$-compatible subset $A\subset\SS_N$,
$$
\VV(A)\leq\UU(A)
\Eq(N.4again)
$$
}

\proof Set $J=\g(A)$. By \eqv(N.4repetee),
$
\VV^{\circ}(J)\leq\UU^{\circ}(J)
$,
and combining with Lemma \thv(B.4),
$$
\VV(A)=\VV^{\circ}(J)\leq\UU^{\circ}(J)=\UU(A)
\Eq(TL.l4.1)
$$
proving \eqv(N.4again).\endproof



We finally conclude this section by comparing Hypothesis $H$ and $H^{\circ}$,
defined respectively in \eqv(P.H) and \eqv(P.H'lump).

\lemma{\TH(B.lemma.5)}{\it Under the assumptions and with the notation of
Lemma \thv(B.lemma.2)
the following holds.
$A\cup\s$ satisfies hypothesis $H(A\cup\s)$ if and only if
$J\cup y$ satisfies hypothesis $H^{\circ}(J\cup y)$.
}

\proof This is again an immediate consequence of Lemma \thv(N.lemma.3).\endproof

Note that in the statement above we allow for the possibility that $y\notin\SS_d$.



\vfill\eject
\bigskip
\line{\bf 5.2. Main results.\hfill}

\noindent{\bf The harmonic measure.} We begin by giving a general result from which
Theorem \thv(I.theo.1) and Corollary \thv(I.cor.theo.1) will be derived.

\theo{\TH(B.theo.00)} {\it
Let $d\leq d_0(N)$.
Given a $d$-lumping $\g$ and a $\g$-compatible
subset $A\subset\SS_N$ we have,
for all $\eta\in A$, and all $\s\in\SS_N\setminus A$,
$$
\eqalign{
H_A(\s,\eta)\geq &
\frac{c_N^-}{|A|}\Bigl[1-(1+O(\sfrac{1}{N}))V_{N,d}(\eta,A)\Bigr]\left(1-V_{N,d}(\s,A)\right)\cr
H_A(\s,\eta)\leq &
\frac{c_N^+}{|A|}\Bigl[1-\max_{\eta'\in A} V_{N,d}(\eta',A)\Bigr]^{-1}+\phi_{\eta}(\Dist(\s,\eta))
}
\Eq(B.10)
$$
where $c_N^{\pm}$ are defined in \eqv(P.l9.2).
Moreover \eqv(B.10) remains true
with either of the following changes:
\item{1)}  replacing $V_{N,d}(\,.\,,A)$ by $U_{N,d}(\,.\,,A)$ and $\phi_{\eta}$ by $F_{N,d}$;
\item{2)} replacing $V_{N,d}(\eta\,,A)$ and $\max_{\eta'\in A} V_{N,d}(\eta',A)$ by $\UU_{N,d}(A)$,
$V_{N,d}(\s\,,A)$ by $U_{N,d}(\s\,,A)$, and $\phi_{\eta}$ by $F_{N,d}$.
}

Note that in case 2), the expressions of the bounds \eqv(B.10)  become independent of $\g$.

\proof With the notation of Theorem \thv(B.theo.00) set $J=\g(A)$, $x=\g(\eta)$, and
$y=\g(\s)$. By assumption $x\in J\subset \SS_d$ and  $y\in\G_{N,d}\setminus J$. Next,
by Lemma  \thv(L.5), $H_A(\s,\eta)=H^{\circ}_J(y,x)$; by Lemma \thv(B.lemma.2),
$V_{N,d}(\,.\,,A)=V^{\circ}_{N,d}(\,.\,,J)$; and, by Lemma \thv(N.lemma.3),
$\Dist(\s,\eta)=\dist(y,x)$. The bounds \eqv(B.10) now follow from
Theorem \eqv(P.theo.00).
From this, lemma \thv(B.lemma.3),  \eqv(I.12repetee)-\eqv(I.12bis),
and Lemma \thv(TL.lemma.4), Assertion 1) and 2) follow.
\endproof

\proofof{Theorem \thv(I.theo.1)}
Let the notation be as in  Theorem \thv(B.theo.00) and let $0\leq\rho\leq N$ be given. Consider \eqv(B.10).
Using Lemma \thv(N.lemma.1) and Lemma \thv(N.lemma.2) successively we have,
for all $\s$ satisfying $\Dist(\s,A)>\rho$,
$$
\phi_{\g(\eta)}(\Dist(\s,\eta))
\leq
\phi_{\g(\eta)}(\rho+1)\leq F_{N,d}(\rho+1)
\Eq(B.10quart)
$$
Moreover this and the definition of $V_{N,d}(\,.\,,A)$ yields
$V_{N,d}(\s,A)\leq |A| F_{N,d}(\rho+1)$  for all $\s$ satisfying $\Dist(\s,A)>\rho$.
We may thus replace $V_{N,d}(\s\,,A)$ by $|A| F_{N,d}(\rho+1)$ and
$\phi_{\g(\eta)}$ by $F_{N,d}(\rho+1)$ in \eqv(B.10). Now, by assertion 2)
of Theorem \thv(B.theo.00) we also may replace
$V_{N,d}(\eta\,,A)$ and $\max_{\eta'\in A} V_{N,d}(\eta',A)$ by
$\UU_{N,d}(A)$ in \eqv(B.10).
Doing so yields
$$
\frac{c_N^-}{|A|}\Bigl[1-(1+O(\sfrac{1}{N}))\UU_{N,d}(A)\Bigr]
\left(1-|A| F(\rho)\right)
\leq
H_A(\s,\eta)
\leq
\frac{c_N^+}{|A|}\Bigl[1-\UU_{N,d}(A)\Bigr]^{-1}
+F_{N,d}(\rho)
\Eq(B.10bis)
$$
and, setting
$\vartheta_{N,d}(A,\rho)\equiv \max\left\{\UU_{N,d}(A), |A|F_{N,d}(\rho+1)\right\}$
we obtain,
$$
\frac{1}{|A|}(1-c^-\vartheta_{N,d}(A,\rho))
\leq
H_A(\s,\eta)
\leq
\frac{1}{|A|}(1+c^+\vartheta_{N,d}(A,\rho))
\Eq(B.10ter)
$$
for some finite positive constants $c^+, c^-$.
%
Theorem \thv(I.theo.1) is proven.\endproof

\proofof{Corollary \thv(I.cor.theo.1)}
If $A\subset\SS_N$ is such that $2^{|A|}\leq d_0(N)$ then,
by Corollary \eqv(I.cor.1), there exists a $d$-partition $\L$ with
$d\leq d_0(N)$ such that $A$ is $\L$-compatible
and $\UU_{N,d}(A)\leq c\frac{1}{(\log N)^2}$
for some constant $0<c<\infty$.
Next, since
$F_1(1)\leq\kappa_0/N$, using that by \eqv(A3.2), $F_2(1)\leq\frac{2\kappa^2_0}{(\log N)^3}$,
we get that $|A|F(1)\leq |A|(F_1(1)+F_2(1))\leq \frac{c'}{(\log N)^2}$
where $0<c'<\infty$. Let us thus choose $\rho=0$ in the statement of Theorem \thv(I.theo.1).
First note that this implies that \eqv(I.16)
is satisfied uniformly in $\s$ for $\s\notin A$.
Next, putting together our bounds on $\UU_{N,d}(A)$ and $|A|F(1)$
gives
$
\vartheta_{N,d}(A,\rho)\leq (c+c')\frac{1}{(\log N)^2}
$
Finally, inserting the latter bound in \eqv(I.16) yields \eqv(I.16bis).
The corollary is proven.
\endproof
%


\noindent{\bf `No return before hitting' probabilities.} Let us now consider hitting
probabilities of the form $\P(\t^{\eta}_{A\setminus \eta}< \t^{\eta}_{\eta})$ for
$A\subset\SS_N$.


\theo{\TH(I.theo.2)}{\it Let $d\leq d_0(N)$.
Given a $d$-lumping $\g$ and a $\g$-compatible
subset $A\subset\SS_N$,
the following holds for all $\eta\in A$:

\item{i)} If $H(A)$ is satisfied then, for all $\eta\in A$,
$$
\left(1-\frac{1}{N}-\frac{c}{N^2}-V(\eta,A)\right)\left(1-\frac{1}{|A|}\right)
\leq
\P\left(\t^{\eta}_{A\setminus\eta}<\t^{\eta}_{\eta}\right)
\leq
\left(1-\frac{1}{N}\right)\left(1-\frac{1}{|A|}\right)
\Eq(I.17)
$$
where $0<c<\infty$ is a numerical constant, whereas if $H(A)$ is not satisfied
the lower bound in \eqv(I.17) remains unchanged,
but the term $1-\frac{1}{N}$ in the upper bound is replaced by $1+O(\frac{1}{N})$.

\item{ii)}  In addition assertion i) remains true with either
of the following changes in the lower bound:
\itemitem{1)} replacing  $V(\eta,A)$ by $U(\eta,A)$;
\itemitem{2)} replacing $V(\eta,A)$  by $\UU_{N,d}(A)$.
}

\proof
With the notation of Theorem \thv(P.theo.0) set $J=\g(A)$ and $x=\g(\eta)$ for $\eta\in A$.
By assumption $x\in J\subset \SS_d$. Next,
by Lemma  \thv(L.5),
$\P\left(\t^{\eta}_{A\setminus\eta}<\t^{\eta}_{\eta}\right)=\P^{\circ}(\t^x_{J\setminus x}< \t^x_x)$;
by Lemma \thv(B.lemma.2),
$V(\eta,,A)=V^{\circ}(x,J)$; and, by Lemma \thv(B.lemma.5), $A$ satisfies hypothesis $H(A)$ if and only if
$J$ satisfies hypothesis $H^{\circ}(J)$. Assertion i) of Theorem \thv(I.theo.2) now follows from
Theorem \thv(P.theo.0).
From this, lemma \thv(B.lemma.3),  \eqv(I.12repetee)-\eqv(I.12bis),
and Lemma \thv(TL.lemma.4), Assertion ii) follows.
\endproof

Consider the case ii-2) in Theorem \thv(I.theo.2).
We see that, for sparse enough sets $A$,
the form of the hitting probability  undergoes a change when the
size of $A$ becomes, roughly, of order $N$. The next corollary
shows that
Theorem \thv(I.theo.2) can yield coinciding upper and lower bounds
uniformly in $\eta$ that are either close to $1-\frac{1}{|A|}$
or close to $1-\frac{1}{N}$.

\corollary{\TH(I.cor.theo.2)}{\it Under the assumptions of assertion ii) of
Theorem \thv(I.theo.2),
the following holds:
\item{(i)} if $\frac{|A|}{N}=o(1)$ and $\UU_{N,d}(A)=o\big(\frac{1}{|A|}\big)$ then, for all $\eta\in A$,
$$
1-\frac{1}{|A|}\left(1+o(1)\right)
\leq
\P(\t^{\eta}_{A\setminus \eta}< \t^{\eta}_{\eta})
\leq
1-\frac{1}{|A|}\left(1-o(1)\right)
\Eq(I.18)
$$
\item{(ii)} if $\frac{N}{|A|}=o(1)$ and
$
\UU_{N,d}(A)=o\big(\frac{1}{N}\big)
$
and if $H(A)$ is satisfied
then, for all $\eta\in A$,
$$
1-\frac{1}{N}\left(1+o(1)\right)
\leq
\P(\t^{\eta}_{A\setminus \eta}< \t^{\eta}_{\eta})
\leq
1-\frac{1}{N}\left(1-o(1)\right)
\Eq(I.20)
$$
}

\proof This is an immediate consequence of Theorem \thv(I.theo.2), ii)-2.\endproof

\remark To understand the difference between \eqv(I.18) and \eqv(I.20)
it is useful  to observe that $\P(\t^{\s}_{\s}=2)=\frac{1}{N}$.

\remark When  $\frac{N}{|A|}=o(1)$ and  $H(A)$ is not satisfied, we do not have
coinciding upper and lower bounds, nor do we have reasons to think that either of the
bounds \eqv(I.17) will, in general, be good. As we explained earlier, the behavior of
$\P(\t^{\eta}_{A\setminus \eta}< \t^{\eta}_{\eta})$ will depend on the structure
of the set $A$ locally (see  the proof of \eqv(P.l3.3') of lemma \eqv(P.lemma.3)
where this remark was made precise).
%



\bigskip
\vfill\eject

\chap{6. Mean times.}6

In this chapter we infer some basic estimates for the mean hitting
times in our model, both for the chain on the hypercube and for the lumped chain,
and prove Theorem \thv(I.theo.3)

We begin with the chain on the hypercube. Theorem \thv(T.theo.1)
below is more general than Theorem \thv(I.theo.3).

\theo{\TH(T.theo.1)} {\it Let $d'\leq d_0(N)/2$ and assume that
$A\subset\SS_N$ is compatible with a partition $\L'$ into $d'$ classes.
Then for all $\s\notin A$ there exits a partition $\L$ into $d$ classes,
with $d'<d\leq 2d'$, compatible with $A\cup\s$.
%
%
Let one such partition be fixed and set
$$
c_N^{\pm}=
1\pm \frac{c}{N^2}
\Eq(T.0)
$$
where $0<c<5$ is a numerical constant. Then, if $\VV_{N,d}(A\cup\s)<{c^-_N}/{2}$,
the following holds:
\item{(i)} if $H(\s\cup A)$ is satisfied
$$
\frac{2^{N}}{|A|(1-\frac{1}{N})}
c^-_N
\Bigl[1-(1+O(\sfrac{1}{N}))
\VV_{N,d}(A\cup\s)
\Bigr]
\leq
\E(\t^{\s}_{A})
\leq
\frac{2^{N}}{|A|(1-\frac{1}{N})}
c^+_N\left[c^-_N-2\VV_{N,d}(A\cup\s)\right]^{-1}
\Eq(T.1)
$$
whereas
if $H(\s\cup A)$ is not satisfied,
the term $1-\frac{1}{N}$ in
the lower bound must be replaced by $1+O(\frac{1}{N})$.

\item{(ii)} for all $\eta\in A$, if
$H(\s\cup A)$ is satisfied,
$$
\eqalign{
\frac{2^{N}}{|A|(1-\frac{1}{N})}c^-_N\Bigl[1-(1+O(\sfrac{1}{N}))\VV_{N,d}(A\cup\s)&\Bigr]^4\cr
\leq
\E\left(\t^{\s}_{\eta}\mid\t^{\s}_{\eta}<\t^{\s}_{A\setminus\eta}\right)
=&
\E\left(\t^{\s}_{A}\mid\t^{\s}_{\eta}<\t^{\s}_{A\setminus\eta}\right)
\cr
\leq&
\frac{2^{N}}{|A|(1-\frac{1}{N})}
c^+_N\left[c^-_N-2\VV_{N,d}(A\cup\s)\right]^{-4}
}
\Eq(T.2)
$$
whereas if $H(\s\cup A)$ is not satisfied,
the term $1-\frac{1}{N}$ in
the lower bound must be replaced by $1+O(\frac{1}{N})$.

Moreover statements (i) and (ii) remain true with $\VV_{N,d}(A\cup\s)$
replaced by $\UU_{N,d}(A\cup\s)$ (see \eqv(I.12repetee) and \eqv(I.12new)).
}


\remark
The only quantity in Theorem \thv(T.theo.1) that depends on
the choice of $d$-lumping $\g$ (or, equivalently, on the
$d$-partition $\L$) is $\VV_{N,d}(A\cup\s)$.
As usual passing from  $\VV_{N,d}(A\cup\s)$
to $\UU_{N,d}(A\cup\s)$ we get rid of this dependence.


We now state the lumped-chain version of Theorem \thv(T.theo.1).


\theo{\TH(T.theo.1.lumped)} {\it  Let $d\leq d_0(N)$.
Then, for all $d$-lumping $\g$ (or equivalently for all
$d$-partition $\L$), the following holds:
for all $I\subset\SS_d$, all $x\in\SS_d\setminus I$, and all $y\in I$,
$$
\E^{\circ}(\t^{x}_{I})
\text{and}
\E^{\circ}\left(\t^{x}_{y}\mid\t^{x}_{y}<\t^{x}_{I\setminus y}\right)
\Eq(T.lump1)
$$
obey the bounds obtained for
$$
\E(\t^{\s}_{A})
\text{and}
\E\left(\t^{\s}_{\eta}\mid\t^{\s}_{\eta}<\t^{\s}_{A\setminus\eta}\right)
\Eq(T.lump2)
$$
in statements (i) and (ii) of Theorem \thv(T.theo.1), with $|A|$,
$H(\s\cup A)$, $\VV_{N,d}(A\cup\s)$, and $\UU_{N,d}(A\cup\s)$
replaced, respectively, by $|I|$, $H^{\circ}(I\cup x)$, $\VV^{\circ}_{N,d}(I\cup x)$,
and $\UU^{\circ}_{N,d}(I\cup x)$, and with $c_N^{\pm}$ given by \eqv(T.0).
}

We will see in Section 7 that, for more detailed investigations of
the distributions of  hitting times (for both the chain on the
hypercube and the lumped chain) we  need to control some further
mean times in the lumped chain. This is the main motivation for our next theorem.


\theo{\TH(T.theo.2)} {\it Under the assumptions of Lemma \thv(P.lemma.7)
the following holds for all $y\in\G_{N,d}$:
For $d>1$,
$$
\E\t^y_0 \leq CN^{2}\prod_{k=1}^{d}|\L_k|\leq  CN^{d+2}
\Eq(T.3)
$$
for some constant $0<C<\infty$. If $d$ is finite
and independent of $N$, and if $\L$ is an equipartition, then
\eqv(T.3) can be refined to
$$
\E\t^y_0 \leq CN^{\frac{d+1}{2}}\log N
\Eq(T.3')
$$
for some constant $0<C<\infty$. Furthermore, if $d=1$,
$$
\E\t^y_0 \leq \E\t^1_0=\frac{N}{4}\log N (1+o(1))
\Eq(T.4)
$$
}

\remark The level of precision of \eqv(T.3') and  \eqv(T.4) is not needed in the sequel.


For later reference we set
$$
\wh\Th(d)
=\cases  CN^{2}\prod_{k=1}^{d}|\L_k|,&\hbox{if}\,\,\,d>1\cr
\frac{N}{4}\log N (1+o(1)),&\hbox{if}\,\,\,d=1\cr
\endcases
\Eq(T.4bis)
$$


Theorems \thv(T.theo.1), \thv(T.theo.1.lumped), and \thv(T.theo.2) follow from estimates of the
previous section and the following well-known
formulas
from potential theory, which hold for any discrete Markov chain (see e\.g\. \cite{So}),
but that we express here for the chain
on the hypercube:
for all subset
$A\subseteq\SS_N$ and all
$\s\in\SS_N$ such that $\s\notin A$,
$$
\E(\t^{\s}_{A})=
\frac{1}{\mu_{N}(\s)
\P(\t^{\s}_{A}<\t^{\s}_{\s})}
\left[\mu_{N}(\s)+
\sum_{\eta\in (A\cup\s)^c}\mu_{N}(\eta)
\P(\t^{\eta}_{\s}<\t^{\eta}_{A})
\right]
\Eq(T.5)
$$
and for all subsets $A,B\subseteq\SS_N$, and all
$\s\in\SS_N$ such that $\s\notin A\cup B$,
$$
\eqalign{
&
\E\left(\t^{\s}_{A}
\mid\t^{\s}_{A}\leq\t^{\s}_{B}
\right)
\cr
=&
\frac{1}{\mu_{N}(\s)
\P(\t^{\s}_{A\cup B}<\t^{\s}_{\s})}
\left[\mu_{N}(\s)+
\sum_{\eta\in (A\cup B\cup\s)^c}\mu_{N}(\eta)
\P(\t^{\eta}_{\s}<\t^{\eta}_{A\cup B})
\frac{\P(\t^{\eta}_{A}\leq\t^{\eta}_{B})}
{\P(\t^{\s}_{A}\leq\t^{\s}_{B})}
\right]
\cr
}
\Eq(T.6)
$$

We prove Theorem \thv(T.theo.1) and  \thv(T.theo.1.lumped) simultaneously.

\proofof{Theorem \thv(T.theo.1) and Theorem \thv(T.theo.1.lumped)}
As in Theorem \thv(T.theo.1) let $d'\leq d_0(N)/2$ and assume that
$A\subset\SS_N$ is compatible with a partition $\L'$ into $d'$ classes.
We first want to see that for all $\s\notin A$ there exits a partition $\L$ into $d$ classes
that satisfies $d'<d\leq 2d'$ and is compatible with $A\cup\s$. This is simple.
Given $\s\notin A$ let $\L$ be the partition obtained as follows: split each
class $\L'_k$, $1\leq k\leq d'$, into two non-empty classes
$\L^+_k$ and $\L^-_k$, where $\L^{\pm}_k=\{i\in\L'_k\mid \s_i=\pm 1\}$
if and only if none of these sets is empty, and if one of them
is empty then leave $\L'_k$ unchanged. Clearly this partition is compatible
with $A\cup\s$ and $d$ satisfies $d'<d\leq 2d'$. Now choose one such $d$-partition $\L$
and let $\g$ be the $d$-lumping generated by $\L$. Set
%
%
$x=\g(\s)$, $I=\g(A)$ and, for $\eta\in A$, $y=\g(\eta)$.
By virtue of Lemma \thv(L.5),
$$
\E(\t^{\s}_{A})=\E^{\circ}(\t^{x}_{I})
\Eq(T.7ter)
$$
$$
\E\left(\t^{\s}_{\eta}\mid\t^{\s}_{\eta}<\t^{\s}_{A\setminus\eta}\right)=
\E^{\circ}\left(\t^{x}_{y}\mid\t^{x}_{y}<\t^{x}_{I\setminus y}\right)
\Eq(T.20bis)
$$
Moreover, since $\L$ is compatible with $A\cup\s$,
$x\in\SS_d$, $I\subset\SS_d$, and $y\in I$. Finally it
follows from \eqv(TL.l4.5bis) of Lemma \thv(B.4) that $\VV^{\circ}_{N,d}(I)=\VV_{N,d}(A)$,
while by Lemma \thv(N.lemma.3), $\Dist(\s,A)=\dist(x,I)$.
From this we conclude that  Theorem \thv(T.theo.1) and Theorem
\thv(T.theo.1.lumped) are equivalent.

We now prove Theorem \thv(T.theo.1). Since the proofs of the two assertions
are very similar we will prove the first assertion in detail, but only sketch the second.
We start with the proof of assertion (i).
By definition of $\mu_{N}$, \eqv(T.5) reads
$$
\E(\t^{\s}_{A})=
\frac{1}{\P(\t^{\s}_{A}<\t^{\s}_{\s})}
\left[1+\sum_{\s'\in (A\cup\s)^c}
\P(\t^{\s'}_{\s}<\t^{\s'}_{A})
\right]
\Eq(T.7)
$$
Setting $J=I\cup x=\g(A\cup\s)\in\SS_d$
and using Lemma \thv(L.6) to pass to the lumped chain
we obtain
(recall the notation \eqv(N.10) for the harmonic measure),
$$
\E(\t^{\s}_{A})=\frac{1}{\P^{\circ}(\t^x_{J\setminus x}< \t^x_x)}
\left[1+\sum_{\s'\in (A\cup\s)^c}H^{\circ}_J(\g(\s'),x)\right]
\Eq(T.8)
$$
Using Theorem \thv(P.theo.0) to express $\P^{\circ}(\t^x_{J\setminus x}< \t^x_x)$,
we have to evaluate the sum appearing in the r\.h\.s. of \eqv(T.8). From the
upper and lower bounds of Theorem \thv(P.theo.00), setting
$$
\eqalign{
R_1(x) & \equiv \frac{|J|}{2^{N}}\sum_{\s'\in (A\cup\s)^c}\phi_x(\dist(\g(\s'),x))
\cr
R_2(J) & \equiv \frac{1}{2^{N}}\sum_{\s'\in (A\cup\s)^c}V^{\circ}(\g(\s'),J)
}
\Eq(T.100)
$$
we deduce that
$$
\eqalign{
\sum_{\s'\in (A\cup\s)^c}H^{\circ}_J(\g(\s'),x)
\leq & c_N^+\frac{2^{N}-|J|}{|J|}
\Bigl[1-\max_{z\in J} V^{\circ}(z,J)\Bigr]^{-1} +\frac{2^{N}}{|J|}R_1(x)
\cr
}
\Eq(T.9)
$$
$$
\eqalign{
\sum_{\s'\in (A\cup\s)^c}H^{\circ}_J(\g(\s'),x)
\geq & c_N^-\frac{2^{N}}{|J|}
\Bigl[1-(1+O(\sfrac{1}{N}))V^{\circ}(x,J)\Bigr]
\left(1-R_2(J)\right)
\cr
}
\Eq(T.10)
$$
To evaluate $R_1(x)$ note that, by Theorem  \thv(P.prop.2),
$$
\eqalign{
\sum_{\s'\in (A\cup\s)^c}\phi_x(\dist(\g(\s'),x))
\leq&
\sum_{\s'\in (A\cup\s)^c}F(\dist(\g(\s'),x))
\cr
\leq&
\sum_{\s'\neq\s}F(\dist(\g(\s'),x))
\cr
=&
\sum_{n=1}^N\binom{N}{n} F(n)
\cr
}
\Eq(T.11)
$$
where, as defined in \eqv(P.p2.02), $F(n)=F_1(n)+F_2(n)$.
By  \eqv(P.p2.03), $\sum_{n=1}^N\binom{N}{n} F_1(n)\leq CN$ for
some constant $C<\infty$, and by Lemma \thv(A3.lemma.1), using that
$$
\eqalign{
\sum_{n=1}^N\binom{N}{n}\left(\frac{n}{N}\right)^{\frac{n}{2}}
\leq&
N a^{\frac{N}{2}}e^{-\frac{N}{2}(1-1/a)}\,,\quad 0\leq a(\approx 1.82)<2
\cr
\sum_{n=1}^{d-2}\binom{N}{n}\left(\frac{d}{N}\right)^{n}
\leq&
be^{d}\,,\quad 0\leq b<\infty
}
\Eq(T.12)
$$
we obtain
$$
\sum_{n=1}^N\binom{N}{n} F_2(n)\leq N^6(a^{\frac{N}{2}}+be^{d})
\Eq(T.12bis)
$$
Since by construction $J\in\SS_d$, and since $|\SS_d|=2^d$ where,
by assumption
$d\leq 2d'\leq d_0(N)\leq \a_0\frac{N}{\log N}$ for some
constant $0<\a_0<1$,
$$
R_1(x)\leq
\frac{|J|}{2^{N}}\sum_{n=1}^N\binom{N}{n} F(n)
\leq
\frac{2^d}{2^{N}}  N^6(a^{\frac{N}{2}}+be^{d}+1)\leq e^{-N/2}
\Eq(T.13)
$$
Similarly, by definition of $V^{\circ}(\g(\s'),J)$ (see \eqv(N.3bis)),
$$
\eqalign{
R_2(J)
=&
\frac{1}{2^{N}}\sum_{\s'\in (A\cup\s)^c}\sum_{\eta\in A\cup\s}\phi_{\g(\eta)}(\dist(\g(\s'),\g(\eta)))
\cr
\leq &
\frac{1}{2^{N}}\sum_{\eta\in\A\cup\s}\sum_{\s'\neq\eta}F(\dist(\g(\s'),\g(\eta)))
\cr
\leq & \frac{|J|}{2^{N}}\sum_{n=1}^N\binom{N}{n} F(n)
\cr
\leq & e^{-N/2}
}
\Eq(T.14)
$$
where the one before last line is obtained just as the last line of \eqv(T.11).
Inserting the previous two bound in \eqv(T.9) and \eqv(T.10) yields
$$
\eqalign{
\sum_{\s'\in (A\cup\s)^c}H^{\circ}_J(\g(\s'),x)
\leq & \frac{2^{N}}{|J|}
\left(
c_N^+\Bigl[1-\max_{z\in J} V^{\circ}(z,J)\Bigr]^{-1} + e^{-N/2}
\right)
\cr
\leq &
(c_N^+ + e^{-N/2})\Bigl[1-\max_{z\in J} V^{\circ}(z,J)\Bigr]^{-1}
\cr
}
\Eq(T.16)
$$
and
$$
\eqalign{
\sum_{\s'\in (A\cup\s)^c}H^{\circ}_J(\g(\s'),x)
\geq & \frac{2^{N}}{|J|}c_N^-(1-e^{-N/2})
\Bigl[1-(1+O(\sfrac{1}{N}))V^{\circ}(x,J)\Bigr]
}
\Eq(T.17)
$$
Finally, plugging \eqv(T.16) and \eqv(T.17) in \eqv(T.8), and using
Theorem \thv(P.theo.0)
to bound $\P^{\circ}(\t^x_{J\setminus x}< \t^x_x)$,
we obtain, for some constant $0<c'<5$,
$$
\eqalign{
\E(\t^{\s}_{A})
\leq &
\frac{2^{N}}{|J|-1}
\frac{c_N^+ + e^{-N/2}}
{\left[1-\frac{1}{N}-\frac{c}{N^2}-V^{\circ}(x,J)\right]
\Bigl[1-\max_{z\in J} V^{\circ}(z,J)\Bigr]}
\cr
\leq &
\frac{2^{N}}{(|J|-1)(1-\frac{1}{N})}
\frac{1+\frac{c'}{N^2}}
{\left[1-\frac{c'}{N^2}-2\max_{z\in J} V^{\circ}(z,J)\right]}
\cr
}
\Eq(T.18)
$$
(which is meaningful whenever
$
1-\frac{c'}{N^2}-2\max_{z\in J} V^{\circ}(z,J)>0
$)
and, for all $\s$ such that $H(\s\cup A)$ is satisfied,
$$
\eqalign{
\E(\t^{\s}_{A})
\geq &  \frac{2^{N}}{(|J|-1)(1-\frac{1}{N})}
(1-\frac{c'}{N^2})
\Bigl[1-(1+O(\sfrac{1}{N}))V^{\circ}(x,J)\Bigr]
}
\Eq(T.19)
$$
whereas if
$H(\s\cup A)$ is not satisfied,
the term $1-\frac{1}{N}$ in
the r\.h\.s\. of \eqv(T.19) must be replaced by $1+O(\frac{1}{N})$. Since
$|J|=|\g(A\cup\s)|=|A|+1$
this proves the first assertion of the theorem.

We now turn to the second assertions of Theorem \thv(T.theo.1).
Here, we first use \eqv(T.6) to write
$$
\E\left(\t^{\s}_{\eta}\mid\t^{\s}_{\eta}<\t^{\s}_{A\setminus\eta}\right)
=
\frac{1}{\P(\t^{\s}_{A}<\t^{\s}_{\s})}
\left[1+
\sum_{\s'\in (A\cup\s)^c}
\P(\t^{\s'}_{\s}<\t^{\s'}_{A})
\frac{\P(\t^{\s'}_{\eta}<\t^{\s'}_{A\setminus\eta})}
{\P(\t^{\s}_{\eta}<\t^{\s}_{A\setminus\eta})}
\right]
\Eq(T.20)
$$
Using Lemma \thv(L.6) to pass to the lumped chain,
\eqv(T.20) becomes
$$
\E(\t^{\s}_{A})=\frac{1}{\P^{\circ}(\t^x_{J\setminus x}< \t^x_x)}
\left[1+\sum_{\s'\in (A\cup\s)^c}H^{\circ}_J(\g(\s'),x)
\frac{H^{\circ}_I(\g(\s'),y)}{H^{\circ}_I(x,y)}\right]
\Eq(T.21)
$$
where $x=\g(\s)$, $y=\g(\eta)$, $I=\g(A)$ and $J=\g(A\cup\s)$. Just as in the proof
of the first assertion the bounds \eqv(T.2) are obtained by inserting the estimates of
Theorem \thv(P.theo.0)
to express  $\P^{\circ}(\t^x_{J\setminus x}< \t^x_x)$,
and those of Theorem \thv(P.theo.00) to evaluate the sum
in the r\.h\.s\.
The only appreciable difference is that, in addition to terms of the form \eqv(T.100),
we now also have to deal with the terms
$$
\eqalign{
R'_1(x,y) & \equiv\frac{|J||I|}{2^{N}}
\sum_{\s'\in (A\cup\s)^c}\phi_x(\dist(\g(\s'),x))\phi_y(\dist(\g(\s'),y))
\cr
R'_2(J,I) & \equiv \frac{1}{2^{N}}\sum_{\s'\in (A\cup\s)^c}V^{\circ}(\g(\s'),J)V^{\circ}(\g(\s'),I)
}
\Eq(T.100bis)
$$
Note however that since $\phi_y(\dist(\g(\s'),y))\leq 1$ we easily get, proceeding
as we did to bound $R_1(x)$ and $R_2(J)$, that
$R'_1(x,y)\leq e^{-N/4}$ and $R'_2(J,I)\leq e^{-N/4}$. We leave the details to the reader.
This concludes the proofs of  Theorem \thv(T.theo.1) and thus, of Theorem \thv(T.theo.1.lumped).\endproof



\proofof{Theorem \thv(I.theo.3)} Theorem \thv(I.theo.3) is an immediate consequence
of Theorem \thv(T.theo.1) when replacing $\VV_{N,d}(A\cup\s)$ by $\UU_{N,d}(A\cup\s)$
in the latter. Note that the condition  $\UU_{N,d}(A\cup\s)\leq\frac{1}{4}$
guarantees that
$
c^-_N-2\UU_{N,d}(A\cup\s)\geq 1/3
$.
The constant $\frac{1}{4}$ has no special significance:
This choice is made for simplicity only.
\endproof

\proofof{Theorem \thv(T.theo.2)}
This Theorem is proven just as Lemma 3.1 of [BEGK1]. The idea is to
evaluate $\E\t^y_0$ using the lumped chain version of \eqv(T.5)
(see (3.12) in [BEGK1]). In the case $d>1$, the main difference between the proof of our bound \eqv(T.4) and
the bound (3.7) of [BEGK1] is that the bound (3.16) in the latter has here
to be replaced by the bound \thv(P.l7.5') from Lemma \thv(P.lemma.7), i\.e\.,
$$
\P^{\circ}(\t^x_0<\t^x_x)
\geq
\frac{c}{N}\left[\frac{1}{d}\sum_{\nu=1}^{d}\frac{1}{\sqrt{|\L_{\nu}|}}\right]^{-1}
\geq
\frac{c}{N}
\Eq(P.l7.5'')
$$
where  $0<c<\infty$ is a numerical constant. Proceeding as in the proof of Lemma 3.1 of [BEGK1]
we then get
$$
\textstyle{
\E\t^y_0 \leq CN^{2}\left(1+\sum_{x\in\G_{N,d}\setminus\{y,0\}} 1\right)
\leq CN^{2}|\G_{N,d}|= CN^{2}\prod_{k=1}^{d}|\L_k|\leq  CN^{d+2}
}
\Eq(P.l7.5''')
$$
The case $d=1$ is of course well known (see e\.g\. (4.34) page 28 in [K]).
Let us mention that the bound \eqv(T.4) can be obtained along the same
lines as above, but using the explicit one dimensional formula \eqv(A.25)
of Appendix A2 to evaluate carefully the right hand side of (the lumped chain version of) \eqv(T.5).
The bound \eqv(T.3') also results from a more careful evaluation of \eqv(T.5), using Lemma
\thv(A.lemma.2) from Appendix A1 to bound $\P^{\circ}(\t^x_0<\t^x_x)$ by a
sum of one dimensional quantities,
and using Lemma \thv(A.lemma.3) from Appendix A2 to bound each term of this sum.
Since the proof of this bound is a simple though lengthy procedure
we leave it out.
\endproof

\bigskip
\vfill\eject

\bigskip
\chap{7. Laplace transforms.}7

In this chapter we compute the Laplace transforms of hitting
times for the chain on the hypercube and prove Theorem \thv(I.theo.4nonsparse)
of Chapter 1 together with its corollary.
In the same spirit as for hitting times these results will
be deduced (in Section 7.3) from their lumped chain counterparts (proved
in Section 7.2). In the first section we collect the statements of the main
results for both chains.


\bigskip
\line{\bf 7.1. Statement of the main results.\hfill}

We will see that Theorem \thv(I.theo.4nonsparse) of Section 1 is a direct consequence of the
following result  for the lumped chain.


\theo{\TH(TLC.theo.1)}{\it  Let $d'\leq d_0(N)/2$ and assume that
$A\subset\SS_N$ is compatible with a partition $\L'$ into $d'$ classes.
Then for all $\s\notin A$ there exits a partition $\L$ into $d$ classes,
with $d'<d\leq 2d'$, compatible with $A\cup\s$.
Let one such partition be fixed. If there exists $0<\d<1$ such that
$$
\VV_{N,d}(A\cup\s)\leq\frac{\d}{4}
\Eq(TLC.t.0)
$$
then for all $\eta\in A$ the following holds:
for all $\e\geq\d$, there exists a constant $0<c_\e<\infty$
(independent of $\s, A, N$, and $d$) such that, for all $s$ real satisfying
$-\infty<s<1-\e$, for all $N$ large enough,
$$
\left|
\E\left(
e^{s{\t^{\s}_{A}}/\E\t^{\s}_{A}
}
\1_{\{\t^{\s}_{\eta}< \t^{\s}_{A\setminus\eta}\}}\right)
-\frac{1}{|A|}\frac{1}{1-s}
\right|\leq \frac{c_\e}{|A|}
\varepsilon_{N,d}(A,\eta,\s)
\Eq(TLC.t.1)
$$
where
$$
\varepsilon_{N,d}(A,\eta,\s)=\wt\varepsilon_{N,d}(A,\s)+|A|\phi_{\g(\eta)}(\Dist(\s,\eta))
\Eq(TLC.t.2)
$$
and
$$
\wt\varepsilon_{N,d}(A,\s)
=
\max\Big\{\VV_{N,d}(A\cup\s),\frac{1}{N^k}\Big\}
\Eq(TLC.2')
$$
for $k$ defined as in \eqv(I.25), i.e.
$$
k=
\cases 2, &\text{if} H(A\cup\s) \text{is satisfied}\cr
1, &\text{if} H(A\cup\s) \text{is not satisfied.}
\endcases
\Eq(I.25repetee)
$$

Moreover the above statement remains true with
$\phi_{\g(\eta)}(\Dist(\s,\eta))$ replaced by
$F(\Dist(\s,\eta))$ in \eqv(TLC.t.2), and
with $\VV_{N,d}(A\cup\s)$ replaced by $\UU_{N,d}(A\cup\s)$
in \eqv(TLC.t.0) and \eqv(TLC.2')
(see \eqv(I.12repetee) and \eqv(I.12new)).
}


As was already proved by  Matthews [M1] (see (3.5) and (3.6) p.138) a sharper result
can be obtained in the special case where $A$ consists of a single point.
This result, which we state below for the sake of completeness,
will be derived from our more general Theorem \thv(TL.theo.2).

\theo{\TH(TLC.theo.2)}{\it For any pair of distinct points $\s,\eta\in\SS_N$
the following holds: for all $\e>0$, there exists a constant $0<c_\e<\infty$
(independent of $\s, \eta$, and $N$) such that, for all $s$ real satisfying
$-\infty<s<1-\e$, for all $N$ large enough,
$$
\left|\E\left(e^{s{\t^{\s}_{\eta}}/2^N}\right)-
\frac{1-s\frac{1}{N}}{1-s(1+\frac{1}{N})}\right|
\leq
\frac{c_{\e}}{N^2}\,\text{if} \Dist(\eta,\s)=1
\Eq(TLC.t2.1)
$$
and
$$
\left|\E\left(e^{s{\t^{\s}_{\eta}}/2^N}\right)-
\frac{1}{1-s(1+\frac{1}{N})}\right|
\leq
\frac{c_{\e}}{N^2}\,\text{if} \Dist(\eta,\s)>1
\Eq(TLC.t2.2)
$$
}


Our next corollary states two key consequences of  Theorem \thv(TLC.theo.1).

\corollary{\TH(TLC.cor.1)}{\it
Under the assumptions and with the notation of Theorem \thv(TLC.theo.1),
the following holds:

\item{i)} For all $\e>\d$ there exists a constant $0<c_\e<\infty$ such that,
for all $s$ real satisfying $-\infty<s<1-\e$ and all $N$ large enough we have
$$
\left|
\E\left(e^{s{\t^{\s}_{A}}/\E\t^{\s}_{A}}\right)
-\frac{1}{1-s}
\right|\leq c_\e\wt\varepsilon_{N,d}(A,\s)
\Eq(TLC.2)
$$
If $\VV_{N,d}(A\cup\s)\rightarrow 0$ as $N\rightarrow\infty$
this implies that
${\t^{\s}_{A}}/\E\t^{\s}_{A}$
converges in distribution to an exponential random variable of mean value one.

\item{ii)} Let $A_1, A_2,\dots,A_n$ be a finite collection of non empty disjoint subsets of
$A$.
Then, for all $\e>\d$, for all $s_i$ real
satisfying $-\infty<s_i<1-\e$, $1\leq i\leq n$, and all $N$ large enough,
$$
\left|
\E\left(e^{\sum_{i=1}^n s_i{\t^{\s}_{A_i}}/\E\t^{\s}_{A_i}}\right)
-\prod_{i=1}^n\left(\E e^{s_i{\t^{\s}_{A_i}}/\E\t^{\s}_{A_i}}\right)
\right|\leq c_{n,\e}\wt\varepsilon_{N,d}(A,\s)
\Eq(TLC.3)
$$
for some constant $0<c_{n,\e}<\infty$.
Thus,
If $\VV_{N,d}(A\cup\s)\rightarrow 0$ as $N\rightarrow\infty$,
the random variables $({\t^{\s}_{A_i}},1\leq i\leq n)$ become
asymptotically independent in the limit. 
}


As we will prove in Section 7.3,
Theorem \eqv(TLC.theo.1), Theorem \eqv(TLC.theo.2), and Corollary \thv(TLC.cor.1)
are direct consequences  of their lumped chain counterparts, namely,
Theorem \eqv(TL.theo.1), Theorem \eqv(TL.theo.2), and
Corollary \thv(TL.cor.1), which we now state.

\theo{\TH(TL.theo.1)}{\it  Let $d\leq d_0(N)$ and
let $\g$ be any $d$-lumping (or equivalently let $\L$ be any
$d$-partition).
%
Let  $I\subset\SS_d$ and $y\in\SS_d\setminus I$ be such that,
for some $0<\d<1$,
$$
\VV^{\circ}_{N,d}(I\cup y)\leq\frac{\d}{4}
\Eq(TL.t.0)
$$
Then, for all $x\in I$
the following holds: for all $\e>\d$, there exists a constant $0<c_\e<\infty$
(independent of $y, I, N$, and $d$) such that, for all $s$ real satisfying
$-\infty<s<1-\e$, for all $N$ large enough,
$$
\left|
\E^{\circ}\left(
e^{s{\t^{y}_{I}}/
\E^{\circ}\t^{y}_{I}
}
\1_{\{\t^{y}_{x}< \t^{y}_{I\setminus x}\}}\right)
-\frac{1}{|I|}\frac{1}{1-s}
\right|\leq \frac{c_\e}{|I|}
\varepsilon^{\circ}_{N,d}(I,x,y)
\Eq(TL.t.1)
$$
where
$$
\varepsilon^{\circ}_{N,d}(I,x,y)=\wt\varepsilon^{\circ}_{N,d}(I,y)+|I|\phi_x(\dist(x,y))
\Eq(TL.t.2)
$$
and
$$
\wt\varepsilon^{\circ}_{N,d}(I,y)
=
\max\Big\{\VV^{\circ}_{N,d}(I\cup y),\frac{1}{N^k}\Big\}
\Eq(TL.2')
$$
where
$$
k=
\cases
2, &\text{if} H^{\circ}(I\cup y) \text{is satisfied}\cr
1,&\text{if} H^{\circ}(I\cup y)  \text{is not satisfied}\cr
\endcases
\Eq(TL.2'')
$$

Moreover the above statement remains true with $\phi_x(\dist(x,y))$ replaced by
$F(\dist(x,y))$ in \eqv(TL.t.2), and
with $\VV^{\circ}_{N,d}(I\cup y)$ replaced by $\UU^{\circ}_{N,d}(I\cup y)$ in
\eqv(TL.t.0) and \eqv(TL.2') (see \eqv(I.12ter) and \eqv(I.12bis)).
}

\remark If
$\VV^{\circ}_{N,d}(I)=o(1)$ and $|I|\max_{z\in I}\phi_z(\dist(z,y))=o(1)$
then \eqv(TL.t.0) holds true with $\d\equiv\d(N)=o(1)$.
Moreover, by \eqv(TL.l4.1) of Lemma \thv(TL.lemma.4),
$\VV^{\circ}_{N,d}(I)=o(1)$  whenever
$\UU^{\circ}_{N,d}(I)=o(1)$.


As announced earlier \eqv(TL.t.1) can be (partially) improved when
$I$ consists of a single point.
Theorem \thv(TL.theo.2) can be seen as a
$d$-dimensional lumped version of the result obtained by Matthews [M1]
for the chain on the hypercube (see Theorem \thv(TLC.theo.2)).

\theo{\TH(TL.theo.2)}{\it Assume that $d^2=O(N)$. Let $x\in\SS_d$
and
$y\in\G_{N,d}\setminus x$.
Then, for all $\e>0$, there exists a constant $0<c_\e<\infty$
(independent of $y, x, N$, and $d$) such that, for all $s$ real satisfying
$-\infty<s<1-\e$, for all $N$ large enough,
$$
\left|\E^{\circ}\left(e^{s{\t^{y}_{x}}/2^N}\right)-
\frac{1-s\frac{1}{N}}{1-s(1+\frac{1}{N})}\right|
\leq
c_{\e}\frac{d}{N^2}\,\text{if} \dist(x,y)=1
\Eq(TL.t2.1)
$$
and
$$
\left|\E^{\circ}\left(e^{s{\t^{y}_{x}}/2^N}\right)-
\frac{1}{1-s(1+\frac{1}{N})}\right|
\leq
\frac{c_{\e}}{N^2}\,\text{if} \dist(x,y)>1
\Eq(TL.t2.2)
$$
}

\remark Note that Theorem \thv(TL.theo.2) is valid not only for
$y\in\SS_d\setminus x$ but for all $y\in\G_{N,d}\setminus x$.
When $y\in\SS_d\setminus x$ then \eqv(TL.t2.1) and \eqv(TL.t2.2)
simply are  reformulations of \eqv(TLC.t2.1) and \eqv(TLC.t2.2).



As a corollary to Theorem \eqv(TL.theo.1), we have:

\corollary{\TH(TL.cor.1)}{\it
Under the assumptions and with the notation of Theorem \thv(TL.theo.1),
the following holds:

\item{i)} For all $\e>\d$ there exists a constant $0<c_\e<\infty$ such that,
for all $s$ real satisfying $-\infty<s<1-\e$ and all $N$ large enough we have
$$
\left|
\E^{\circ}\left(e^{s{\t^{y}_{I}}/\E^{\circ}\t^{y}_{I}}\right)
-\frac{1}{1-s}
\right|\leq c_\e\wt\varepsilon^{\circ}_{N,d}(I,y)
\Eq(TL.2)
$$
%
If $\VV^{\circ}_{N,d}(I\cup y)\rightarrow 0$ as $N\rightarrow\infty$
this implies that
${\t^{y}_{I}}/\E^{\circ}\t^{y}_{I}$
converges in distribution to an exponential random variable of mean value one.

\item{ii)} Let $I_1, I_2,\dots,I_n$ be a finite collection of non empty disjoint subsets of
$I$.
%
%
Then, for all $\e>\d$, for all $s_i$ real
satisfying $-\infty<s_i<1-\e$, $1\leq i\leq n$, and all $N$ large enough,
$$
\left|
\E^{\circ}\left(e^{\sum_{i=1}^n s_i{\t^{y}_{I_i}}/\E^{\circ}\t^{y}_{I_i}}\right)
-\prod_{i=1}^n\left(\E^{\circ} e^{s_i{\t^{y}_{I_i}}/\E^{\circ}\t^{y}_{I_i}}\right)
\right|\leq c_{n,\e}\wt\varepsilon^{\circ}_{N,d}(I,y)
\Eq(TL.3)
$$
for some constant $0<c_{n,\e}<\infty$.
Thus, if $\VV^{\circ}_{N,d}(I\cup y)\rightarrow 0$ as $N\rightarrow\infty$,
the random variables $({\t^{y}_{I_i}},1\leq i\leq n)$ become asymptotically
independent in the limit.
}


The rest of this section is organized as follows. We will
first show how Theorem \thv(TL.theo.1) implies Corollary
\thv(TL.cor.1); doing this will explain
the role and usefulness of the special form of the Laplace transform
appearing in \eqv(TL.t.1). Theorem \thv(TL.theo.1) and Theorem \thv(TL.theo.2)
are themselves  specializations of a more general results, namely
Proposition \thv(TL.prop.1) and Corollary \thv(TL.cor.3), which we next
state and prove. Lastly, we prove Theorem \thv(TL.theo.1) and Theorem \thv(TL.theo.2).

\smallskip
\bigskip
\line{\bf 7.2.  Laplace transforms of hitting times for the lumped chain.\hfill}

Let us fix the notation for the Laplace transforms
of interest. If $I$ and $J$ are disjoint subsets of
$\G_{{N},d}$, and if $y$ is any point in $\G_{{N},d}$
(we include the possibility that $y\in I\cup J$), we define
$$
G^y_{I}(u)\equiv \E^{\circ} e^{u\t^y_I}, \quad
G^y_{I,J}(u)\equiv \E^{\circ} e^{u\t^y_I}\1_{\{\t^y_I<\t^y_J\}}
\Eq(TL.4)
$$
for $u\in D\subset\C$, where $D$ is chosen in a such a way that
the right hand sides of \eqv(TL.4) exist.
Note that
$$
G^y_{I}(u)=\sum_{x\in I}G^y_{x,I\setminus x}(u)
\Eq(TL.5)
$$
(which of course is useful only when $I$ does not consist of a single point)
and
$$
G^y_{I,J}(u)=\sum_{x\in I}G^y_{x, (I\setminus x)\cup J}(u)
\Eq(TL.6)
$$
The study of the Laplace transforms \eqv(TL.4) thus reduces to that of the
basic quantities
$$
G^y_{x, J}(u),\text{for} J\subset\G_{{N},d},\, x\in\G_{{N},d}\setminus J,  \text{and}
y\in\G_{{N},d}
\Eq(TL.7)
$$
to which we must add
$$
G^y_{x}(u), \text{for}\, x\in\G_{{N},d}  \text{and}
y\in\G_{{N},d}
\Eq(TL.8)
$$
if we want to cover the case where $I$ consists of a single point.




\proofof{Corollary \thv(TL.cor.1)} Note that for $\varepsilon^{\circ}_{N,d}(I,x,y)$
and $\wt\varepsilon^{\circ}_{N,d}(I,y)$ defined in \eqv(TL.t.2) and \eqv(TL.2'),
by \eqv(N.3bis) and \eqv(I.12bis),
$$
\frac{1}{|I|}\sum_{x\in I}\varepsilon^{\circ}_{N,d}(I,x,y)
=\wt\varepsilon^{\circ}_{N,d}(I,y)+V^{\circ}_{N,d}(y,I)
\leq C\wt\varepsilon^{\circ}_{N,d}(I,y)
\Eq(TL.9)
$$
for some positive finite constant $C$.
Thus, \eqv(TL.2) of assertion (i) is a direct
consequence of \eqv(TL.t.1) and \eqv(TL.5); the fact that it implies convergence
in distribution when $\wt\varepsilon^{\circ}_{N,d}(I,y)=o(1)$ is a classical result
(see e\.g\. [Fe], Chapter XIII, Section 1, Theorem 2).
Let us turn to assertion (ii). In what follows $I_1, I_2,\dots,I_n$ is a finite
collection of non-empty disjoint subsets of $I$ of cardinality $|I_i|=M_i$,
and we assume that  $-\infty<s_i<1-\e$, $1\leq i\leq n$. Let us observe that,
for all $z\in I\cup y$,
$$
\left|
\prod_{i=1}^n G^{z}_{I_i}\Big(t_i(1-\sfrac{1}{N})/2^N\Big)
-\prod_{i=1}^n (1-\sfrac{t_i}{M_i})^{-1}
\right|
\leq c_{n,\e}\wt\varepsilon^{\circ}_{N,d}(I,z)
\Eq(TL.70)
$$
for some constant $0<c_{n,\e}<\infty$.
Indeed, by \eqv(TL.2), since
$\sum_i\wt\varepsilon^{\circ}_{N,d}(I_i,z)\leq \wt\varepsilon^{\circ}_{N,d}(\cup_i I_i,z)$,
$$
\left|
\prod_{i=1}^n G^{z}_{I_i}(s_i/\E^{\circ}\t^{z}_{I_i})
-
\prod_{i=1}^n (1-s_i)^{-1}
\right|
\leq c'_{n,\e}\wt\varepsilon^{\circ}_{N,d}(I,z)
\Eq(TL.70')
$$
for some $0<c'_{n,\e}<\infty$. Since the underlying $d$-lumping $\g$
is assumed to be generated by a log-regular $d$-partition, we may use
Theorem \thv(T.theo.1.lumped) to write the quantities
$\E^{\circ}\t^{z}_{I_i}$ in the form
$$
\E^{\circ}\t^{z}_{I_i}=\frac{2^{N}}{M_i}
\Big(1+\frac{1}{N}\Big)(1+O(\wt\varepsilon^{\circ}_{N,d}(I_i,z)))
\Eq(TL.71)
$$
where as before $\wt\varepsilon^{\circ}_{N,d}(I_i,z)$ is given by \eqv(TL.2').
Then, making the change of variable
$
t_i=s_iM_i(1+O(\wt\varepsilon^{\circ}_{N,d}(I_i,z)))
$,
\eqv(TL.70') yields \thv(TL.70)
(recall that by assumption $\VV^{\circ}_{N,d}(I\cup y)\leq\d/4$, and this implies that
$\VV^{\circ}_{N,d}(I_i\cup z)\leq\d/4$ for all $1\leq i\leq n$; consequently,
$\wt\varepsilon^{\circ}_{N,d}(I_i,z)\leq\d/2$ for all $1\leq i\leq n$,
and since $\e\geq\d$, this guarantees that
$\sfrac{t_i}{M_i}<(1-\e)(1+\d/2)<
(1-\e)(1+\e/2)<1-\e'$ for some $0<\e'<\e<1$).
As a consequence, \eqv(TL.3) is equivalent to
$$
\left|
\E^{\circ}\left(e^{\sum_{i=1}^n t_i{\t^{y}_{I_i}}(1-\frac{1}{N})/2^N}\right)
-\prod_{i=1}^n (1-\sfrac{t_i}{M_i})^{-1}
\right|\leq c''_{n,\e}\wt\varepsilon^{\circ}_{N,d}(I,y)
\Eq(TL.3bis)
$$
for some constant $0<c''_{n,\e}<\infty$.

We now proceed to prove \eqv(TL.3bis) using an inductive argument.
To start the induction just observe that, by \eqv(TL.70),
\eqv(TL.3bis)
is true when the collection $I_1, I_2,\dots,I_n$
is reduced to just one of its elements; more precisely,
for each $1\leq i\leq n$ and arbitrary $z\in I\cup y$
$$
\left|
\E^{\circ}\left(e^{t_i{\t^{z}_{I_i}}(1-\frac{1}{N})/2^N}\right)-
(1-\sfrac{t_i}{M_i})^{-1}
\right|
\leq c_{1,\e}\wt\varepsilon^{\circ}_{N,d}(I_i,z)
\Eq(TL.72)
$$
Let now $1<m\leq n$ and choose $m$ elements in the collection
$I_1, I_2,\dots,I_n$; without loss of generality we may take
$I_1, I_2,\dots,I_m$. We will next establish that, if for each
$1\leq j\leq m$ and any $z\in (\cup_{{i=1};i\neq j}^m I_i)\cup y$
$$
\left|
\E^{\circ}\left(e^{\sum_{{i=1}\atop{i\neq j}}^m t_i{\t^{z}_{I_i}}(1-\frac{1}{N})/2^N}\right)
-\prod_{{i=1}\atop{i\neq j}}^m
(1-\sfrac{t_i}{M_i})^{-1}
\right|\leq c_{m-1,\e}\wt\varepsilon^{\circ}_{N,d}(\cup_{{i=1}\atop{i\neq j}}^m I_i,z)
\Eq(TL.73)
$$
holds true for some $0<c_{m-1,\e}<\infty$, then \eqv(TL.3bis)
holds true with $n=m$. To do this we set
$B_m=\cup_{j=1}^m I_i$ and write
$$
\E^{\circ}\left(e^{\sum_{i=1}^m t_i{\t^{y}_{I_i}}(1-\frac{1}{N})/2^N}\right)
=\sum_{j=1}^m\sum_{x\in I_j}\E^{\circ}\left(e^{\sum_{i=1}^m t_i{\t^{y}_{I_i}}(1-\frac{1}{N})/2^N}
\1_{\{\t^{y}_{x}=\t^{y}_{B_m\setminus x}\}}\right)
\Eq(TL.74)
$$
Next, for each $1\leq j\leq m$ and each $x\in I_j$,
$$
\eqalign{
&\E^{\circ}\left(e^{\sum_{i=1}^m t_i{\t^{y}_{I_i}}(1-\frac{1}{N})/2^N}
\1_{\{\t^{y}_{x}=\t^{y}_{B_m\setminus x}\}}\right)
\cr
=&
\E^{\circ}\left(e^{
t_j\t^{y}_{x}(1-\frac{1}{N})/2^N
+\sum_{{i=1}\atop{i\neq j}}^m t_i{\t^{y}_{I_i}}(1-\frac{1}{N})/2^N
}
\1_{\{\t^{y}_{x}=\t^{y}_{B_m\setminus x}\}}
\1_{\{\t^{y}_{x}=\t^{y}_{I_j\setminus x}\}}\right)
\cr
=&
\E^{\circ}\left(e^{
(\sum_{i=1}^m t_i){\t^{y}_{x}}(1-\frac{1}{N})/2^N
+\sum_{{i=1}\atop{i\neq j}}^m t_i{\t^{x}_{I_i}}(1-\frac{1}{N})/2^N
}
\1_{\{\t^{y}_{x}=\t^{y}_{B_m\setminus x}\}}
\1_{\{\t^{y}_{x}=\t^{y}_{I_i\setminus x}\}}\right)
\cr
=&
\E^{\circ}\left(e^{
(\sum_{i=1}^m t_i){\t^{y}_{B_m}}(1-\frac{1}{N})/2^N
+\sum_{{i=1}\atop{i\neq j}}^m t_i{\t^{x}_{I_i}}(1-\frac{1}{N})/2^N
}
\1_{\{\t^{y}_{x}=\t^{y}_{B_m\setminus x}\}}\right)
\cr
=&
G^{y}_{x, B_m\setminus x}
\left(\t^{y}_{B_m}\sum_{i=1}^m t_i(1-\sfrac{1}{N})/2^N\right)
\E^{\circ}\left(
e^{\sum_{{i=1}\atop{i\neq j}}^m t_i{\t^{x}_{I_i}}(1-\frac{1}{N})/2^N}
\right)
\cr
}
\Eq(TL.75)
$$
Setting $\overline M=\sum_{i=1}^m M_i$, we then define
$$
\eqalign{
V_m(x)&\equiv
G^{y}_{x, B_m\setminus x}
\left(\t^{y}_{B_m}\sum_{i=1}^m t_i(1-\sfrac{1}{N})/2^N\right)
-\sfrac{1}{\overline M}
(1-\sfrac{1}{\overline M}
{\scriptstyle{
\sum_{i=1}^m
}}
t_i)^{-1}
\cr
W_m(x)&\equiv
\E^{\circ}\left(e^{\sum_{{i=1}\atop{i\neq j}}^m t_i{\t^{x}_{I_i}}(1-\frac{1}{N})/2^N}\right)
-\prod_{{i=1}\atop{i\neq j}}^m (1-\sfrac{t_i}{M_i})^{-1}
\cr
}
\Eq(TL.76)
$$
and rewrite \eqv(TL.75) as
$$
\E^{\circ}\left(e^{\sum_{i=1}^m t_i{\t^{y}_{I_i}}(1-\frac{1}{N})/2^N}
\1_{\{\t^{y}_{x}=\t^{y}_{B_m\setminus x}\}}\right)
=
\sfrac{1}{\overline M}(1-\sfrac{1}{\overline M}
{\scriptstyle{\sum_{i=1}^m }}t_i)^{-1}
\prod_{{i=1}\atop{i\neq j}}^m (1-\sfrac{t_i}{M_i})^{-1}+\RR(x)
\Eq(TL.77)
$$
where
$$
\RR(x)\equiv
W_m(x)\sfrac{1}{\overline M}(1-\sfrac{1}{\overline M}
{\scriptstyle{\sum_{i=1}^m }}t_i)^{-1}
+V_m(x) \prod_{{i=1}\atop{i\neq j}}^m (1-\sfrac{t_i}{M_i})^{-1}
+V_m(x)W_m(x)
\Eq(TL.78)
$$
Of course we want to make use of \eqv(TL.77) in \eqv(TL.74):
observing first that
$$
\eqalign{
\sum_{j=1}^m\sum_{x\in I_j}
\sfrac{1}{\overline M}(1-\sfrac{1}{\overline M}
{\scriptstyle{\sum_{i=1}^m }}t_i)^{-1}
\prod_{{i=1}\atop{i\neq j}}^m (1-\sfrac{t_i}{M_i})^{-1}
=&\prod_{{i=1}}^m (1-\sfrac{t_i}{M_i})^{-1}
\cr
}
\Eq(TL.79)
$$
we arrive at
$$
\E^{\circ}\left(e^{\sum_{i=1}^m t_i{\t^{y}_{I_i}}(1-\frac{1}{N})/2^N}\right)
=\prod_{{i=1}}^m (1-\sfrac{t_i}{M_i})^{-1}
+\sum_{j=1}^m\sum_{x\in I_j}\RR(x)
\Eq(TL.80)
$$
and it remains to bound the sum appearing in the r\.h\.s\..
By \eqv(TL.t.1) and an appropriate change of variable,
$
\left|V_m(x)\right|\leq \frac{1}{\overline M}
\tilde c_{\e}\varepsilon^{\circ}_{N,d}(B_m,x,y)
$
and, reasoning as in the proof of \eqv(TL.2),
$
\sum_{j=1}^m\sum_{x\in I_j}\left|V_m(x)\right|
\leq \frac{1}{\overline M}\tilde c_{\e}\wt\varepsilon^{\circ}_{N,d}(B_m,y)
$.
Next, by \eqv(TL.73), for $x\in I_j$,
$
\left|W_m(x)\right|\leq
\hat c_{m-1,\e}\wt\varepsilon^{\circ}_{N,d}(B_m\setminus I_j,x)
$
and
$
\frac{1}{\overline M}\sum_{j=1}^m\sum_{x\in I_j}\left|W_m(x)\right|
\leq
\frac{1}{\overline M}\hat c_{m-1,\e}\wt\varepsilon^{\circ}_{N,d}(B_m\setminus I_j,x)
\leq
\hat c_{m,\e}\wt\varepsilon^{\circ}_{N,d}(B_m,y)
$
In this way, one easily checks that
$
\left|\sum_{j=1}^m\sum_{x\in I_j}\RR(x)\right|
\leq c''_{m,\e}\wt\varepsilon^{\circ}_{N,d}(B_m,y)
$
(all the constants $\tilde c_{\e}, \hat c_{m,\e}, c''_{m,\e}$
above being positive and finite.) Now, this proves \eqv(TL.3bis) with $n=m$.
Note that we can prove in exactly the same way that \eqv(TL.3bis) holds with
$n=m$ and $x$ replaced by any $z\in B_m\cup y$. This completes the
inductive argument that started in \eqv(TL.73),
and concludes the proof of assertion (ii) of Corollary \thv(TL.cor.1).
\endproof



We now turn to the proof of Theorem \thv(TL.theo.1). It will heavily rely on a
detailed analysis of the basic Laplace Transforms $G^y_{x,J}(u)$ introduced
in \eqv(TL.7).
We summarize the results of this analysis in Proposition \thv(TL.prop.1), which
we now state. We will then immediately proceed to its proof,
give next the proofs of Theorem \thv(TL.theo.1) and  Theorem \thv(TL.theo.2),
and close this section with the proofs of
Theorem  \thv(I.theo.4) and  Theorem \thv(I.theo.4nonsparse) and
Corollary \thv(I.cor.theo.4) of Chapter 1.


\proposition{\TH(TL.prop.1)}{\it Let $d\leq d_0(N)$ and
let $\g$ be any $d$-lumping (or equivalently let $\L$ be any
$d$-partition).
Let $J\subset\SS_d$
and $x\in\SS_d\setminus J$ be such that, for some $0<\d<1$,
$$
\VV^{\circ}_{N,d}(J\cup x)\leq\frac{\d}{4}
\Eq(TL.12bis)
$$
Then,  for all $y\in\G_{{N},d}$, the following holds.
Set
$$
\barunder u(d)^{-1}\equiv \wh\Th^2(d)/\E^{\circ} \t^0_0\,,\quad
\bar u^{-1}\equiv\frac{2^N}{|J\cup x|}
\Bigl(1+\frac{1}{N}\Bigr)\,,
\Eq(TL.12)
$$
where $\wh\Th(d)$ was defined in \eqv(T.4bis), and define
$$
s(u)= u/\bar u 
\Eq(TL.13)
$$

\item{i)}  For all $u$ real satisfying $-\rho\barunder u(d)<u< \bar u$ for some
$0<\rho<1$, we have:

\itemitem{}if $J\neq\emptyset$,
$$
\eqalign{
G^y_{x,J}(u)
&=
\R^\circ\left(\t^y_x<\t^y_{J}\right)
\frac{1}{1-s(u)}
+
\P^{\circ}\left(\t^{y}_{x}<\t^{y}_{J\cup 0}\right)
\frac{-s(u)}{1-s(u)}
+
\RR_{0}(u)
}
\Eq(TL.14)
$$
\itemitem{}if  $J=\emptyset$,
$$
\eqalign{
G^y_{x}(u)
&=
\frac{1}{1-s(u)}
+
\P^{\circ}\left(\t^{y}_{x}<\t^{y}_{0}\right)
\frac{-s(u)}{1-s(u)}
+
\RR_{\emptyset}(u)
}
\Eq(TL.15)
$$
\item{} where
$$
\eqalign{
\RR_0(u)&=\frac{\P^{\circ}\left(\t^{0}_{x}<\t^{0}_{J}\right)}{1-s(u)}
\left[\RR_1(u)+\P^{\circ}\left(\t^{y}_{0}<\t^{y}_{J\cup x}\right)\RR_2(u)\right]+\RR_3(u)
\cr
\RR_{\emptyset}(u)&=\frac{1}{1-s(u)}
\left[\RR_1(u)+\P^{\circ}\left(\t^{y}_{0}<\t^{y}_{x}\right)\RR_2(u)\right]+\RR_3(u)
\cr
}
\Eq(TL.16)
$$
and (uniformly in $x,y$ and $J$)
$$
\eqalign{
\RR_1(u)&=O(|u|\wh\Th)\cr
\RR_3(u)&=O(|u|\wh\Th)\cr
\RR_2(u)&=O\left(\max\left\{
\VV^{\circ}_{N,d}(J\cup x)\left|\frac{-s(u)}{1-s(u)}\right|,\frac{1}{N^2}\left|\frac{-s(u)}{1-s(u)}\right|,\frac{|u|}{\barunder u(d)}
\right\}\right)\cr
}
\Eq(TL.17)
$$

\item{ii)}  Let $\ell_y=\dist(0,y)$
For all $u$ real satisfying $u<-\rho\barunder u(d)$ for some
$0<\rho<1/9$,
$$
G^y_{x,J}(u)\leq
G^y_{x,J\cup 0}(u)+
|J\cup x|\frac{2^{-N+2}e^{-|u|(\ell_x+\ell_y)}}{\barunder u(d)\rho(1-9\rho)}
\Bigl(1-\frac{1}{N}+3\max\Big\{\VV^{\circ}_{N,d}(J\cup x),\frac{4}{N^2}\Big\}\Bigr)
\Eq(TL.19)
$$
}


\remark Since $x\in \SS_d$, $ \ell_x=N/2$ (where we assumed to simplify that
$|\L_k|$ is even for all $1\leq k\leq d$). Thus for $-u$ large enough,
more precisely for $-u$ such that $e^{-|u|N/2}/\barunder{\it u}(d)=o(1)$,
the  coefficient of $2^{-N}$ in \eqv(TL.19) tends to zero, and thus the
second term of \eqv(TL.19) decays faster than $2^{-N}$.

\remark One might expect that $\barunder {\it u}(d)\sim 1/\E^{\circ} \t^0_0$, or
at least $1/\wh\Th(d)$. We are however not able to prove this. This is
due to rather coarse estimates on $G^0_0(u)$ for $u>0$.

\remark The only place where we will make  use of condition \eqv(TL.12bis)
is \eqv(TL.28).
It is used to ensure that $1/G^x_{0,x\cup J}(0)=O(1)$ in \eqv(TL.23).
(We see that though \eqv(TL.12bis) will do no arm we could have asked less.)


Our aim in  Proposition \thv(TL.prop.1) was to make statements that are valid without assumptions on $y$.
This explains the special form of \eqv(TL.14) and \eqv(TL.15), where we kept the term
$\P^{\circ}\left(\t^{y}_{x}<\t^{y}_{J\cup 0}\right)$  explicit.
This enables us in particular to deduce the following result,
which is tantamount to the statement of Theorem \thv(TL.theo.2).

\corollary{\TH(TL.cor.3)}{\it Assume that $d^2=O(N)$. Let $x\in\SS_d\setminus J$
and $y\in\G_{N,d}\setminus x$. Then
$$
\VV^{\circ}_{N,d}(y\cup x)\leq\frac{1}{N}\left(1+O\Bigl(\frac{d}{N}\Bigr)\right)
$$
Moreover, with the notation of Proposition \thv(TL.prop.1),
\item{i)}  For all $u$ real satisfying $-\rho\barunder u(d)<u< \bar u$ for some
$0<\rho<1$, we have:
$$
G^y_{x}(u)
=
\frac{1-\frac{1}{N}(1+O(\frac{d}{N}))s(u)}{1-s(u)}+\RR_{\emptyset}(u)\,\text{if} \dist(x,y)=1
\Eq(TL.c3.1)
$$
and
$$
G^y_{x}(u)
=
\frac{1}{1-s(u)}+\RR_{\emptyset}(u)\,\text{if} \dist(x,y)>1
\Eq(TL.c3.1')
$$
where
$$
\RR_{\emptyset}(u)=O\left(\max\left\{
\frac{1}{N^2}\left|\frac{-s(u)}{1-s(u)}\right|,\frac{|u|}{\barunder u(d)}, \frac{2-s(u)}{1-s(u)}|u|\wh\Th
\right\}\right)
\Eq(TL.c3.2)
$$
\item{ii)}
For all $u$ real satisfying $u<-\rho\barunder u(d)$ for some
$0<\rho<1/9$,
$$
G^y_{x}(u)\leq
G^y_{x,0}(u)+
|J\cup x|\frac{2^{-N+2}e^{-|u|(\ell_x+\ell_y)}}{\barunder u(d)\rho(1-9\rho)}
\Bigl(1+\frac{3}{N}\Bigr)
\Eq(TL.c3.3)
$$
}

\proofof{Corollary \thv(TL.cor.3)} Note that when $\dist(x,y)=1$,
$\P^{\circ}\left(\t^{y}_{x}<\t^{y}_{0}\right)\geq r_N(y,x)=\frac{1}{N}$.
Together with the upper bound of Theorem  \thv(P.prop.2) and
Lemma \thv(A3.lemma.1) of Appendix A3, under the assumption that $d^2=O(N)$, this yields
$\P^{\circ}\left(\t^{y}_{x}<\t^{y}_{0}\right)=\frac{1}{N}(1+O(\frac{d}{N}))$
which in turn implies that
$
\VV^{\circ}_{N,d}(J\cup x)\leq \frac{1}{N}(1+O(\frac{d}{N}))
$
.
Corollary \thv(TL.cor.3) is now an immediate consequence of \eqv(TL.15)
and \eqv(TL.19) of Proposition \thv(TL.prop.1). \endproof

\proofof{Proposition  \thv(TL.prop.1)}
It is rather simple to see that the minimal eigenvalue of  the generator
$1-P_N$ of the simple random walk with Dirichlet boundary conditions on a
finite set of points is of the order of $2^{-N}$; thus the Laplace transforms
$G^y_{x,J}(u)$  defined in \eqv(TL.7) will have poles at distance $2^{-N}$
from zero on the positive real axis. This makes it rather hard analytically to
get precise information on their behavior near zero directly via e.g. expansions.
On the other hand, if we consider the generator of the lumped chain with Dirichlet
conditions at zero, it turns out that the minimal eigenvalue is polynomial in $N$,
so that the corresponding Laplace transforms have their first pole much farther away
from zero. Thus our strategy will be to decompose all processes
at visits at zero, and to express the full Laplace transforms as functions
of Laplace transforms of processes that are killed at zero.
In practice, this yields:
$$
G^y_{x,J}(u)=G^y_{x,J\cup 0}(u)+G^y_{0,J\cup x}(u)
\frac{G^0_{x,J\cup 0}(u)}{1-G^0_{0,J\cup x}(u)}
\Eq(TL.10)
$$
and
$$
G^y_{x}(u)=G^y_{x, 0}(u)+G^y_{0,x}(u)
\frac{G^0_{x, 0}(u)}{1-G^0_{0,x}(u)}
\Eq(TL.11)
$$
To prove Proposition \thv(TL.prop.1) we will estimate each of the Laplace Transforms
appearing in the right hand side of \eqv(TL.10) and \eqv(TL.11) separately.
We begin by the proof of assertion (i) for $J\neq\emptyset$.
Our starting point here is \eqv(TL.10). Using reversibility,
\eqv(TL.10) can be rewritten as
$$
G^y_{x,J}(u)=G^y_{x,J\cup 0}(u)+G^y_{0,J\cup x}(u)
\frac{G^x_{0,J\cup x}(u)}
{\frac{\Q(0)}{\Q(x)}\left(1-G^0_{0,J\cup x}(u)\right)}
\Eq(TL.20)
$$

Let us first consider the second term in the r\.h\.s\. of \eqv(TL.20).
We call this term $h(u)$.

\lemma{\TH(TL.lemma.3)}{\it Under the assumptions of Proposition \thv(TL.prop.1),
for all $u$ satisfying $-\rho\barunder{\it u}(d)<u< \bar u$ for some $0<\rho<1$,
$$
\eqalign{
h(u)&=\frac{g}{1-s(u)}(f+\wt\RR_1(u))\left(1+\wt\RR_2(u)\right)
\cr
}
\Eq(TL.48)
$$
where $f$ and $g$ can be written as
$$
\eqalign{
f
=&\R\left(\t^y_0<\t^y_{x\cup J}\right)
=1-\frac{1}{N}+Z_{N,d}(J\cup x)
\cr
g
=&\R\left(\t^0_x<\t^0_{J}\right)
=\frac{1}{|J\cup x|}\Bigl(1+Z_{N,d}(J\cup x)
\Bigr)
}
\Eq(TL.l3.1)
$$
where $Z_{N,d}(J\cup x)$ obeys the bound
$$
\left|Z_{N,d}(J\cup x)\right|< 3\max\Big\{\VV^{\circ}_{N,d}(J\cup x),\frac{4}{N^2}\Big\}
\Eq(TL.l3.0)
$$
and
$$
\eqalign{
\wt\RR_1(u)&=O(|u|\wh\Th)
\cr
\wt\RR_2(u)&=O\left(\max\left\{
\VV^{\circ}_{N,d}(J\cup x)\left|\frac{-s(u)}{1-s(u)}\right|,\frac{1}{N^2}\left|\frac{-s(u)}{1-s(u)}\right|,\frac{|u|}{\barunder u(d)}
\right\}\right)\cr
}
\Eq(TL.l3.2)
$$
}

\remark Note that by \eqv(TL.t.0), for large enough $N$, $\left|Z_{N,d}(J\cup x)\right|< 3/4$.

\proofof{Lemma \thv(TL.lemma.3)} Using a second order Taylor expansion around $u=0$
to express $G^0_{0,J\cup x}(u)$,
and a first order expansion everywhere else, we obtain
$$
h(u)=
\frac {\left[G^y_{0,x\cup J}(0)+u\frac d{du}G^y_{0,x\cup J}(u')\right]
\left[G^x_{0,x\cup J}(0)+u\frac d{du}G^x_{0,x\cup J}(u'')\right]}
{\frac {\Q(0)}{\Q(x)}\left[1-G^0_{0,x\cup J}(0)-
u\frac d{du}G^0_{0,x\cup J}(0)-\frac {u^2}2\frac {d^2}{du^2}
G^0_{0,x\cup J}
(u''')\right]}
\Eq(TL.21)
$$
which can be rewritten as
$$
h(u)=
g\frac{f+uf'+u^2f''}{1-ug'\left(1+\frac{u}{2}g''\right)}
\Eq(TL.21')
$$
where
$$
f=G^y_{0,x\cup J}(0)\,,\quad
g=\frac {G^x_{0,x\cup J}(0)}
{\frac {\Q(0)}{\Q(x)}\left[1-G^0_{0,J\cup x}(0)\right]}
\Eq(TL.22)
$$
and, for some $0< u',u'',u'''<u$,
$$
\eqalign{
f'=&
\frac{\frac d{du}G^y_{0,x\cup J}(u')+\frac d{du}G^x_{0,x\cup J}(u'')}
{G^x_{0,x\cup J}(0)}
\,,\quad
f''=
\frac{\frac d{du}G^y_{0,x\cup J}(u')\frac d{du}G^x_{0,x\cup J}(u'')}
{G^x_{0,x\cup J}(0)}
\cr
g'=&
\frac{\frac d{du}G^0_{0,x\cup J}(0)}{1-G^0_{0,x\cup J}(0)}
\,,\quad
g''=
\frac{\frac {d^2}{du^2}G^0_{0,x\cup J}(u''')}{\frac d{du}G^0_{0,x\cup J}(0)}
\cr
}
\Eq(TL.23)
$$

Both $f$ and $g$ in \eqv(TL.22) are probabilities, namely
$$
g=\frac{\R\left(\t^0_x<\t^0_{0\cup J}\right)}
{\R\left(\t^0_{x\cup J}<\t^0_0\right)}
=\R\left(\t^0_x<\t^0_{J}\right)\quad\text{and}\quad
f=\R\left(\t^y_0<\t^y_{x\cup J}\right)\,,
\Eq(TL.31)
$$
that are well controlled through the results of Corollary \thv(P.cor.1)
and Lemma \thv(P.lemma.9) of Section 4. This yields \eqv(TL.l3.1).

The terms $f', f''$ and $g', g''$ will require some extra work.
While we will clearly need to get precise control on $g'$, rather rough bounds
on $f', f'', g''$ will suffice. To this aim the next lemma collects estimates
on the Laplace transforms appearing in \eqv(TL.20), together with estimates
and on their derivatives.

\lemma{\TH(TL.lemma.1)}{\it
Let $\phi(u)$ denote any of the Laplace transforms
$G^y_{x,J\cup 0}(u), G^y_{0,J\cup x}(u)$, $G^x_{0,J\cup x}(u)$, or
$G^0_{0,J\cup x}(u)$.
Let $\wh\Th(d)$ be given by \eqv(T.4bis). Then, for all $0<\e<1$ and all
real $u$ satisfying $u<(1-\e)/\wh\Th(d)$,
$$
\phi(u)\leq \frac 1{1-u\wh\Th(d)}\leq 1/\e
\Eq(TL.25)
$$
Therefore, $\phi(u)$ is analytic
for $u\in\C$ with $\Re(u)<(1-\e)/\wh\Th(d)$,
and,
if $|u|\leq (1-\e)/\wh\Th(d)$,
$$
\left|\phi(u)\right|\leq \frac {1}{\e}
\Eq(TL.26bis)
$$
$$
\eqalign{
\left|\frac d{du}\phi(u)\right|\leq &\frac {\wh\Th(d)}{\e(1-\e)}\cr
}
\Eq(TL.26)
$$
and
$$
\left|\frac {d^2}{du^2}\phi(u)\right|\leq \frac {2\wh\Th^2(d)}{\e(1-\e)^2}
\Eq(TL.27)
$$
}

\proofof{Lemma \thv(TL.lemma.1)} The proof of \eqv(TL.25) follows from the arguments used in
\cite{BEGK1}
(see
Section 3 of \cite{BEGK1};
see also Lemma 3.4 of \cite{BBG2})
for bounding Laplace transforms of positive random variables, together
with the bounds from Theorem \thv(T.theo.2).
The bound \eqv(TL.26bis) is then obvious since $\t^y_A$ is a positive
random variable, and \eqv(TL.26) and \eqv(TL.27)  result from the
Cauchy bound for derivatives of analytic functions.
\endproof

To control the term $g''$ we further need the following result.

\lemma{\TH(TL.lemma.2)}{\it
With the notation of Proposition \thv(TL.prop.1), we have:
$$
\E^{\circ}\t^{0}_{0}
\geq
\frac d{du}G^0_{0,I}(0)
\geq
\E^{\circ}\t^{0}_{0}(1-
C|I|2^{-N}N^{d+2})
\Eq(TL.30)
$$
for some constant $0<C<\infty$.
}

\proof From the identity
$$
\1_{\{\t^{0}_{0}<\t^{0}_{I}\}}=1 -\sum_{y\in I} \1_{\{\t^{0}_{y}<\t^{0}_{I\setminus y}\}}
\Eq(TL.l2.1)
$$
we  deduce that
$$
\eqalign{
G^0_{0,I}(u)
&=G^0_{0}(u)-\sum_{y\in I}G^y_{0}(u)G^0_{y,(I\setminus y)\cup 0}(u)
\cr
&=G^0_{0}(u)-\E^{\circ}\t^{0}_{0}
\sum_{y\in I}\Q(y)G^y_{0}(u)G^y_{0,I}(u)
}
\Eq(TL.l2.2)
$$
where the last line follows from reversibility together with the fact that
$\Q(0)\E^{\circ}\t^{0}_{0}=1$ (see the proof of Lemma \thv(L.4)).
Taking the derivative with respect to $u$, evaluated at $u=0$,
$$
\frac d{du}G^0_{0,x\cup J}(0)
=\E^{\circ}\t^{0}_{0}
-
\E^{\circ}\t^{0}_{0}
\sum_{y\in I}\Q(y)\left[
\E^{\circ}\t^{y}_{0}\P^{\circ}(\t^{y}_{0}<\t^{y}_{I})
+\E^{\circ}\t^{y}_{0}\1_{\{\t^{y}_{0}<\t^{y}_{I}\}}
\right]
\Eq(TL.l2.3)
$$
Now
$$
\E^{\circ}\t^{y}_{0}\P^{\circ}(\t^{y}_{0}<\t^{y}_{I})
+\E^{\circ}\t^{y}_{0}\1_{\{\t^{y}_{0}<\t^{y}_{I}\}}
\leq
2\E^{\circ}\t^{y}_{0}
\Eq(TL.l2.4)
$$
Hence
$$
\sum_{y\in I}\Q(y)\left[
\E^{\circ}\t^{y}_{0}\P^{\circ}(\t^{y}_{0}<\t^{y}_{I})
+\E^{\circ}\t^{y}_{0}\1_{\{\t^{y}_{0}<\t^{y}_{I}\}}
\right]
\leq
2\Q(I)\max_{y\in I}\E^{\circ}\t^{y}_{0}
\leq
|I|2^{-N+1}CN^{d+2}
\Eq(TL.l2.5)
$$
for some finite constant $C>0$, where the last inequality follows
from Theorem \thv(T.theo.2). Plugging this bound in \eqv(TL.l2.3) proves \eqv(TL.30).
\endproof

From now on
we assume that $u$ lies on the real
half line $u<(1-\e)/\wh\Th(d)$ for some fixed $0<\e<1$. Then, Lemma \thv(TL.lemma.1)
together with the probability estimates \eqv(P.c1.0) of Corollary \thv(P.cor.1)
immediately gives, assuming \eqv(TL.12bis) (which in fact implies that
$1-\VV^{\circ}_{N,d}(J\cup x)\geq 1/3$),
$$
\eqalign{
f'=&O(\wh\Th)\cr
f''=&O(\wh\Th^2)\cr
}
\Eq(TL.28)
$$
and by Lemma \thv(TL.lemma.1) and
Lemma \thv(TL.lemma.2), with $\barunder {\it u}(d)$ defined in \eqv(TL.12),
$$
0\leq g''\leq  3\bigl(\e(1-\e)^2\barunder {\it u}(d)\bigr)^{-1}
\Eq(TL.39)
$$
We now bound $g'$. Observe that by \eqv(P.c1.00) of
Corollary \thv(P.cor.1), with $Z_{N,d}(J\cup x)$ defined as in \eqv(TL.l3.0),
$$
1-G^0_{0,x\cup J}(0)
=\P^{\circ}\left(\t^{0}_{J\cup x}<\t^{0}_{0}\right)
=\frac{|J\cup x|}{2^N\Q(0)}\Bigl(1-\frac{1}{N}+Z_{N,d}(J\cup x)\Bigr)
\Eq(TL.40)
$$
Combining \eqv(TL.40) with Lemma \thv(TL.lemma.2) then yields,
$$
g'
=\frac {2^N\Q(0)}{|J\cup x|}\E^{\circ}\t^0_0
\Bigl(1+\frac{1}{N}+Z_{N,d}(J\cup x)\Bigr)
\Eq(TL.42)
$$
and since
$
\E\t^0_0{\Q^{}_N(0)}=1
$
(see the proof of Lemma \thv(L.4)),
$$
\eqalign{
g'
=&\frac {2^N}{|J\cup x|}
\Bigl(1+\frac{1}{N}+Z_{N,d}(J\cup x)\Bigr)
\cr
=&\bar u^{-1}\Bigl(1+Z_{N,d}(J\cup x)\Bigr)
\cr
}
\Eq(TL.43)
$$

We may now combine our estimates for $f',f'', g',g''$  with \eqv(TL.21').
Setting $s(u)\equiv u/\bar u$, $\bar h\equiv \bar u g'$,  and
$\barunder{\it h}\equiv \barunder{\it u}g''$ in \eqv(TL.21'), $h(u)$  can
be brought into the form \eqv(TL.48)
with
$$
\eqalign{
\wt\RR_1(u)&\equiv uf'+u^2f''
\cr
\wt\RR_2(u)&\equiv
\frac
{\displaystyle{
-\frac{u}{\bar u}(1-\bar h)
-\frac{1}{2}\frac{u}{\barunder{\it u}}\barunder{\it h}
-\frac{1}{2}\frac{u}{\barunder{\it u}}\frac{u}{\bar u}\barunder{\it h}\bar h
}}
{\displaystyle{
1-\frac{u}{\bar u}
\left(1+\frac{1}{2}\frac{u}{\barunder{\it u}}\barunder{\it h}\right)
}}
}
\Eq(TL.29)
$$
Now, by \eqv(TL.28),
$$
\wt\RR_1(u)=O(|u|\wh\Th)
\Eq(TL.29')
$$
and knowing from \eqv(TL.39) and \eqv(TL.43) that
$$
\eqalign{
&\bar h=1+Z_{N,d}(J\cup x)
\cr
&0<\barunder{\it h}<const.\cr
}
\Eq(TL.46)
$$
one easily checks that
$$
\wt\RR_2(u)
=
\RR_2(u)=O\left(\max\left\{
\VV^{\circ}_{N,d}(J\cup x)\left|\frac{-s(u)}{1-s(u)}\right|,\frac{1}{N^2}\left|\frac{-s(u)}{1-s(u)}\right|,\frac{|u|}{\barunder u(d)}
\right\}\right)
\Eq(TL.47)
$$
for all $u$ satisfying $-\rho\barunder{\it u}(d)<u< \bar u$ for some $0<\rho<1$.
This concludes the proof of Lemma \thv(TL.lemma.3).\endproof

Let us now turn to the first term in the r\.h\.s\. of \eqv(TL.20). Here,
we will simply write
$$
G^y_{x,J\cup 0}(u)=G^y_{x,J\cup 0}(0)+u\frac d{du}G^y_{x,J\cup 0}(\tilde u)
=\R\left(\t^y_x<\t^y_{J\cup 0}\right)+\wt\RR_3(u)
\Eq(TL.49)
$$
where, for some $0\leq \tilde u\leq u$,
$$
\wt\RR_3(u)=u\frac d{du}G^y_{x,J\cup 0}(\tilde u)=O(|u|\wh\Th)
\Eq(TL.50)
$$
the last equality above being Lemma \thv(TL.lemma.1) again.

We may now collect our estimates.
Adding \eqv(TL.48) and \eqv(TL.49) yields, 
$$
G^y_{x,J}(u)=\R\left(\t^y_x<\t^y_{J\cup 0}\right)+\frac{fg}{1-s(u)}+\wt\RR_0(u)
\Eq(TL.51)
$$
where
$$
\wt\RR_0(u)=
\frac{g}{1-s(u)}\left[\wt\RR_1(u)\left(1+\wt\RR_2(u)\right)+f\wt\RR_2(u)\right]
+\wt\RR_3(u)
\Eq(TL.52)
$$
To arrive at \eqv(TL.15) first observe that by  \eqv(TL.31), \eqv(TL.52) becomes
$$
\wt\RR_0(u)=
\frac{\R\left(\t^0_x<\t^0_{J}\right)}{1-s(u)}
\left[\wt\RR_1(u)\left(1+\wt\RR_2(u)\right)
+\R\left(\t^y_0<\t^y_{x\cup J}\right)\wt\RR_2(u)\right]
+\wt\RR_3(u)
\Eq(TL.53)
$$
and, in view of \eqv(TL.29), \eqv(TL.47), and \eqv(TL.50), we may choose
$
\RR_1(u)=\wt\RR_1(u)\left(1+\wt\RR_2(u)\right)$, $\RR_2(u)=\wt\RR_2(u)$,
 $\RR_3(u)=\wt\RR_3(u),
$
and set
$$
\RR_0(u)=\wt\RR_0(u)
\Eq(TL.55)
$$
On the other hand
$$
\eqalign{
fg&=\R\left(\t^y_0<\t^y_{x\cup J}\right)\R\left(\t^0_x<\t^0_{J}\right)\cr
&=\R\left(\t^y_0<\t^y_x<\t^y_{J}\right)\cr
&=\R\left(\t^y_x<\t^y_{J}\right)-\R\left(\t^y_x<\t^y_{J\cup 0}\right)\cr
}
\Eq(TL.56)
$$
which, inserted in \eqv(TL.51) yields \eqv(TL.14).
This concludes the proof of assertion (i) for $J\neq\emptyset$.
The proof of the case $J=\emptyset$ is a straightforward rerun of
the case $J\neq\emptyset$, taking \eqv(TL.11) rather than \eqv(TL.11)
for starting point. The first assertion of Proposition \thv(TL.prop.1) is proven.

\remark Note that
\eqv(TL.20) implies that $G^y_{x,J}(u)$ has a pole at the point
$u^*>0$ defined as the smallest real number that solves the equation
$G^0_{0,J\cup x}(u)=1$.  Now our estimates imply that  $u^*\approx \bar u$
and, from its boundedness at $u=\bar u$,
that $G^y_{x,J}(u)$ is analytic  for all
 for $u\in\C$ satisfying $\Re(u)<\bar u$.
One then checks that assertion (i) remains
valid 
in the region
of the complex plane given by
$|u|\leq (1-\e)/\wh\Th(d)$ intersected with
$-\rho\barunder {\it u}(d)<\Re(u)< \bar u$.


We now turn to the proof of assertion (ii).
Again we start with \eqv(TL.20) and call the second summand $h(u)$.
Clearly, if $u\leq 0$,
$$
G^y_{0,x\cup J}(u) \equiv
\E^{\circ} e^{u\t^y_0}\1_{\{\t^y_0<\t^y_{x\cup J}\}}
\leq \E^{\circ} e^{u\ell_y}\1_{\{\t^y_0<\t^y_{x\cup J}\}}
= e^{u\ell_y}G^y_{0,x\cup J}(0)
\Eq(TL.60)
$$
and in the same way
$$
G^0_{x,0\cup J}(u) \leq e^{u\ell_y}G^0_{x,0\cup J}(u)
\Eq(TL.61)
$$
Moreover, all Laplace transforms and their derivatives are positive monotone
increasing functions of $u$;
thus, for $u\leq -\rho\barunder {\it u}(d)$,
$$
1-G^0_{0,x\cup J}(u)\geq 1-G^0_{0,x\cup J}(-\rho\barunder {\it u}(d))
\Eq(TL.62)
$$
and, using \eqv(TL.60), \eqv(TL.61) and \eqv(TL.62) in $h(u)$,
$$
h(u)\leq
\frac {G^y_{0,x\cup J}(0)G^x_{0,x\cup J}(0)}
{\frac {\Q(0)}{\Q(x)}\left[1-G^0_{0,J\cup x}(-\rho\barunder {\it u}(d))\right]}
=\frac{fg}{1+\rho\barunder {\it u}(d)g'
\left(1-\frac{1}{2}\rho\barunder {\it u}(d)g''\right)}
\Eq(TL.63)
$$
where $f,g, g'$ and $g''$ are as in \eqv(TL.22) and  \eqv(TL.23)
for some $0\leq u'''\leq \rho\barunder {\it u}(d)$. By \eqv(TL.l3.1), \eqv(TL.39),
 \eqv(TL.43), and \eqv(TL.12bis),
$$
h(u)\leq
\frac{e^{-|u|(\ell_x+\ell_y)}}{1+s(\barunder{\it u}(d))
(\rho/4)(1-9\rho)}
\Bigl(1-\frac{1}{N}+Z_{N,d}(J\cup x)\Bigr)
\Eq(TL.64)
$$
and inserting \eqv(TL.64) in \eqv(TL.20) yields  \eqv(TL.19).
The proof of Proposition \thv(TL.prop.1) is done.\endproof

We are now ready to prove Theorem \thv(TL.theo.1).

\proofof{Theorem \thv(TL.theo.1)} By \eqv(TL.4), setting $J=I\setminus x$
and $u={s}/{\E^{\circ}\t^{y}_{I}}$,
$$
\E^{\circ}\left(
e^{s{\t^{y}_{I}}/\E^{\circ}\t^{y}_{I}}\1_{\{\t^{y}_{x}< \t^{y}_{I\setminus x}\}}\right)
=G^y_{x,J}(u)
\Eq(TL.t.3)
$$
Under the assumptions of Theorem \thv(TL.theo.1) we may use assertion (i) of
Proposition \thv(TL.prop.1) to express the Laplace transform \eqv(TL.t.3). We will only
treat the case $J\neq\emptyset$, namely use \eqv(TL.14).
(The case $J=\emptyset$ is similar but simpler since it relies on the use of \eqv(TL.15).)
We first have to verify that \eqv(TL.14) is valid on the domain $-\infty<s<1-\e$,
for all $\e>0$. Recall that \eqv(TL.14) was established for
$-\rho\barunder u(d)<u< \bar u$ for some $0<\rho<1$ thus,
making the change of variable  $u={s}/{\E^{\circ}\t^{y}_{I}}$, for
$-\rho\barunder u(d)\E^{\circ}\t^{y}_{I}<s< \bar u\E^{\circ}\t^{y}_{I}$.
Now, as in \eqv(TL.71), we may write
$$
\E^{\circ}\t^{y}_{I}=\frac{2^{N}}{|I|}
\Big(1+\frac{1}{N}\Big)
(1+O(\wt\varepsilon^{\circ}_{N,d}(I,y)))
\Eq(TL.t.4)
$$
Note here that $|I|\leq |\SS_d|=2^{d}$,
and since by assumption $d\leq d_0(N)=o(N)$,
$$
\frac{2^{N}}{|I|}\geq 2^{N(1-o(1))}
\Eq(TL.t.15)
$$
Moreover, by \eqv(TL.t.0), $\wt\varepsilon^{\circ}_{N,d}(I,y)\leq\d/2$.
Thus, together with \eqv(TL.12), \eqv(TL.t.4) yields
$$
\bar u\E^{\circ}\t^{y}_{I}
=1+O(\wt\varepsilon^{\circ}_{N,d}(I,y))
\geq 1-\d/2\geq 1-\e/2\,,\quad
\Eq(TL.t.16)
$$
and
$$
\barunder u(d)\E^{\circ}\t^{y}_{I}
=\barunder u(d)\frac{2^{N}}{|I|}
\Big(1+\frac{1}{N}\Big)(1+O(\wt\varepsilon^{\circ}_{N,d}(I,y)))
\geq 2^{N(1-o(1))}
\Eq(TL.t.17)
$$
where we used that $\barunder u(d)$ is polynomial in $N$
together with \eqv(TL.t.15). Clearly, for all $\e>0$, for all
$-\infty<s<1-\e$, choosing $N$ large enough guarantees that
$
-\rho2^{N(1-o(1))}<s< 1-\e/2
$

Let us next consider the terms $\frac{1}{1-s(u)}$ and $\frac{-s(u)}{1-s(u)}$ in \eqv(TL.14).
Using again \eqv(TL.t.4) we have, by  \eqv(TL.12) and \eqv(TL.13), for $u={s}/{\E^{\circ}\t^{y}_{I}}$,
$$
s(u)=
\frac{u}{\bar u}=\frac{s}{\bar u\E^{\circ}\t^{y}_{I}}
=s(1+O(\wt\varepsilon^{\circ}_{N,d}(I,y)))
\Eq(TL.t.5)
$$
and
$$
\frac{1}{1-s(u)}
=
\frac{1}{1-s}\big(1+O(\wt\varepsilon^{\circ}_{N,d}(I,y))\big)
\Eq(TL.t.6)
$$
Let us now consider the two probabilities $\R^\circ\left(\t^y_x<\t^y_{J}\right)$ and
$\P^{\circ}\left(\t^{y}_{x}<\t^{y}_{J\cup 0}\right)$:
on the one hand Theorem \thv(P.theo.00) gives
$$
\R^\circ\left(\t^y_x<\t^y_{J}\right)
=
\frac{1}{|I|}\Big(1+Z_{N,d}(J\cup x)+|I|\phi_x(\dist(x,y))\Big)
\Eq(TL.t.7)
$$
while on the other hand
$$
0\leq \P^{\circ}\left(\t^{y}_{x}<\t^{y}_{J\cup 0}\right)\leq
\P^{\circ}\left(\t^{y}_{x}<\t^{y}_{0}\right)
\leq\frac{1}{|I|}\Big(|I|\phi_x(\dist(x,y))\Big)
\Eq(TL.t.8')
$$
where the rightmost inequality follows from \eqv(N.2). At this stage we see,
inserting the estimates \eqv(TL.t.5)-\eqv(TL.t.8') in  \eqv(TL.14), that
for all $\e>0$ and all $-\infty<s<1-\e$, for $N$ large enough,
and for $\varepsilon^{\circ}_{N,d}(I,x,y)$ defined in
\eqv(TL.t.2),
$$
\left|G^y_{x,J}(u)-\frac{1}{|I|}\frac{1}{1-s}\right|
\leq\frac{ c_\e}{|I|}
\varepsilon^{\circ}_{N,d}(I,x,y)
+\RR_{0}({s}/{\E^{\circ}\t^{y}_{I}})
\Eq(TL.t.8)
$$
for some constant $0<c_\e<\infty$ that does not depend on $N, I$, or $d$,
but on $\e$ only. It remains to bound $\RR_{0}(u)$ for  $u={s}/{\E^{\circ}\t^{y}_{I}}$.
To deal with $\RR_1(u)$ and $\RR_3(u)$ (see \eqv(TL.16),  \eqv(TL.17))
note that by \eqv(TL.t.4),
$$
|u|\wh\Th(d)
=\frac{|s|}{|I|}\left(\frac{|I|^2}{2^{N}}\wh\Th(d)\right)
\Big(1-\frac{1}{N}\Big)(1+O(\wt\varepsilon^{\circ}_{N,d}(I,y)))
\Eq(TL.t.9)
$$
Reasoning as in \eqv(TL.t.15)
we get, for $\wh\Th(d)$ defined in \eqv(T.4),
$$
\frac{|I|^2}{2^{N}}\wh\Th(d)\leq 2^{-N(1-o(1))}
\Eq(TL.t.10)
$$
To bound the term ${|u|}/{\barunder u(d)}$ in $\RR_2(u)$ observe that,
reasoning again as in \eqv(TL.t.15),
$$
\frac{|u|}{\barunder u(d)}\leq |u|\wh\Th(d)^2\leq\frac{|s|}{|I|}2^{-N(1-o(1))}
\Eq(TL.t.11)
$$
where the leftmost inequality follows from \eqv(TL.12). Hence, for $1\leq i\leq 3$,
for some constant $0<c<\infty$,
$$
\RR_{i}(u)\leq c\frac{|s|}{|I|}2^{-N(1-o(1))}\wt\varepsilon^{\circ}_{N,d}(I,y)
\Eq(TL.t.12)
$$
Combining  \eqv(TL.16) with \eqv(TL.t.12)  and the estimates on
$\R\left(\t^y_0<\t^y_{x\cup J}\right)$
and $\R\left(\t^y_0<\t^y_{x\cup J}\right)$ from \eqv(TL.l3.1), we obtain
$$
\RR_{0}(u)\leq
c'\frac{|s|}{|I|}2^{-N(1-o(1))}\wt\varepsilon^{\circ}_{N,d}(I,y)
\Eq(TL.t.13)
$$
for some constant $0<c'<\infty$.
Putting \eqv(TL.t.8) and \eqv(TL.t.13) together yields
$$
\left|G^y_{x,J}(u)-\frac{1}{|I|}\frac{1}{1-s}\right|
\leq c_\e\frac{1}{|I|}\varepsilon^{\circ}_{N,d}(I,x,y)
\Eq(TL.t.14)
$$
and proves Theorem \thv(TL.theo.1). \endproof

\proofof{Theorem \thv(TL.theo.2)} Theorem \thv(TL.theo.2) follows from
Corollary \thv(TL.cor.3) in the same way that Theorem \thv(TL.theo.1) follows
from Proposition \thv(TL.prop.1). We skip the details.\endproof.

\bigskip
\line{\bf 7.3.  Back to the hypercube.\hfill}

To go from the lumped chain back to the chain on the hypercube we proceed as in
Chapter 5, but rely on Lemma \thv(L.7) rather than Lemma \thv(L.6).

\proofof{Theorem \thv(TLC.theo.1)} We will show that
Theorem \eqv(TLC.theo.1) is a consequence of
Theorem \eqv(TL.theo.1). Under the assumptions
and with the notation of Theorem \eqv(TLC.theo.1) set
$I=\g(A)$, $x=\g(\eta)$, and $y=\g(\s)$.
Here $x\in I$, $y\notin I$, and $I\cup y\subset\SS_d$ i\.e\.
$A\cup\s$ is $\g$-compatible (see Definition \thv(I.def.1bis)).
Now by Lemma \thv(B.4), $\VV_{N,d}(A\cup\s)=\VV^{\circ}_{N,d}(I\cup y)$,
implying that the conditions \eqv(TLC.t.0) and  \eqv(TL.t.0) are equivalent.
Next, by Lemma \thv(L.7)
$$
\E\left(
e^{s{\t^{\s}_{A}}/\E\t^{\s}_{A}}
\1_{\{\t^{\s}_{\eta}< \t^{\s}_{A\setminus\eta}\}}\right)
=
\E^{\circ}\left(
e^{s{\t^{y}_{I}}/
\E^{\circ}\t^{y}_{I}
}
\1_{\{\t^{y}_{x}< \t^{y}_{I\setminus x}\}}\right)
\Eq(TLC.10)
$$
It thus remains to see that
$
\varepsilon_{N,d}(A,\eta,\s)=\varepsilon^{\circ}_{N,d}(I,x,y)
$.
But this holds true since by Lemma \thv(N.lemma.3),
$
|A|\phi_{\g(\eta)}(\Dist(\s,\eta))=|I|\phi_x(\dist(x,y))
$,
and by Lemma \thv(B.lemma.5) and Lemma \thv(B.4),
$
\wt\varepsilon_{N,d}(A,\s)=\wt\varepsilon^{\circ}_{N,d}(I,y)
$.
Hence \eqv(TLC.t.1) and \eqv(TL.t.1) also are equivalent.
Theorem \eqv(TL.theo.1) thus implies Theorem \eqv(TLC.theo.1).

That \eqv(TLC.t.0) and \eqv(TLC.t.1) remain true with
$\phi_{\g(\eta)}(\Dist(\s,\eta))$ replaced by
$F(\Dist(\s,\eta))$ in \eqv(TLC.t.2), and
with $\VV_{N,d}(A\cup\s)$ replaced by $\UU_{N,d}(A\cup\s)$
in \eqv(TLC.t.0) and \eqv(TLC.2')
simply follows from Lemma \thv(TL.lemma.4) and Lemma \thv(N.lemma.2).
\endproof

\proofof{Theorem \thv(TLC.theo.2)} Theorem \thv(TLC.theo.2)
is deduced from the
special case of Theorem \thv(TL.theo.2)
obtained by taking $d=1$.
To see this choose
$\g(\s')=\frac{1}{N}\sum_{i=1}^{N}\s'_i\eta_i$,
$\s'\in\SS_N$. Setting $y=\g(\s)$ and $x=\g(\eta)$, we
of course have $y\in\G_{N,1}$ and $x\in\SS_N$. Then, by lemma \thv(L.7)
with $A=\{\eta\}$, using that $\E\t^{\s}_{\eta}=\E^{\circ}\t^{x}_{y}$
(see \eqv(T.7ter)), we get that
$
\E e^{s{\t^{\s}_{\eta}/2^N}}
=
\E^{\circ}e^{s{\t^{y}_{x}}/2^N}
$.
Finally by Lemma \thv(N.lemma.3), $\dist(x,y)=\Dist(\s,\eta)$.
The proof is done.
\endproof

\proofof{Corollary \thv(TLC.cor.1)}
Let the assumptions and the notation be those of
Theorem \thv(TLC.theo.1) and its proof. Note that
$$
\E e^{s{\t^{\s}_{A}}/\E\t^{\s}_{A}}
=
\sum_{\eta\in A}\E e^{s{\t^{\s}_{A}}/\E\t^{\s}_{A}}
\1_{\{\t^{\s}_{\eta}< \t^{\s}_{A\setminus\eta}\}}
=
\sum_{x\in I}\E^{\circ}
e^{s{\t^{y}_{I}}/\E^{\circ}\t^{y}_{I}}
\1_{\{\t^{y}_{x}< \t^{y}_{I\setminus x}\}}
\Eq(TLC.11)
$$
The first equality in \eqv(TLC.11) is the analogue of \eqv(TL.5) for the
chain on the hypercube, and the last follows from Lemma \thv(L.7).
Corollary \thv(TLC.cor.1)
is then deduced from Corollary \thv(TL.cor.1)
in the same way that Theorem \thv(TLC.theo.1)
was deduced from Theorem \eqv(TL.theo.1).
\endproof

As announced in the Section 1 we now specialize Theorem \thv(I.theo.4nonsparse)
to the case where the starting point $\s$ is chosen in $\WW(A,|A|)$ and the
condition \eqv(I''.24) is replaced by
$
\UU_{N,d}(A)=o(1)
$.

%
%

\theo{\TH(I.theo.4)}{\it  Let $d'\leq d_0(N)/2$  and let $\L'$ be a log-regular
$d'$-partition.  Assume that  $A\subset\SS_N$ is compatible with $\L'$.
Then for all $\s\in\WW(A,|A|)$ there exists
an integer $d$ with $d'<d\leq 2d'$ such that if
$$
\eqalign{
&\UU_{N,d}(A)=o(1)\,, \quad N\rightarrow\infty\cr
}
\Eq(I.24')
$$
the following holds for all $\eta\in A$:
for all $\e>0$, there exists a constant $0<c_\e<\infty$
(independent of $\s, |A|, N$, and $d$) such that, for all $s$ real satisfying
$-\infty<s<1-\e$, we have, for all $N$ large enough,
$$
\left|
\E\left(
e^{s{\t^{\s}_{A}}/\E\t^{\s}_{A}}
\1_{\{\t^{\s}_{\eta}<\t^{\s}_{A\setminus\eta}\}}\right)
-\frac{1}{|A|}\frac{1}{1-s}
\right|\leq
\frac{1}{|A|}c_{\e}\max\Big\{\UU_{N,d}(A),\frac{1}{N^k},
|A|F_{N,d}(\rho(|A|)+1)\Big\}
\Eq(I.24)
$$
where
$$
k=
\cases 2, &\text{if} H(A\cup\s) \text{is satisfied}\cr
1, &\text{if} H(A\cup\s) \text{is not satisfied.}
\endcases
\Eq(I.25)
$$
}

We finally prove Theorem \thv(I.theo.4), and Theorem \thv(I.theo.4nonsparse)
and Corollary \thv(I.cor.theo.4) of Section 1.

\proofof{Theorem \thv(I.theo.4)}
We want to show that when restricting the starting point $\s$ to
sets of the form $\WW(A,|A|)$ (see \eqv(I.14)) and when replacing the assumption
\eqv(TLC.t.0) by the stronger assumption
$$
\VV_{N,d}(A\cup\s)=o(1)\,, \quad N\rightarrow\infty
\Eq(TLC.t.0bis)
$$
Theorem \thv(TLC.theo.1)
entails Theorem \thv(I.theo.4) for all large enough $N$.
Thus let the assumptions and the notation be those of
Theorem \thv(TLC.theo.1) but take $\s\in\WW(A,|A|)$, assume \eqv(TLC.t.0bis)
instead of \eqv(TLC.t.0)
(i\.e\. assume that $\d\equiv\d(N)\rightarrow 0$ as $N\rightarrow\infty$), and
let $\g$ be any $d$-lumping compatible with $A\cup\s$
(recall that this in particular implies that $\g(A\cup\s)\in\SS_d$).

Proceeding as we did in the proof of Theorem \thv(I.theo.1) to obtain
\eqv(B.10quart)
we get, for all $\s\in\WW(A,|A|)$ and all $\eta\in A$,
$$
\phi_{\g(\eta)}(\Dist(\s,\eta))
\leq
\phi_{\g(\eta)}(\rho(|A|)+1)\leq F_{N,d}(\rho(|A|)+1)
\Eq(TLC.12)
$$
$$
\phi_{\g(\s)}(\Dist(\eta,\s))
\leq
\phi_{\g(\s)}(\rho(|A|)+1)\leq F_{N,d}(\rho(|A|)+1)
\Eq(TLC.12')
$$
Eq\. \eqv(TLC.12) implies in particular that
$$
\sum_{\eta\in A}\phi_{\g(\eta)}(\Dist(\s,\eta))
\leq
|A|F_{N,d}(\rho(|A|)+1)
\Eq(TLC.13)
$$
Now by \eqv(I.12new) with $|A|>1$,
$$
\eqalign{
\VV_{N,d}(A\cup\s)
&=
\max_{\eta\in A\cup\s}\sum_{\eta'\in (A\cup\s)\setminus\eta}
\phi_{\g(\eta)}(\Dist(\eta,\eta'))
\cr
&=
\max\left\{
\sum_{\eta'\in A}\phi_{\g(\eta')}(\Dist(\s,\eta'))
+
\max_{\eta\in A}
\left[
\sum_{\eta'\in A\setminus\eta}\phi_{\g(\eta')}(\Dist(\eta,\eta'))
+\phi_{\g(\s)}(\Dist(\eta,\s))
\right]
\right\}
}
\Eq(TLC.14)
$$
Together with  \eqv(TLC.12') and \eqv(TLC.13), \eqv(TLC.14) yields
$$
\VV_{N,d}(A)
\leq
\VV_{N,d}(A\cup\s)
\leq
(|A|+1)\max_{\eta\in A} F_{N,d}(\rho(|A|)+1)+\VV_{N,d}(A)
\Eq(TLC.15)
$$
In view of \eqv(I.14)-\eqv(I.15), \eqv(TLC.15) implies that
$$
\VV_{N,d}(A\cup\s)=o(1)\text{iff} \VV_{N,d}(A)=o(1)
\Eq(TLC.16)
$$
Thus, for all $\s\in\WW(A,|A|)$, assumption \eqv(TLC.t.0bis)
is equivalent to $\VV_{N,d}(A)=o(1)$.
Let us now consider the term $\varepsilon_{N,d}(A,\eta,\s)$
in \eqv(TLC.t.1). By \eqv(TLC.12) and \eqv(TLC.15),
$
\varepsilon_{N,d}(A,\eta,\s)\leq\wh\varepsilon_{N,d}(A,\eta,\s)
$
where
$$
\wh\varepsilon_{N,d}(A,\eta,\s)
=
\frac{1}{|A|}O\Big(\max\Big\{\VV_{N,d}(A),\frac{1}{N^k},
|A|F_{N,d}(\rho(|A|)+1)\Big\}\Big)
\Eq(TLC.17)
$$
where
$$
k=
\cases 2, &\text{if} H(A\cup\s) \text{is satisfied}\cr
1, &\text{if} H(A\cup\s) \text{is not satisfied.}
\endcases
\Eq(TLC.18)
$$

Gathering the previous results we conclude that,
under the assumptions of Theorem \thv(TLC.theo.1),
restricting $\s$ to the set $\WW(A,|A|)$,
we have, for all $A\subset \SS_N$ such that
$\VV_{N,d}(A)=o(1)$ and all $\eta\in A$,
that \eqv(TLC.t.1) holds true  with
$\varepsilon_{N,d}(A,\eta,\s)$ replaced by
$\wh\varepsilon_{N,d}(A,\eta,\s)$. This would be the
statement of Theorem \thv(I.theo.4) if we could
replace $\VV_{N,d}(A\cup\s)$ by $\UU_{N,d}(A\cup\s)$.
But this is made possible by
Lemma \thv(TL.lemma.4). Once this replacement done,
all dependence on the
choice of the underlying $d$-lumping $\g$ is suppressed.
Theorem \thv(I.theo.4) is thus proven.\endproof

\proofof{Theorem \thv(I.theo.4nonsparse)} Theorem \thv(I.theo.4nonsparse)
follows from Theorem \thv(TLC.theo.1) in the same way that
Theorem \thv(I.theo.1) follows from  Theorem \thv(B.theo.00).
We skip the details.\endproof

\proofof{Corollary \thv(I.cor.theo.4)} Corollary \thv(I.cor.theo.4)
is deduced from Theorem \thv(I.theo.4) just as
Corollary \thv(I.cor.theo.1) is deduced from Theorem \thv(I.theo.1).
Again we skip the details.\endproof



\medskip
\thanks
The authors thank A\. Bovier who participated in an early version of this work.
V.G. thanks the Centre Interdisciplinaire Bernoulli of the EPFL, Lausanne,
for hospitality, and the Weierstrass Institute, Berlin, for hospitality and
financial support.

\bigskip

\bigskip
\vfill\eject


\bigskip
\chap{8. Appendix A1}8

We state here the two simple lemmata that are used in Chapter 3 to
bound `no return before hitting' probabilities, i\.e\. probabilities
of the form
$$
\P^{\circ}(\t^x_J<\t^x_x) \text{for} J\subset\G_{N,d} \text{and} x\in \G_{N,d}\setminus x.
$$
The first of these lemmata is a partial analogue, for our reversible Markov chains, of the
classical Dirichlet principle from potential theory.
Let $\HH^x_J$ be the space of functions
$$
\HH^x_J\equiv\left\{
h\: \G_{N,d}\rightarrow [0,1] \mid h(x)=0,
\,\,\hbox{\rm and} \,\,h(y)=1 \,\,\hbox{\rm for}\,\, y\in J
\right\}
\Eq(A.8)
$$
and define the Dirichlet form
$$
\Phi_{N,d}(h)\equiv
\frac 12\sum_{y',y''\in\G_{N, d}} \Q_N(y')r_N(y',y'')
[h(y')-h(y'')]^2
\Eq(A.6)
$$
Note that the function $h^x_J$ defined in \eqv(A.9) below is in
$\HH^x_J$:
$$
h^x_J(y)=\cases 1,&\text{if}\,y\in J\cr
            0,&\text{if}\,y=x\cr
            \P^{\circ}(\t^y_J<\t^y_x),  &\text{if}\,y\notin J\cup x
\endcases
\Eq(A.9)
$$
The following lemma can be found e\.g\. in Liggett's book ([Li], pp 99, Theorem 6.1).

\lemma{\TH(A.lemma.1)}{\it Let $J\subset\G_{N,d}$ and $x\in \G_{N,d}\setminus x$.
Then
$$
\Q_N(x)\P^{\circ}(\t^x_J<\t^x_x) =
\inf_{h\in \HH^x_J}\Phi_{N,d}(h)=\Phi_{N,d}(h^x_J)
\Eq(A.10)
$$
}

\remark Note that `no return before hitting' probabilities are closely related to the notion of
capacity since, in potential theoretic language, the capacitor $(x,J)$
has capacity
$$
\capa_x(J)=\left[\Q_N(x)\P^{\circ}(\t^x_J<\t^x_x)\right]^{-1}
\Eq(A.11)
$$

Clearly, guessing the minimizing function $h^x_J$  in \eqv(A.10) yields an
upper bound on $\P^{\circ}(\t^x_J<\t^x_x)$.
To get a lower bound we use that the $d$-dimensional variational problem
\eqv(A.10) can be compared to a sum of (hopefully easier to handle)
one-dimensional ones.
This idea was heavily exploited in [BEGK1,2] and [BBG1];
the next lemma is quoted from [BBG1] (Lemma 4.1 of the appendix).


\lemma{\TH(A.lemma.2)}{\it Let $\D_{k}\subset\G_{N,d}$, $1\leq k\leq K$,
be a collection subgraphs of $\G_{N,d}$ and let $\wt\P^{\circ}_{\D_{k}}$
denote the law of the Markov chain with transition rates
$$
\wt r_{\D_{k}}(x',x'')
=\cases r_N(x',x'') ,&\text{if} x'\neq x''
                              \,,\hbox{and}\,\, (x',x'')\in E({\D_{k}})\cr
                           0,&\text{otherwise}\endcases
\Eq(A.1)
$$
and invariant measure
$$
\wt\Q^{\circ}_{\D_{k}}(y)=\Q_N(y)/\Q_N(\D_{k})\,, \quad y\in\D_{k}\,.
\Eq(A.1bis)
$$
Assume that
$$
E(\D_{k})\cap E(\D_{k'})=\emptyset,\,\,\,\,\forall k,k'\in\{1,\dots,K\},
k\neq k'
\Eq(A.2)
$$
and that
$$
y,x\in\bigcap_{k=1}^K V(\D_{k})
\Eq(A.3)
$$
Then
$$
\P^{\circ}\left(\t^y_x<\t_y^y\right)
\geq\sum_{k=1}^K\wt\P^{\circ}_{\D_{k}}\left(\t^y_x<\t_y^y\right)
\Eq(A.4)
$$
}

\bigskip
\chap{9. Appendix A2}9

As we just saw it is crucial in our approach to get sharp lower bounds on one
dimensional `no return before hitting probabilities'. The next lemma provide such bounds for
$\P^{\circ}(\t^x_0<\t^x_x)$.

\lemma{\TH(A.lemma.3)}{\it Let $d=1$ and let
$x\equiv x(N)\in
\G_{N,1}\equiv
\left\{1-\frac{2k}{N}, 0\leq k\leq N\right\}
$.
Then, setting
$$
\varrho_{N,1}(x):=\P^{\circ}(\t^x_0<\t^x_x)\,
\Eq(A.19)
$$
the following holds: $\varrho_{N,1}(x)=\varrho_{N,1}(-x)$ and

\item{i)} if $x(N)=1$,
$$
\varrho_{N,1}(x)
=1-\sfrac{1}{N}+O(\sfrac{1}{N^2})
\Eq(A.20)
$$

\item{ii)} if $\lim_{N\rightarrow\infty} x(N)=x_{\infty}>0$
then there exists constants $c_0, c_1>0$ and $c_3>1$ such that
$$
\varrho^{-1}_{N,1}(x)
\leq
\frac{1}{x}\left[1+\frac{c_0}{N}\frac{\log N}{\left|\log\left(\frac{1-x}{2}\right)\right|}\right]
+\frac{c_1}{N^{c_3}}
\Eq(A.21)
$$

\item{iii)} if $\lim_{N\rightarrow\infty} x(N)=0$ and
$\lim_{N\rightarrow\infty}  x(N)\sqrt{N}=\infty$,
$$
\varrho^{-1}_{N,1}(x)
\leq\frac{1}{x}(1+o(1))
\Eq(A.22)
$$

\item{iv)} if $\lim_{N\rightarrow\infty} x(N)=0$ and
$x(N)\sqrt{N}=O(1)$,
$$
\varrho^{-1}_{N,1}(x)
\leq Nx(1+o(1))
\Eq(A.23)
$$
Gathering the previous bounds,
$$
\inf_{x\in\G_{N,1}}\varrho^{-1}_{N,1}(x)
\leq C\sqrt N
\Eq(A.24)
$$
for some constant $0<C<\infty$.
}

\proof Without loss of generality we may assume that $N$ is even.
Given $x\in\G_{N,1}$ set $m=\frac{N}{2}(1-x)$, $L=\frac{N}{2}-m$, and
$$
\o_0(n)=1-\frac{2}{N}(m+n)\,,\quad 0\leq n\leq L
$$
Then, formula \eqv(P.l5.13) (or equivalently  \eqv(P.l5.14)) shows that,
$$
\P^{\circ}(\t^x_0<\t^x_x)
=
\left[
\sum_{n=0}^{L-1}
\frac{\Q_N(\o_0(m))}{\Q_N(\o_n(m))}
\frac{1}{\frac{1}{2}(1+\o_n(m))}
\right]^{-1}
=
\left[
\sum_{n=0}^{L-1}
\frac{\binom{N}{m}}{\binom{N}{n+m}}
\frac{N}{N-(n+m)}
\right]^{-1}
\Eq(A.25)
$$
From this explicit formula, it is a routine (though tedious task)
to prove Lemma \thv(A.lemma.3).
(Note that from this same formula we can of course also derive
lower bounds on $\varrho^{-1}_{N,1}(x)$.)
\endproof

\bigskip
\vfill\eject

\bigskip
\chap{10. Appendix A3}{10}


We now focus on the function $F(n)=F_1(n)+F_2(n)$ of Definition \thv(P.def.1).
This function is used to control the smallness of (hopefully) sub-leading terms in
virtually all our estimates (see e\.g\. \eqv(I.17), \eqv(I.22bis) or \eqv(I''.24));
it is in particular used through \eqv(I.12) to define the sparseness of sets $A\in\SS_N$
(see Definition \thv(I.def.3)). For practical purposes however the complexity of
the function $F_2(n)$ is a serious hindrance.
Our main aim in this appendix is to provide simpler, workable, expressions for $F_2(n)$,
for all $d\leq N$ and $N$ large enough.

Our main result is Lemma \thv(A3.lemma.1). It contains a collection of upper bounds
on $F_2(n)$ that suggest very strongly that for `small $n$', namely for $n+2\leq d$,
$F_2(n)$ has two distinct asymptotic behaviors, depending on whether the ratio $\frac{d^2}{N}$
goes to zero or not as $N$ diverges. Indeed our upper bound on $F_2(n)$ is
essentially independent of $d$ when $\frac{d^2}{N}=o(1)$, but this ceases to be true
as soon as $d^2\geq c N$ for any $c>0$. This reflects the fact that when $\frac{d^2}{N}=o(1)$ the
discreteness of the state space $\G_{N,d}$ is washed out in the limit (the limit is diffusive),
whereas when $d^2\geq c N$ the discrete nature of $\G_{N,d}$ is retained.

In  contrast, for larger values of $n$, i\.e\. for $n+2\leq d$, our upper bound on $F_2(n)$ is
uniform in $d$. Simplifying this bound further we show in Corollary \thv(A3.cor.3) that,
for $n+2\leq d$ and for large enough $d$, $F_2(n)$ is bounded above by a decreasing
function. This feature will be extremely useful in applications.

Finally, in Corollary \thv(A3.cor.3), we compare the functions $F_1(n)$ and $F_2(n)$.

\lemma{\TH(A3.lemma.1)}{\it With the notation of Definition \thv(P.def.1) we have:
$$
F_2(n)\leq \kappa^2(n+2) \frac{(n+2)!}{N^{(n+2)}}
\sum_{m\in I(n)}\frac{N^{(n+2-m)/2}}{[(n+2-m)/2]!}
\binom{d+m-1}{m}
\Eq(A3.1)
$$
In particular, 
\item{a)}for all $d$ and all large enough $N$,
$$
\eqalign{
\text{for} n=1\,,\quad
&F_2(n)
\leq
\kappa^2(3)
\left\{\sfrac{3!}{N}(\sfrac{d}{N})+(\sfrac{d+2}{N})^3\right)\}\cr
\text{for} n=2\,,\quad
&F_2(n)\leq
\kappa^2(4)\left\{\sfrac{4!}{2!N}(\sfrac{d+1}{N})^2+(\sfrac{d+3}{N})^4\right\}\cr
}
\Eq(A3.2)
$$

\item{b)} If $\sfrac{d^2}{N}=o(1)$, for all large enough $N$,
there exists a positive constant $C<\infty$ such that, setting
$$
p^*={
\cases \sfrac{n}{2}&\text{if $n$ is even}\cr
       \sfrac{n+1}{2}&\text{if $n$ is odd}\cr
\endcases}\,,\quad m^*=n+2-2p^*
\Eq(A3.13bis)
$$
\itemitem{$\bullet$} for all fixed $n$ independent of N satisfying $n^2\leq d-1$,
$$
F_2(n)
\leq C \frac{1}{N^{p^*}}
\left(\frac{d}{N}\right)^{m^*}
\Eq(A3.3-1)
$$
\itemitem{$\bullet$} for all $n+2\leq d$,
$$
F_2(n)
\leq C(n+2)^{\frac{3}{2}}\kappa^2(n+2)
\left(\bar\rho_{n,d}\right)^{\frac{n+2}{2}}
\left(\frac{n}{N}\right)^{\frac{n+2}{2}}
\Eq(A3.3-2)
$$
where
$$
\bar\rho_{n,d}=2\exp\left\{-1+\sfrac{n+2}{d-1}
+\sqrt{2\sfrac{d^2}{N}}\left(1+O\Bigl(\sqrt{\sfrac{d^2}{N}}\Bigr)\right)\right\}
\Eq(A3.3-3)
$$
\itemitem{$\bullet$} for all $n+2\geq d$,
$$
F_2(n)
\leq C\kappa^2(n+2)
\left(\rho_{n,d}\right)^{\frac{n+2}{2}}
\left(\frac{n}{N}\right)^{n+2-p^*}
\Eq(A3.3-4)
$$
where
$$
\rho_{n,d}=2e^{-1+2h(d/(n+2))} \text{and}h(x)=|x\log x|+x+x^2/2\,,\quad x\geq 0
\Eq(A3.3-5)
$$
\item{c)} If there exists a constant $c_0>0$ such that,  for all large enough $N$,
$\sfrac{d^2}{N}>c_0$, then
there exists a positive constant $C<\infty$ such that,
\itemitem{$\bullet$} For all $n+2\leq d$,
$$
\eqalign{
F_2(n)
&\leq
C(n+2)^{\frac{3}{2}}e^{\frac{(n+2)^2}{2(d-1)}}
\left(\frac{d}{N}\right)^{n+2}
}
\Eq(A3.4-1)
$$
\itemitem{$\bullet$} For all $n+2\geq d$,
$$
F_2(n)
\leq C\kappa^2(n+2)
\left(\rho_{n,d}\right)^{\frac{n+2}{2}}
\left(\frac{n}{N}\right)^{n+2-p^*}
\Eq(A3.4-2)
$$
where $\rho$ is defined in \eqv(A3.3-5).
}

Obviously our bounds on $F_2(n)$ are useful only if they guarantee that
$F_2(n)\leq 1$. Inspecting \eqv(A3.4-1) and  \eqv(A3.4-2) of assertion (c)
we see that this will always be the case when $d\leq d_0(N)$.


\remark The bound \eqv(A3.4-1) is the worst possible bound we could
derive from \eqv(A3.36). It is expected to be good only for small values
of $n$ (namely for fixed finite $n$ indepdendent of $N)$.
For larger values of $n$ one can improve it by
working directly with \eqv(A3.36).

Although we cannot prove that  $F_2(n)$, and hence $F(n)$, is a decreasing function,
the next corollary shows that this will be the case for suitably chosen $n$ and $d$.

\corollary{\TH(A3.cor.3)}{\it Let $d\geq\frac{\log N}{\log\log N}$. There exists $\varrho<1$
such that for all $n+2\geq d$ and large enough $N$,
$$
F(n)\leq \varrho^n\,
\Eq(A3.4-3)
$$
}

\proof Consider the right hand side of \eqv(A3.4-2). Given $\d<1$, let $C(\d)$ be defined by
$C(\d)=\arg\inf\left\{c>0\mid 2e^{-1+2h(d/(n+2))}\leq \d\right\}$. Next observe that
r\.h\.s\. of \eqv(A3.4-2) can be piecewise bounded above by decreasing functions as
follows: denoting by $C$ a finite positive constant whose value may change from line to line,
$$
\eqalign{
F_2(n)\leq
CN^2\left(2e^2C(\d)\frac{d}{N}\right)^{n/2}
&\text{if}
d\leq n+2\leq C(\d)d
\cr
F_2(n)\leq
CN^2\left(\d\frac{n}{N}\right)^{n/2}\,\,\,\,\,\,\,\,\,\,\,\,\,\,\,\,\,\,\,\,
&\text{if}
C(\d)d< n+2\leq \frac{N}{e}
\cr
F_2(n)\leq
CN^2\left(\frac{2}{e}(1+o(1))\right)^{n/2}
&\text{if}
n+2> \frac{N}{e}
}
\Eq(A3.4-3bis)
$$
By the second assertion of Corollary \thv(A3.cor.2), $F(n)\leq 2F_2(n)$.
Under the assumption that $d\geq\frac{\log N}{\log\log N}$,
the bound \eqv(A3.4-3) now easily follows.
\endproof

%
%


Of course the bound \eqv(A3.4-3) is a very coarse upper bound on \eqv(A3.3-4)
(or \eqv(A3.4-2)). Note that the larger $n$ is and the closer this bound gets to
\eqv(A3.3-4). The next Corollary contains a trite but useful upper bound
on $F_2(n)$ that will be good for very small values of $n$ only.

\corollary{\TH(A3.cor.2)}{\it  For all $d$ such that $\frac{d}{N}=o(1)$
we have, for all $N$ large enough:
\item{i)} $F_2(1)\geq F_2(n)$  for all $n\geq 1$,
\item{ii)}  $F_2(n)\geq F_1(n)$ for all $n\geq 3$.
}

\proof This is  an immediate consequence of
the bounds of Lemma \thv(A3.lemma.1).
\endproof

\proofof{Lemma \thv(A3.lemma.1)} Let $\left|\QQ^{\G}_d(n)\right|$ and
$\left|\QQ_d(n)\right|$ denote, respec.,
the number of solutions of \eqv(P.p2.31) and \eqv(P.p2.31bis). As established in Lemma
\thv(P.lemma.10), $\left|\del_{m}x\right|=\left|\QQ^{\G}_d(n)\right|$. But
$\left|\QQ^{\G}_d(n)\right|\leq\left|\QQ_d(n)\right|\leq\binom{d+m-1}{m}$, proving
\eqv(A3.1). The bounds of \eqv(A3.2) immediately follow from \eqv(A3.1). To further express
\eqv(A3.1) we will make use of the following lemma.

\lemma{\TH(A3.lemma.2)}{\it
$$
{\textstyle{\binom{d+m-1}{m}}}\leq\frac{m^{d-1}}{(d-1)!}e^{\frac{(d-1)^2}{2m}}\,,\quad\text{If} m\geq d
\Eq(A3.10)
$$
$$
{\textstyle{\binom{d+m-1}{m}}}\leq\frac{(d-1)^m}{m!}e^{\frac{m^2}{2(d-1)}}\,,\quad\text{If} m\leq d
\Eq(A3.11)
$$
$$
\sum_{m\in I(n)}{\textstyle{\binom{d+m-1}{m}}}\leq {\textstyle{\binom{d+(n+2)}{d}}}
\Eq(A3.12)
$$
}

\proofof{Lemma \thv(A3.lemma.2)} \eqv(A3.10) and \eqv(A3.11) are immediate and
\eqv(A3.12) follows from Pascal's recursion formula for the binomial coefficients (see eg [Co])
since
$$
\sum_{m\in I(n)}
{\textstyle{\binom{d+m-1}{m}}}
=\sum_{m\in I(n)}
{\textstyle{\binom{d+m-1}{d-1}}}
\leq \sum_{m=1}^{n+2}
{\textstyle{\binom{d+m-1}{d-1}}}
\leq
{\textstyle{\binom{d+(n+2)}{d}}}
$$
\endproof

Bearing in mind that $N$ and $d$ are fixed parameters, and $n$ the only variable,
let us now distinguish the cases $n+2< d$ et $n+2\geq d$.

\noindent $\bullet$ {\bf The case $n+2\geq d$.} For $m\in I(n)$ set
$p\equiv p(m)=(n+2-m)/2$. Note that the function $\N\ni p\mapsto{N^p}/{p!}$
is strictly increasing on $\{1,\dots,N\}$. Now, setting
$p^*\equiv p^*(n)=\max_{m\in I(n)}$, we have
$$
p^*={
\cases \sfrac{n}{2}&\text{if $n$ is even}\cr
       \sfrac{n+1}{2}&\text{if $n$ is odd}\cr
\endcases}
\Eq(A3.13)
$$
Thus $p\leq p^*<N$ for all $m\in I(n)$ and all $n$,  and
$$
\eqalign{
F_2(n)
&\leq
\kappa^2(n+2) \frac{N^{p^*}}{p^*!} \frac{(n+2)!}{N^{(n+2)}}
\sum_{m\in I(n)}\binom{d+m-1}{m}
\cr
&\leq
\kappa^2(n+2) \frac{N^{p^*}}{p^*!} \frac{(n+2)!}{N^{(n+2)}}
{\textstyle{\binom{d+(n+2)}{d}}}
\cr
&\leq
\kappa^2(n+2) \frac{N^{p^*}}{p^*!} \frac{(n+2)!}{N^{(n+2)}}
\frac{(n+2)^{d}}{d!}e^{\frac{d^2}{2(n+2)}}
\cr
}
\Eq(A3.14)
$$
where the last two lines follow, respectively, from \eqv(A3.12) and \eqv(A3.10).
Using Stirling's formula one then gets that, for some constant $0<C<\infty$,
$$
F_2(n)
\leq \kappa^2(n+2)
C\left(\frac{n}{N}\right)^{n+2-p^*}
\left(\rho_{n,d}\right)^{\frac{n+2}{2}}
\Eq(A3.15)
$$
where
$$
\rho_{n,d}=2e^{-1+2h(d/(n+2))}
\Eq(A3.17)
$$
and
$$
h(x)=|x\log x|+x+x^2/2\,,\quad x\geq 0
\Eq(A3.16)
$$
\remark Note that $h(x)$ is strictly decreasing on $[0,2]$, that $h(0)=0$, and $2h(1)=3$.
Thus for fixed $d$, $\rho_{n,d}$ is a decreasing function of $n$ that satisfies the bounds
$2e^{-1}\leq\rho_{n,d}\leq 2e^2$.
Moreover, one easily sees that there exists $\frac{1}{29}<\a<\frac{1}{30}$
such that, for $d<\a(n+2)$, $\rho_{n,d}< 1$.

Since \eqv(A3.15) is valid for all $d$ and all large enough $N$,
\eqv(A3.3-4) and \eqv(A3.4-2) are proven.

\noindent $\bullet$ {\bf The case $n+2\leq d$.} In this case,
since $m\leq n+2$, the bound \eqv(A3.11) applies for each
$m\in I(n)$ and thus,
$$
F_2(n)\leq \kappa^2(n+2) \frac{(n+2)!}{N^{(n+2)}}
\sum_{m\in I(n)}\frac{N^p}{p!}
\frac{(d-1)^m}{m!}e^{\frac{m^2}{2(d-1)}}
\Eq(A3.18)
$$
where, as before, $p\equiv p(m)=(n+2-m)/2$ for $m\in I(n)$. With $p^*$ as
in \eqv(A3.13), setting
$$
\rho=\frac{d-1}{\sqrt N}
\Eq(A3.19)
$$
\eqv(A3.18) may be rewritten as
$$
F_2(n)
\leq \kappa^2(n+2) e^{\frac{(n+2)^2}{2(d-1)}}
\frac{(n+2)!}{N^{(n+2)/2}}
\sum_{p=0}^{p^*}\rho^{n+2-2p}\frac{1}{p!(n+2-2p)!}
\Eq(A3.20)
$$
Defining
$$
f_{\rho}(x)\equiv [x\log x-x]+[(1-2x)\log(1-2x)-(1-2x)]+(1-x)\log(n+2)-(1-2x)\log \rho
\Eq(A3.21)
$$
we get, using Stirling's formula,
$$
\eqalign{
\sum_{p=0}^{p^*}\rho^{n+2-2p}\frac{1}{p!(n+2-2p)!}
\leq &\sum_{p=0}^{p^*}C\exp\left\{-(n+2)f_{\rho}(p/n)\right\}
\cr
\leq &\frac{n+3}{2}C\exp\left\{-(n+2)\inf_{0\leq x\leq 1/2} f_{\rho}(x)\right\}
}
\Eq(A3.22)
$$
where $0<C<\infty$. It is now easy to see that
$\inf_{0\leq x\leq 1/2} f_{\rho}(x)=f_{\rho}(x^*_{N,d}(n))$ where
$$
\eqalign{
x^*_{N,d}(n)&\equiv(\phi\circ\zeta_{N,d})(n)\cr
\zeta_{N,d}(n)&\equiv\frac{\rho^2}{4(n+2)}
\cr
\phi(z)&\equiv\frac{1}{2}\left\{(1+z)-\sqrt{(1+z)^2-1}\right\}\,,
\quad z\geq 0
\cr
}
\Eq(A3.23)
$$
To simplify the notation we will sometimes write
$x^*\equiv x^*(n)\equiv x^*_{N,d}(n)$. Using that $x^*_{N,d}(n)$
obeys the relation $f'_{\rho}(x^*_{N,d}(n))=0$, \eqv(A3.21) yields
$$
f_{\rho}(x^*)={
\cases
\log(1-2x^*)+\log(n+2)-(1-x^*)-\log\rho\,,
&\text{$0\leq x^*<\frac{1}{2}$}\cr
\frac{1}{2}\log(x^*)+\frac{1}{2}\log(n+2)-(1-x^*)\,,
&\text{$0< x^*\leq\frac{1}{2}$}\cr
\endcases}
\Eq(A3.24)
$$

Note now that the function $\phi(z)$ is strictly decreasing,
that $\phi(z)(0)=\sfrac{1}{2}$, $\lim_{z\rightarrow\infty}\phi(z)=0$,
and that
$$
\eqalign{
\phi(z)&=
\frac{1}{2}\left\{1+z-\sqrt{2z}\sqrt{1+z/2})\right\}=
\frac{1}{2}\left\{1+z-\sqrt{2z}(1+O(z))\right\}\,,
\quad z> 0
\cr
\phi(z)&=
\frac{1}{4(1+z)}\left\{1-O\left(\frac{1}{1+z}\right)^2\right\}\,,
\quad z\rightarrow\infty
\cr
}
\Eq(A3.25)
$$
Since $\zeta_{N,d}(n)$ is a strictly decreasing function of $n$,
$x^*_{N,d}(n)$ is itself strictly increasing and,
recalling that by assumption $2\leq n+2\leq d$,
$$
0\leq x^*_{N,d}(0)=\phi(\sfrac{d^2}{8N})
\leq x^*_{N,d}(n)\leq x^*_{N,d}(d-2)=\phi(\sfrac{d}{4N})
\leq \sfrac{1}{2}
\Eq(A3.26)
$$
On the other hand, by \eqv(A3.24), $f_{\rho}(x^*)$ is a strictly
increasing function of $x^*$, and thus
$$
f_{\rho}(0)\leq f_{\rho}(x^*_{N,d}(0))=f_{\rho}(\phi(\sfrac{d^2}{8N}))\leq
f_{\rho}(x^*_{N,d}(n))\leq
f_{\rho}(x^*_{N,d}(d-2))=f_{\rho}(\phi(\sfrac{d}{4N}))
\leq f_{\rho}(\sfrac{1}{2})
\Eq(A3.27)
$$
Using that $\frac{d}{N}\leq 1$ one easily checks that, by the first line of
\eqv(A3.25) and the second line of \eqv(A3.24), using the series expansion of
$\log(1+u)$, $|u|<1$,
$$
f_{\rho}(x^*_{N,d}(d-2))=f_{\rho}(\sfrac{1}{2})
-\sqrt{\sfrac{d}{2N}}\left(1+O\Bigl(\sqrt{\sfrac{d}{N}}\Bigr)\right),\,
\quad \sfrac{d}{N}< 1
\Eq(A3.28)
$$
From this and \eqv(A3.27) we see that
the range of $f_{\rho}(x^*_{N,d}(n))$
will depend on the behavior of $\frac{d^2}{N}$
and $\frac{d}{N}$.

If $\frac{d^2}{N}=o(1)$ then $\frac{d}{N}=o(1)$. Then just as before we have
$$
f_{\rho}(x^*_{N,d}(0))=f_{\rho}(\sfrac{1}{2})
-\sqrt{2\sfrac{d^2}{N}}\left(1+O(\sqrt{d^2/N})\right),\,
\quad \sfrac{d^2}{N}\rightarrow 0
\Eq(A3.29)
$$
and \eqv(A3.27),  \eqv(A3.28), and  \eqv(A3.29) imply that,
for all $0\leq n\leq d-2$,
$$
\left|f_{\rho}(x^*_{N,d}(n))-f_{\rho}(\sfrac{1}{2})\right|
\leq
\sqrt{2\sfrac{d^2}{N}}\left(1+O(\sqrt{d^2/N})\right)\,,
\quad \sfrac{d^2}{N}\rightarrow 0
\Eq(A3.30)
$$
In other words, $f_{\rho}(x^*_{N,d}(n))$ remains essentially constant
for $n\in\{0,\dots,d-2\}$.
If on the contrary there exist positive finite constants $c_0, N_0$
such that $\sfrac{d^2}{N}>c_0$ for all $N>N_0$,
then $f_{\rho}(x^*_{N,d}(n))$ is no longer constant
when $n$ varies from $0$ to $d-2$. In particular, if
$\sfrac{d^2}{N}\rightarrow \infty$
then, by the second line of \eqv(A3.25) and the first line of \eqv(A3.24),
$$
f_{\rho}(x^*_{N,d}(0))=f_{\rho}(0)
-\frac{1}{4(1+d^2/N)}\left(1+O(\sqrt{N/d^2})\right),\,
\quad \sfrac{d^2}{N}\rightarrow \infty
\Eq(A3.31)
$$
so that $f_{\rho}(x^*_{N,d}(n))$ ranges from $f_{\rho}(0)$
to $f_{\rho}(x^*_{N,d}(d-2))=f_{\rho}(\phi(\sfrac{d}{4N}))$ when $n$ varies from $0$ to $d-2$.
But $\frac{d}{N}$ may only vary from $0$ to $1$ and thus,
$f_{\rho}(0)< f_{\rho}(\phi(\sfrac{1}{4})) \leq f_{\rho}(\phi(\sfrac{d}{4N}))\leq f_{\rho}(\sfrac{1}{2})$.
Now
$$
\frac{(n+2)!}{N^{(n+2)/2}}\exp\{-(n+2)f_{\rho}(\sfrac{1}{2})\}\leq
c\sqrt{n+2}\left(\frac{2}{e}\frac{n}{N}\right)^{\frac{n+2}{2}}
\Eq(A3.32)
$$
and
$$
\frac{(n+2)!}{N^{(n+2)/2}}\exp\{-(n+2)f_{\rho}(0)\}\leq
c\sqrt{n+2}\left(\frac{d}{N}\right)^{n+2}
\Eq(A3.33)
$$
Therefore, if $\sfrac{d^2}{N}=o(1)$, collecting \eqv(A3.20), \eqv(A3.22),
\eqv(A3.30), and \eqv(A3.32), there exists a positive constant $C<\infty$
such that,
$$
F_2(n)
\leq C(n+2)^{\frac{3}{2}}\kappa^2(n+2)
\left(\frac{n}{N}\right)^{\frac{n+2}{2}}
\left(\bar\rho_{n,d}\right)^{\frac{n+2}{2}}\,,
\quad \sfrac{d^2}{N}\rightarrow 0
\Eq(A3.34)
$$
where
$$
\bar\rho_{n,d}=2\exp\left\{-1+\sfrac{n+2}{d-1}
+\sqrt{2\sfrac{d^2}{N}}\left(1+O\Bigl(\sqrt{{d^2}/{N}}\Bigr)\right)\right\}
\Eq(A3.35)
$$
Otherwise, if there exist positive finite constants $c_0, N_0$
such that $\sfrac{d^2}{N}>c_0$ for all $N>N_0$, then by \eqv(A3.20),
\eqv(A3.22), \eqv(A3.27), and \eqv(A3.33),
$$
\eqalign{
F_2(n)
&\leq  C(n+2)\kappa^2(n+2) e^{\frac{(n+2)^2}{2(d-1)}}
\frac{(n+2)!}{N^{(n+2)/2}}
\exp\{-(n+2)f_{\rho}(x^*_{N,d}(n))\}
}
\Eq(A3.36)
$$
In particular, using that $f_{\rho}(x^*_{N,d}(n))\geq f_{\rho}(0)$
(see \eqv(A3.27)) together with \eqv(A3.33), it follows from
\eqv(A3.36) that for all $0\leq n\leq d-2$
$$
\eqalign{
F_2(n)
&\leq
C(n+2)^{\frac{3}{2}}e^{\frac{(n+2)^2}{2(d-1)}}
\left(\frac{d}{N}\right)^{n+2}
}
\Eq(A3.37)
$$
This bound, valid for all $n+2\leq d$, is only reasonable however
for small $n$ (more precisely for fixed finite $n$ independent
of $N$) when $ \sfrac{d^2}{N}\rightarrow \infty$.

Since \eqv(A3.35) proves \eqv(A3.3-2) and \eqv(A3.37) proves
\eqv(A3.4-1) it remains to prove \eqv(A3.3-1). We will treat the case
 $n$ odd only; the case of even $n$ is similar. Here $m, p$, and $n$ all
are fixed and independent of $N$, and
since $\frac{m^2}{2(d-1)}\leq \frac{n^2}{2(d-1)}\leq 1$, it follows from
\eqv(A3.11) that,
$$
{\textstyle{\binom{d+m-1}{m}}}\leq c(d-1)^m
\Eq(A3.38)
$$
for some constant $0<c<\infty$. Thus,
$$
F_2(n) \leq
c\kappa^2(n+2)\sum_{m\in I(n)} \frac{N^{p}}{N^{(n+2)}} (d-1)^m
=c\kappa^2(n+2)\sum_{m\in I(n)} \frac{1}{N^{p}} \left(\frac{d-1}{N}\right)^m
\Eq(A3.39)
$$
We want to show that the leading term in the sum above is given by
$(m, p)=(m^*, p^*)$. But one easily verifies that for all $m>m^*$,
$$
\frac{1}{N^{p}} \left(\frac{d-1}{N}\right)^m N^{p^*} \left(\frac{N}{d-1}\right)^{m^*}=
\left(\frac{(d-1)^2}{N}\right)^{m-1}.
\Eq(A3.41)
$$
Hence, for some constants $c',c'',c'''>0$,
$$
\eqalign{
F_2(n) &\leq
c'
\frac{1}{N^{p^*}} \left(\frac{d-1}{N}\right)^{m^*}
\sum_{m\in I(n)}\left(\frac{(d-1)^2}{N}\right)^{m-1}
\cr
&\leq
c''
\frac{1}{N^{p^*}} \left(\frac{d-1}{N}\right)^{m^*}
\sum_{m\geq 1}\left(\frac{d^2}{N}\right)^{m-1}
\cr
&\leq
c'''
\frac{1}{N^{p^*}} \left(\frac{d-1}{N}\right)^{m^*}
\cr
}
\Eq(A3.42)
$$
which proves  \eqv(A3.3-1) for odd values of $n$.

The proof of Lemma \thv(A3.lemma.1) is now complete.
\endproof

\bigskip
\vfill\eject


\bigskip
\chap{11. Appendix A4}{11}

Let $A\subset\SS_N$ be compatible with some $d$-partition $\L$.
In this appendix we collect a few ad hoc estimates on $\UU_{N,d}(A)$ (see \thv(I.12))
that allow to quantify the sparseness of the set $A$
in two cases: roughly speaking 1) when $|A|$ is small enough and
2) when the elements  of $A$ satisfy a certain minimal distance assumption.
This estimates are derived from elementary observations stemming from
Definition \thv(I.def.1)
and the properties of the function $F_{N,d}$.


\noindent{\bf Case 1).}
Our first three results (Lemma \thv(A4.lemma.1), Lemma \thv(I.lemma.3),
and Corollary \thv(I.cor.1)) are concerned with `small' subsets $A$ of $\SS_N$.
In Lemma \thv(A4.lemma.1) we provide a sufficient condition on the size of $A$
which entails that $A$ is compatible with some
$d$-partition $\L$.

\lemma{\TH(A4.lemma.1)} {\it Let $A\subset\SS_N$ be such that
$2^{|A|}\leq N$. Then there exists a  $d$-partition $\L$ with
$d\leq 2^{|A|}$ such that, for any $\xi\in\SS_N$, $A$ is
$(\L,\xi)$-compatible. If $|A|=1$ one may chose
the trivial partition $\L=\{1,\dots,N\}$. In this case $d=1$.

}


The next lemma allows to quantify the sparseness of sets $A$ of arbitrary size
but is clearly useful for small enough sets only .

\lemma{\TH(I.lemma.3)} {\it Let $A\subset\SS_N$. For all $d$ such that
$\frac{d}{N}=o(1)$ and all $N$ large enough,
$$
\UU_{N,d}(A)\leq C|A|\max\left\{\frac{1}{N},\Bigl(\frac{d}{N}\Bigr)^3\right\}
\Eq(A4.1)
$$
$$
|A|F(n)\leq C|A|\max\left\{\frac{1}{N},\Bigl(\frac{d}{N}\Bigr)^3\right\}\text{for all} n\geq 1
\Eq(A4.1')
$$
for some constant $0<C<\infty$. In particular, if $d\leq\a\frac{N}{\log N}$
for some constant $\a>0$,
$$
\UU_{N,d}(A)\leq C|A|\Bigl(\frac{\a}{\log N}\Bigr)^3
\Eq(A4.2)
$$
$$
|A|F(n)\leq C|A|\Bigl(\frac{\a}{\log N}\Bigr)^3\text{for all} n\geq 1
\Eq(A4.2')
$$
}

Finally, the corollary below is geared to the case $d\leq d_0(N)$ for which
most results in this paper obtain.
Combining Lemma \thv(A4.lemma.1) and \eqv(A4.2) of Lemma \thv(I.lemma.3),
we can conclude that:

\corollary{\TH(I.cor.1)} {\it Let $A\subset\SS_N$ be such that
$2^{|A|}\leq C\frac{N}{\log N}$ for some $0<C<\infty$.
Then there exists a  $d$-partition $\L$ with
$d\leq C\frac{N}{\log N}$ such that, for any $\xi\in\SS_N$, $A$ is
$(\L,\xi)$-compatible and
$$
\UU_{N,d}(A)\leq C'\frac{1}{(\log N)^2}
\Eq(A4.3)
$$
$$
|A|F(n)\leq C'\frac{1}{(\log N)^2}\text{for all}n\geq 1
\Eq(A4.3')
$$
for some constant $0<C'<\infty$.
}


\proofof{Lemma \thv(A4.lemma.1)}
the case  $|A|=1$ is immediate. Let us assume that  $|A|\geq 2$ and
call $\s^1,\dots,\s^{|A|}$ the  elements of $A$, i\.e\. set
$A=\{\s^1,\dots,\s^{|A|}\}$.
We define a partition of the set $\{1,\dots,i,\dots,N\}$ into
$d:=2^{|A|}$ subsets $\L_k$, $1\leq k\leq d$ in the following way.
Let us identify the collection $A$ to the $|A|\times N$
matrix whose row vectors are the configurations $\s^{\mu}$,
$$
\s^{\mu}=(\s^{\mu}_{i})_{i=1,\dots,N}\in\SS_N,\,\,\,\,\mu\in\{1,\dots,|A|\}\,,
\Eq(II.32)
$$
and denote by $\s_{i}$ the column vectors
$$
\s_{i}=(\s^{\mu}_{i})^{\mu=1,\dots,|A|}\in\SS_{|A|},\,\,\,\,i\in\{1,\dots,N\}
\Eq(II.33)
$$
(hence $\s^{\mu}_{i}$ is the element lying at the intersection of the $\mu$-th
row and $i$-th colum).
Observe that, when carrying an index placed as a superscript, the letter $\s$
refers to an element of the cube $\SS_N$ while, when carrying an index placed as a
subscript, it refers to an element of the cube $\SS_{|A|}$.
Next, let $\{e_1,\dots,e_k,\dots,e_{d}\}$ be an arbitrarily chosen
labelling of
all $d=2^{|A|}$ elements of $\SS_{|A|}$. Then, since $d=2^{|A|}\leq N$,
$A$ induces a partition $\L$ of $\{1,\dots,N\}$ into
at most $d$
classes  $\L_k$, $1\leq k\leq d$, defined by
$$
\L_k=\{1\leq i\leq N\mid \s_i=e_k\}\,
\Eq(2.1.6)
$$
if and only if $\L_k\neq\emptyset$.
Now clearly, for any $\xi\in\SS_N$, $A$ is $(\L,\xi)$-compatible.
\endproof

\proofof{Lemma \thv(I.lemma.3)} By Definition \thv(P.def.1)
and Corollary \thv(A3.lemma.2),
$$
F(n)=F_1(n)+F_2(n)\leq F_1(1)+F_2(1)\leq  C\max\left\{\frac{1}{N},\Bigl(\frac{d}{N}\Bigr)^3\right\}
\Eq(A4.4)
$$
where the last inequality, valid for some constant $0<C<\infty$,
follows from the bound \eqv(A3.2) on $F_2(1)$
and the fact that, by definition (see \eqv(P.p2.03)),
$F_1(1)=\kappa(1)\frac{1}{N}$. This immediately yields  \eqv(A4.1').
Moreover, inserting this bound in the definition \eqv(I.12) of $\UU_{N,d}(A)$
yields \eqv(A4.1). From this \eqv(A4.2) and  \eqv(A4.2') are immediate. \endproof

\proofof{Corollary \thv(I.cor.1)} Since $2^{|A|}\leq C\frac{N}{\log N}\leq N$,
it follows from  Lemma \thv(A4.lemma.1) that there exists a $d$-partition
$\L$ with
$d\leq 2^{|A|}\leq C\frac{N}{\log N}$ such that, for any $\xi\in\SS_N$, $A$ is
$(\L,\xi)$-compatible. For such a $d$, the bounds
\eqv(A4.2) and \eqv(A4.2') apply, proving \eqv(A4.3) and \eqv(A4.3').
\endproof

\noindent{\bf Case 2).} In what follows we consider sets $A$ that are compatible with some
$d$-partition $\L$, for $d\geq\frac{\log N}{\log\log N}$.

\lemma{\TH(A4.lemma.2)} {\it Let $d\geq\frac{\log N}{\log\log N}$ and let $A$ be
compatible with some $d$-partition $\L$.  There exists $\varrho<1$
such that for all $C\geq 1$ and all $n\geq C d$ we have, for large enough $N$,
$$
|A|F(n)\leq \varrho^{n(1-\e)} \text{where} \e={\log 2}/{C}
\Eq(A4.5)
$$
}

It is now easy to deduce a bound on $\UU_{N,d}(A)$
when the minimal distance between
the elements of $A$ is larger than $d$.

\lemma{\TH(A4.lemma.3)} {\it  Let $d\geq\frac{\log N}{\log\log N}$ and let $A$ be
compatible with some $d$-partition $\L$. Set
$$
n^*:=\inf_{\eta\in A}\Dist(\eta,A\setminus\eta)
\Eq(A4.6)
$$
There exists $\varrho<1$ such that if $n^*\geq C d$ for some $C\geq 1$ then,
for large enough $N$,
$$
\UU_{N,d}(A)\leq \varrho^{n^*(1-\e)} \text{where} \e={\log 2}/{C}
\Eq(A4.7)
$$
}

\proofof{Lemma \thv(A4.lemma.2)} Since by assumption $A$ is compatible with
some $d$-partition $\L$ then $|A|\leq 2^d$. This observation combined with
Corollary \thv(A3.cor.3) of Appendix A3 proves the lemma.\endproof

\proofof{Lemma \thv(A4.lemma.3)} Assume that  $n^*\geq C d$ for some $C\geq 1$.
By Corollary \thv(A3.cor.3) of Appendix A3,  $F(n)$ restricted to the set
$n\geq n^*$ is decreasing. By \eqv(I.12) we then have,
$$
\eqalign{
\UU_{N,d}(A)
&=
\max_{\eta\in A}\sum_{\s\in A\setminus\eta}F(\Dist(\eta,\s))
\cr
&\leq
\max_{\eta\in A}\sum_{\s\in A\setminus\eta}F(n^*)
\cr
&\leq
|A|F_{N,d}(n^*)
\cr
&\leq
\varrho^{n^*(1-\e)} 
\cr
}
\Eq(A4.8)
$$
where the last line follows from Lemma \thv(A4.lemma.2). \endproof

\bigskip
\vfill\eject

\chap{References.} 0

\frenchspacing

\item{[A1]} D.~Aldous, ``Random walks on finite groups and rapidly mixing Markov chains'',
Seminar on probability, XVII,  243--297, Lecture Notes in Math., 986, Springer, Berlin, 1983.
\item{[A2]} D.~Aldous, ``On the time taken by random walks on finite groups to visit every state'',
Z. Wahrsch. Verw. Gebiete {\bf 62}, 361--374 (1983).
\item{[AD1]} D.~Aldous and P.~Diaconis ``Strong uniform times and finite random walks'',
 Adv. in Appl. Math. {\bf 8},  69--97 (1987).
\item{[AD2]} D.~Aldous and P.~Diaconis ``Shuffling cards and stopping times'',
  Amer. Math. Monthly   {\bf 93}, 333--348  (1986).
\item{[BBG1]} G. Ben Arous, A.~Bovier, and V. Gayrard, ``Glauber dynamics
             of the random energy model. 1. Metastable motion on the extreme
             states.'',  Commun. Math. Phys.  {\bf 235}, 379-425 (2003).
\item{[BBG2]} G. Ben Arous, A.~Bovier, and V. Gayrard, ``Glauber dynamics
             of the random energy model. 2. Aging below the critical
             temperature'', Commun. Math. Phys. {\bf 236}, 1-54 (2003).
\item{[BC]}  G. Ben Arous and J. Cerny, ``The Arcsine law as a universal aging scheme for trap models'', preprint (2006).

\item{[BEGK1]} A\.~Bovier, M\.~Eckhoff, V\.~Gayrard and M\.~Klein,
``Metastability in stochastic dynamics of disordered mean field models'',
{\it Prob\. Theor\. Rel\. Fields} {\bf 119}, 99-161 (2001).
\item{[BEGK2]} A\. Bovier, M\. Eckhoff, V\. Gayrard and M\. Klein,
``Metastability and low lying spectra in reversible Markov chains'',
to appear in {\it Commun\. Math\. Phys\.} (2001).
\item{[BG]} A\. Bovier and  V\. Gayrard, ``An almost sure large deviation principle for the Hopfield model'',
Ann. Probab., {\bf  24},  1444--1475 (1996).
\item{[G1]} V\. Gayrard,
``Thermodynamic limit of the $q$-state Potts-Hopfield
model with infinitely many patterns'', {\it J\. Stat\. Phys\.} {\bf
68}, 977-1011 (1992).
\item{[G2]} V\. Gayrard,  ``Glauber dynamics of the 2-GREM. 1. Metastable motion on the extreme
             states.'', preprint (2006).
\item{[BR]} C.J. Burke and M. Rosenblatt, ``A Markovian function of a
Markov Chain'',
Ann. Math. Statist. {\bf  29}, 1112-1122, (1958).
\item{[Co]} L. Comtet, {\it Analyse combinatoire. Tome I, II},
Presses Universitaires de France, Paris, 1970.
\item{[D]} P.~Diaconis, ``Applications of noncommutative Fourier analysis to
 probability problems'', \'Ecole d'\'Et\'e de Probabilit\'es de
   Saint-Flour XV--XVII, 1985--87, 51--100,
Lecture Notes in Math., 1362, Springer, Berlin, 1988.
\item{[DGM]}
P.~Diaconis, R.L.~Graham, J.A.~Morrison,
Asymptotic analysis of a random walk on a hypercube with many dimensions.
Random Structures Algorithms {\bf 1}, 51--72  (1990).
\item{[KS]} J.G.~Kemeny and J.L.~Snell, Finite Markov chains,
D. van Nostrand Company, Princeton, 1960.
\item{[K]} J.H.B~Kemperman, The passage problem for a stationary Markov chain,
The University of Chicago Press, Chicago, 1961.

\item{[Fe]} W. Feller, {\it An Introduction to Probability Theory and Its Applications, Volume 2}, John Wiley and Sons, Inc. New York, 1970.

\item{[KP]} H. Koch and J. Piasko, ``Some rigorous results on the Hopfield Neural Network model''
{\it J\. Stat\. Phys\.} {\bf55}, 903 (1989).

\item{[Li]} T\.M\. Liggett, {\it Interacting particle systems}, Springer, Berlin, 1985.
\item{[M1]} P.~Matthews,
``Some sample path properties of a random walk on the cube'',
J. Theoret. Probab. 2 (1989), no. 1, 129--146.
\item{[M2]} P.~Matthews, ``Covering problems for Markov chains'' Ann. Probab. {\bf 16}, 1215--1228 (1988).
\item{[M3]} P.~Matthews, ``Mixing rates for a random walk on the cube'' SIAM J. Algebraic Discrete Methods {\bf 8}, 746--752  (1987).
\item{[SaCo]} L.~Saloff-Coste, {\it Lectures on finite Markov
   chains}, Lectures on probability theory and statistics (Saint-Flour,
   1996), 301--413, Lecture Notes in Math., 1665, Springer, Berlin,
   1997.
\item{[So]}  P.M. Soardi, {\it Potential theory on infinite networks},
LNM 1590, Springer, Berlin-Heidelberg-New York, 1994.
\item{[Sp]} F. Spitzer, {\it Principles of random walks}, Second edition.
Graduate Texts in Mathematics, Vol. 34. Springer-Verlag,
New York-Heidelberg, 1976.

\item{[V]} M\. Voit, ``Asymptotic distributions for the Ehrenfest Urn and related walks"
Journal of Applied Probability {\bf 33}, 340-356 (1996).

\end

\end